\documentclass[a4paper]{amsbook}   

\usepackage{amssymb} \usepackage{amscd} \usepackage{times}

\def\zbb{\mathbb{Z}}  
  
  \def\phi{\varphi}
 \def\p1{{\mathbb{P}^1_\zbb}}

\newtheorem{Theorem}{\quad Th\'eor\`eme}[section] 

\newtheorem{Definition}[Theorem]{\quad Definition} 

\newtheorem{Corollary}[Theorem]{\quad Corollaire}

\usepackage[latin1]{inputenc}
\usepackage[francais]{babel}

\begin{document}

%\tableofcontents

\title{ {\vspace*{10cm}  Sur les estimations a priori du type sup+inf et sup $\times$ inf.}}

\author{ Samy Skander Bahoura }

%\address{Departement de Mathematiques, Universite Pierre et Marie Curie, 2 place Jussieu, 75005, Paris, France.}
              
%\email{samybahoura@yahoo.fr} 

\date{}

\maketitle
\newpage

\vspace*{3cm}

\tableofcontents{}

\newpage

\section{Introduction} 

On part de l'\'equation de Yamabe et celle plus g\'en\'erale de la Courbure Scalaire Prescrite en dimension $ \geq 2 $. Nous pr\'esentons quelques in\'egalit\'es caract\'eristiques de ces \'equations.

\bigskip

Dans la suite, nous notons $ \Delta =\nabla^i(\nabla_i) $.

\bigskip

\underbar { \bf { 1. Cas des fonctions harmoniques:}}

\bigskip

L'exemple fondamental est celui des c\'el\`ebres in\'egalit\'es de Harnack pour les fonctions harmoniques positives ou nulles:

\bigskip

Consid\'erons $ \Omega $ un ouvert de $ {\mathbb R}^n, n\geq 2 $ et soit $ u \geq 0 $ telle que:

$$ -\Delta u = -\sum_{k=1}^n \partial_{kk} u = 0. $$

Alors, pour tout compact $ K $ de $ \Omega $, il existe une constante positive $ c = c(K, \Omega, n) $ telle que:

$$ \dfrac{ \sup_K u}{\inf_K u} \leq c. $$

Nous pouvons avoir le m\^eme type d'in\'egalit\'es pour des op\'erateurs elliptiques d'ordre 2 plus g\'en\'eraux.

\bigskip

\underbar { \bf { 2. Cas des fonctions sous-harmoniques:}}

\bigskip

On suppose que $ B_{2R}(0) \subset \Omega $ et que:

$$ -\Delta u \geq 0. $$

Alors, il existe  une constante positive $ c $ telle que:

$$ \left (\int_{B_R(0)} u^p \right )^{1/p} \leq c \min_{B_R(0)} u. $$

\underbar { \bf { 3. Cas des fonctions sur-harmoniques:}}

\bigskip

On suppose toujours $ B_{2R}(0) \subset \Omega $ et:

$$ -\Delta u \leq 0. $$

Alors, il existe une constante positive $ c $ telle que:

$$ \max_{B_R(0)} u \leq c \left ( \int_{B_R(0)} u^p \right )^{1/p}. $$

Dans ce qui suit, on s'int\'eresse aux estimations a priori pour les solutions d'\'equations provenant de la g\'eom\'etrie conforme. Les \'equations types sont de la forme:

$$ -\Delta u+ R_g=V(x)e^u, \,\,\, {\rm en \,\, dimension \,\, 2 } \,\,\, \qquad (E_1) $$ 

$$ -\Delta_g u+R_gu=Vu^{N-1},\,\, u >0, \,\, {\rm et } \,\, N=\dfrac{n+2}{n-2}. ,\,\, {\rm en \,\, dimension \,\, n \geq 3 } \,\,\,\qquad (E_2) $$

Ici, $ R_g $ est la courbure scalaire de la vari\'et\'e Riemannienne $ (M, g) $ (surface en dimenion 2). On peut remplacer la fonction $ R_g $ par une autre fonction $ h $ et ces equations deviennent plus g\'en\'erales. Dans le cas $ h=R_g $, ces \'equations expriment le changement de m\'etrique conforme et $ V $ est la courbure scalaire prescrite (la courbure scalaire pour la nouvelle metrique $ \tilde g= e^u g $).

\bigskip

Les estimations a priori recherch\'ees, sont de la forme:

$$ \sup_K u +C_1\inf_M u <C_2,  \,\,\, {\rm ou } \,\,\, \sup_M u + C_1 \inf_K u > C_2 \,\,\, {\rm en \,\, dimension \,\, 2 } \,\,$$
et 

$$ \sup_K u  \times \inf_M u <C_2,  \,\,\, {\rm ou } \,\,\, \sup_M u \times \inf_K u > C_2 \,\,\, {\rm en \,\, dimension \,\,  n \geq 3 } \,\,$$

Une question interm\'ediaire est de savoir s'il est possible d'avoir l'estimation a priori suivante, en supposant une minoration du minimum,

$$ \inf_M u \geq m >0 \,\Rightarrow^?  \,\,\, \sup_K u \leq C(V, m, K, M, g) $$

Comme ce que Brezis et Merle on fait en dimension 2. La constante $ C(V, m, K, M, g) $ depend des bornes a priori de $ V $, par exemple, en supposant:

$$  0 < a \leq V(x) \leq b < + \infty $$

ou, si de plus on rajoute la condition sur le gradient (par exemple) ou les d\'eriv\'ees successives de $ V $,

$$ ||\nabla V||_{\infty}\leq A $$

Par exemple:

$$ C(V, m, K, M, g)=C(a,b, A, m, K, M, g). $$

\bigskip
 
Notons que des estimations de ce type ont ete prouv\'ees par Siu et Tian pour l'\'equation de Monge-Ampere complexe avec condition sur la classe de Chern (positive).

\bigskip

\begin{center}

$ \begin{cases} 

(\omega_g+\partial \bar \partial \phi)^n=e^{f-t\phi} \omega_g^n, \\

\omega_g+\partial \bar \partial \phi >0 \,\, {\rm on } \,\,  M  \\
 
\end{cases}  $ 
    
\end{center} 

\bigskip

Ils prouvent des estimations du type:

 $$ \sup_M(\phi-\psi)+m \inf_M(\phi-\psi) \leq C(t) \,\, {\rm ou } \,\, \sup_M(\phi-\psi)+m \inf_M(\phi-\psi) \geq C(t), $$
 
quand la premiere classe de Chern est positive.

\bigskip

La fonction  $ \psi $  est  de classe $ C^2 $ et telle que:

$$ \omega_g+\partial \bar \partial \psi \geq 0 \,\, {\rm and } \,\,  \int_M e^{f-t\psi}\omega_g^n = Vol_g(M), $$

Voir [1-56]

\newpage

On commence par exposer quelques r\'esultats concernant les in\'egalit\'es du type $ \sup + \inf  $ en dimension 2. Pour obtenir son r\'esultat, Shafrir  utilise des inegalit\'es isoperim\'etriques g\'en\'erales et dont la preuve est bas\'ee sur des in\'egalit\'es sur les longueur des paralelles, on considere les coordonn\'ees g\'eodesiques paralelles. Ces coordonn\'ees localisent un point par son angle par rapport a une courbe donn\'ee et sa distance a cette m\^eme courbe. Il faut d\'eriver la fonction "longueur d'une paralelle" et utiliser la formule de Gauss-Bonnet. L'autre m\'ethode est celle de C. Bandle se basant sur les lignes de niveaux, le but alors est d'obtenir une inegalit\'e g\'eom\`etrique liant la longueur (longueur riemannienne) d'une courbre ferm\'ee et l'aire (volume riemannien) de la surface qu'elle d\'elimite. Cette m\'ethode utilise l'in\'egalit\'e isoperim\'etrique de Nehari (qui concerne les fonctions harmoniques) et une sorte de symm\'etrisation de Schwarz. (cette preuve est bas\'ee sur la preuve de l'in\'egalit\'e isoperimetrique classique, par les series de Fourier, dans le cas d'une composante connexe, puis pour le cas d'un nombre fini de composante connexes).

ON s'intersse aussi aux r\'esultats de CC.Chen et CS.Lin, ils utilisent les deux principaux faits suivants:  un r\'esultat de Suzuki et une in\'egalit\'e g\'eom\`etrique concernant les courbures int\'egrales pour am\'eliorer le r\'esultat de Shafrir et obtenir une in\'egalit\'e du type $ \sup + \inf $ optimale, lorsque les fonctions potentiel (courbure scalaire prescrite) sont uniform\'ement holderiennes. Le r\'esultat de Suzuki am\'eliore celui de C. Bandle, celle-ci prouve une sorte de principe du maximum pour une classe d'\'equations elliptiques non-lin\'eaire avec une non-lin\'earite exponentielle. Plus exactement, on a une information sur le maximum des solutions de cette \'equation en fonction de ces valeurs au bord  en supposant le volume strictement plus petit que $ 8 \pi $. Le r\'esultat de Suzuki am\'eliore ce dernier en remplacant les valeurs aux bord par une int\'egrale de la fonction au bord, c'est une in\'egalit\'e de Harnack. Quant a l'in\'egalit\'e g\'eom\`etrique elle est obtenue en int\'egrant une in\'equation diff\'erentielle.

Le travail de CC.Chen et C.S. Lin a \'et\'e fait pour les fonctions blow-ups. Ils raisonnent par l'absurde et prouvent des estimations asymptotiques pour les fonctions blow-up, grace a l'in\'egalit\'e de Suzuki puis  aboutissent a une contradiction en utilisant l'in\'egalit\'e g\'eom\'etrique.

Dans la continuit\'e du travail de ces deux auteurs (CC.Chen et C.S. Lin) et c'est l'objet d'un article publi\'e(cas ou les courbures prescrites sont holderiennes), nous obtenons une nouvelle d\'emonstartion bas\'ee sur la m\'ethode moving-plane.

La m\'ethode moving-plane est bas\'ee sur le principe du maximum et lemme de Hopf. Le fait d'avoir consid\'erer les courbures prescrites non constantes, empeche d'utiliser le principe du maximum directement, car il y a un terme qui perturbe l'\'equation, on corrige alors ce terme on consid\'erant une fonction auxili\`ere. On applique alors le principe du maximum \`a la somme des deux fonctions (fonction principale+fonction auxili\`ere). (Pour cela, il faut verifer plusieurs hypotheses, positivt\'e de la fonctions, conditions aux bords, et s'il existe un rang, qui permettent de commencer la "methode moving plane", ceci est valable pour les resultats de CC.Chen et C-S; Lin, en dimensions $ \geq 3 $).

Pour ce qui est des r\'esultats de Suzuki et Brezis-Merle, c'est  la quantisation, c'est le comportement des points blow-ups et des fonctions autour de ces blow-ups.  Ces derniers, en supposant que pour une suite de fonctions solution d'une EDP elliptique diverge, alors,  il y a une sous-suite qui v\'erifie  une alternative:  soit,  elles convergent sur tout compact, soit, il y a un nombre fini de points-blow ups limite, elles se concentrent, soit, elle divergent completement sur tout compact. La particularit\'e du r\'esultat de Suzuki est qu'il est obtenu grace aux r\'esultats sur les zeros des fonctions holomorphes et que les masses (le co\'efficient devant la masse de Dirac), en cas de concentrations sont egales \`a $ 8 \pi $. Quant a ceux de Brezis-Merle, ils sont la consequence de la convergence faible dans $ L^1 $ d'une suite de fonctions vers une mesure qui peut avoir des points de concentrations.

Notons qu'une des consequences du travail de quantisation est l'obtension d'une estimation du supremum des solutions en fonctions de l'infimum, d'ou la question de savoir comment varie ce maximum en fonction du minimum et c'est le resultat obtenu par Shafrir, puis am\'elior\'e par CC.Chen et CS.Lin. Et dont on donne une nouvelle preuve.

Dans notre travail, nous \'etendons ce type de ph\'enomene au bord. De plus nous formulons une preuve d'un r\'esultat de compacit\'e par un argument de convergence faible au sens $ L^2 $ (sur le bord). Ce resultat est valable pour des domaines lisses ou $ C^{2,\alpha},\alpha >0 $. 

Pour obtenir une minoration de la somme $ \sup + \inf $, on utilise la formule de repr\'esentation de Green et une in\'egalit\'e de Sobolev  pour la sph\`ere de dimension 2. D'autre part nous estimons le minorant grace \`a une in\'egalit\'e de Sobolev optimale. Pour ce qui concerne les surfaces de Riemann compacte sans bord, on utilise un blow-up, la nature de la fonction de Green (=noyau de Newton + fonction continue) en dimension 2 et une repr\'esentation int\'egrale des solutions.

S'agissant des dimension $ \geq 3 $, Nous signalons entre autres, des r\'esultats connus de YY.Li et, CC. Chen et CS.Lin, le premier fait appel au proc\'ed\'es de s\'el\'ection permettant d'aboutir a la notion de blow-ups isol\'es et isol\'es simples. Gr\^ace a la formule de Pohozaev et des conditions de platitudes des courbures prescrites les auteurs prouvent certaines esimees a priori: ils mettent en evidence les notions de blow-ups isoles simples, une certaine suite  de fonctions converge vers une fonction semblable a la fonction de Green, la formule de Pohozaev permet de conclure. Notons que il y a une autre methode, non citee ici, et basee sur la methode moving-plane ainsi que des estimees asymptotiques des fonctions blow-up, c'est celle de YY.Li

La methode de blow-ups isol\'es et isol\'es simples est d'une grande utilit\'e dans la comprehension du probl\`eme de compacit\'e des solutions de l'\'equation de Yamabe.  

Concernant la minoration du produit  $ \sup \times  \inf  $, nous utilisons le proc\'ed\'e d'iteration de Nash-Moser combin\'e \`a une estimation de la fonction de Green de l'operateur inversible $ -\Delta + \epsilon $, par valeurs inferieures, ici les in\'egalit\'es  de Sobolev et celles de Harnack sont utiles, ainsi que le principe du maximum. Pour ce dernier point  on se place sur des  vari\'et\'es compactes sans bord.

 Le cas des vari\'et\'es compactes avec bord est trait\'e dans deux directions, la premi\`ere concerne les \'equations du type Yamabe sans condition de Dirchelet et la seconde, avec condtion de Dirichelet. Pour la premi\`ere, une minoration de la fonction de Green est necessaire et un proc\'ed\'e d'iteration de Nash-Moser est appliqu\'e. Notons qu'ici on a impose une condition de saut \`a l'interieur du domaine. Cela revient a dire que toute suite de fonctions converge vers 0 uniformement sur tout compact, soit on a une minoration du produit $ \sup  \times \inf $. Pour ce qui concerne le probleme avec condition de Dirichelet, nous utilisons un blow-up combin\'e \`a l'utilisation du principe du maximum pour comparer les fonctions au noyau de Green.

\bigskip

Il y a d'autres resultats importants, qu'on a pas cit\'e. On donne ici un aperçu des estim\'ees uniformes pour l'equation de Yamabe, de type Yamabe, de la courbure prescrite, de type courbure prescrite. De Liouville et de type Gauss.

\bigskip

\chapter{Le Cas de la dimension 2: la quantisation.} 

%\underbar {\bf I) Cas de la dimension 2:} 

\bigskip

On s'int\'eresse ici aux estimations a priori du type:

 $$ \sup_K u +C_1\inf_M u <C_2,  \,\,\, {\rm et } \,\,\, \sup_M u + C_1 \inf_K u > C_2 $$

ou $ u $ est solution de l'\'equation de la courbure prescrite, sur une surface $ M $,  \`a savoir,

 $$ -\Delta u+ R=V(x)e^u, \qquad (E) $$ 

avec, $ R  $ la courbure scalaire et $ V $ la courbure scalaire prescrite.

\bigskip

Dans le cas $ M =  \Omega $ ouvert de $ {\mathbb R}^2 $, l'\'equation se r\'eduit  :

$$ -\Delta u=V(x)e^u, \qquad (E') $$

Concernant cette \'equation, de nombreux r\'esultats existent  et nous nous interessons aux estimations a priori qui en d\'ecoulent.

\bigskip

\section{In\'egalit\'es de Bandle-Suzuki: une sorte de principe du maximum.} 

~~\\

1-Un des r\'esultats de C. Bandle est le suivant: (une sorte de principe du maximum)

\smallskip
 
 \begin{Theorem} Si $ u $ est solution de $ (E') $ sur $ \Omega= B=B_R(0) $, avec :
 
 $$ V(x)\leq \lambda, $$
 
 $$ \Sigma= \int_B  \lambda e^u dx  < 8 \pi, $$
 
 Alors,
 
 $$ u(0) \leq \max_{B} u \leq \max_{\partial B}  u -2 \log(1- \dfrac{\Sigma}{8\pi}). $$
 
 \end{Theorem}

 2-Un des r\'esultat de T. Suzuki est le suivant: (une in\'egalit\'e de Harnack)

\smallskip

 \begin{Theorem} Si $ u $ est solution de $ (E') $ sur $ \Omega= B=B_R(0) $, avec :
 
 $$ V(x)\leq \lambda, $$
 
 $$ \Sigma= \int_B  \lambda e^u dx  < 8 \pi, $$
 
 Alors,
 
 $$ u(0) \leq \dfrac{1}{2\pi R} \int_{\partial B} u d\sigma -2 \log(1- \dfrac{\Sigma}{8\pi}). $$
 
  \end{Theorem}

Voir [15-49]
 
 3- Sur le resultat de quantization de T. Suzuki:
 
 \smallskip
 
 Commencons par donner la defintion d'un point blow-up pour une suite de fonctions.
  
  \smallskip
   
 \begin{Definition}. On dit qu'un point $ x_0 $ est un point blow-up de la suite de fonctions $ (u_k) $ s'il existe une suite de points $ (x_k)_k $ telle que :
 
 a) $ x_k \to x_0 $,
 
 b) $ u_k(x_k) \to + \infty $. 
 
\end{Definition} 
 
 \bigskip
  
 \begin{Theorem} Soit $ u $ solution du probleme suivant:
 
 $$ \begin{cases}

-\Delta u=\lambda f(u), \,\, {\rm dans } \,\, \Omega,\,\,\lambda >0\\
             u=0 \,\, {\rm sur } \,\, \partial \Omega, 
             
 \end{cases}  $$
 
 On suppose que $ \Omega $ est regulier et  $ f >0 $ de la forme:
 
 $$ f(t)=e^t+g(t),\,\, |g(t)|=o(e^t),  \,\, t \to +\infty $$
 
 et,
 
 $$  |g'(t)-g(t)|\leq G(t),\,\,  G(t)+|G'(t)|=O(e^{\gamma t}), \,\, t \to +\infty, \,\,  \gamma < 1/4. $$
 
 Alors:
 
 \smallskip
 
  La quantit\'e $ \Sigma= \int_B  \lambda f(u) dx   $ tend vers $ 8 \pi k$, avec $ k=0,1, 2,..., +\infty $ quand $  \lambda \to 0 $. Les solutions $ u $ verifient l'une des conditions suivantes (quand $ \lambda \to 0 $):
 
 \smallskip
 
 1-Si $ k=0 $, $ ||u||_{\infty}  \to 0 $
 
 \bigskip
 
  2-Si $ 0 < k < +\infty  $, il existe un nombre fini de points blow-up, en dehors des quels  $ u $ est localement uniformement  borne.
 
 \smallskip
 
  3-Si $ k=+ \infty $, alors, $ u $ diverge vers  $ + \infty $.
 \end{Theorem}
\bigskip
 
 Esquisse de la preuve:
 
  \bigskip 
 
 Le cas  1 est clair par les estimations elliptiques, l'inegalite de Suzuki (qu'on a presente avec l'inegalite de Bandle) et les estimations au bord de De Figueiredo-Lions-Nusbaum. Le cas 3, on utilise la comparaison de la fonction de Green  la premeire fonction prorpre et la formule de representation de Green. Pour ce qui est du cas 2, l'auteur procede comme suit:
  
  \bigskip 
  
   Etape 1:
  
  \smallskip 
    
Le fait que $ \Omega $ soit regulier, un resultat de De Figueiredo-Lions-Nussbaum (methode moving-plane) implique que $ u $ est uniformement bornee sur un voisinage du bord. S'il y a blow-up, il est a l'interieur de l'ouvert.
 
  \smallskip 
  
   Etape 2: 
  
   Soit $ v=e^{-u/2} $ et $ S  $ la fonction complexe, $ \partial_{zz}u-(1/2)(\partial_z u)^2 $, alors $ v $ est solution de $ \partial_{zz}v+ S v = 0 $ avec $  |\partial_{ \bar z} S|_{\infty} \in O(1) $. S converge alors vers une fonction holomorphe et asymptoiquement $ v $ s'ecrit comme combinaison lin\'eaire de deux fonctions voisines de fonctions holomorphes  (donc, dont les zeros sont isol\'es, c'est la localisation des points blow-ups):
 
 $$ v(x)=v(x_0)u_1+v(x_0)\lambda f(u(x_0)) u_2 $$
 
 avec, $ u_1 $ et $ u_2 $ voisines en norme uniforme de deux fonctions analytiques  et $ x_0 $ le point blow-up. Les co\'efficients $ v(x_0) $ tend vers 0 alors que le deuxieme non, sinon $ u $ tenderait vers l'infini partout, ce qui n'est pas possible au voisinage du bord. En dehors du point blow-up, mais au voisinage de ce point $ v $ est uniformement bornee donc $ u $.

\smallskip 

Voir [19] des references supplementaires et, [28] et[49] de la bibliographie.

 \bigskip
 
{\bf Remark 1:} Le passage en complexe est motive par le fait que lorsque $ f=e^t $, on a la celebre solution fondamentale de l'equation de Liouville.

 \smallskip
 
{\bf Remark 2:} Concernant les estimations de De Figueiredo-Lions-Nussbaum: ce sont des estimations uniformes au voisingage du bord, on les prouve grace a la m\'ethode moving-plane. C'est essentiellement une transformationde Kelvin, pour rendre le voisinage du bord convexe, et une application du principe du maximum (suivant la direction normale).

\bigskip

\section{In\'egalit\'es  d'Alexandrov-Bol-Fiala-Bandle: in\'egalit\'es isop\'erimetriques.}

~~\\

1-L'inegalite isoperimetrique classique:

\begin{Theorem} Soit $ C $ une courbe $ C^1$ fermee de $ {\mathbb R}^2 $ delimitant un domaine $ \Omega $. On note $ L $ la longeur de $ C $ et $ A $ l'aire de  $ \Omega $. Alors:

$$ L^2 \geq 4\pi A. $$

L'egalite a lieu si et seulement si $ C $ est un cerle.
\end{Theorem} 

\bigskip

Il exsite plusieurs preuve de l'inegalite isoperimetrique. Nous renvoyant a celle d'Hurwitz qui utilise l'inegalite de Wirtinger.

\bigskip

2-Sur une boule $ B = B_R $, si on note,

$$ M_1= \int_B  e^u dx \,\, {\rm ,} \,\,  L_1=\int_{\partial B}   e^{u/2} d\sigma, $$

et,

$$ \omega_{K_0}^+(B)=\int_{ \{ x \in B, K(x) > K_0 \} } (K(x)-K_0) e^{u} dx, $$

alors on a,

\bigskip

\begin{Theorem} Pour $ K_0 $ reel tel que $ \alpha:=2\pi -\omega_{K_0}^+(B) >0 $ et $ K_0M_1 <2 \alpha $, on a,

$$ L_1^2 \geq ( 2 \alpha- K_0M_1)M_1, $$

Dans le cas paticulier ou $ K(x) \leq K_0 $ sur $ B $, on a:

$$ L_1^2 \geq ( 4 \pi- K_0M_1)M_1, $$

 \end{Theorem} 
 
 Voir [15]

1- La methode Fiala:

\bigskip

La construction que fait Fiala, pour prouver l'in\'egalit\'e isop\'erim\'etrique,se base sur une d\'efinition des coordonnees geodesiques paralleles. Il considere une courbre $ C $ d\'elimtant une domaine simplement connexe. Cette courbre $ C $ est suppos\'ee r\'eguliere (analytique) dans un espace de Riemann complet ( Hopf -Rinow qui permet d'avoir entre deux points une g\'eodesique minimsante et donc definir ces coordonn\'ees.) puis s'occupe de d\'efinir des points extremaux, points ou les g\'eodesiques ne sont plus minimisantes, points focaux, afocaux, qui sont en nombre fini sur tout eensemble born\'e. A partir de la, on peut parcourir notre espace  partir de la courbe, en \'etant  une distance r de celle-ci et ne rencontrer qu'un nombre fini de points extremes stationnaires.

\bigskip

Si on note $ L(p) $ la longeur de la vraie parallele alors le but est de prouver que:

$$  \dfrac{dL(p)}{dp}  \leq  \int_0^{L(C)} k(q) dq-C(p), \,\, si \,\, p >0 $$

$$  \dfrac{dL(p)}{dp}  \geq  \int_0^{L(C)} k(q) dq + C(p), \,\, si  \,\, p_{min} < p < 0 $$

ou, $ C(p) $ est l'integrale de la courbure totale sur le domaine compris entre la courbe $ C $ et la vraie parallele.

\bigskip

Ceci, se fait grace  l'expression de $ L(p) $ (qui est continue) en coordonnes g\'eodesiques paralleles (qui sont semblables aux coordonnees de Fermi) et la formule de Gauss-Bonnet (ici, la caracteristique d'Euler-Poincare est 1, car $ C $ delimite un domaine simplement connexe $ F $. Si on note $ K $ la courbure de Gauss, on obtient:

$$ 2\pi =  \int_0^{L(C)} k(q) dq+ \int_F K(x) dx $$

Les deux in\'egalit\'es precedentes, peuvent etre \'ecrites sous la forme suivante:

$$  \dfrac{dL(p)}{dp}  \leq  2\pi - \int_{F_p} K(x) dx, \,\, si \,\, p >0 $$

$$  \dfrac{dL(p)}{dp}  \geq  2\pi - \int_{F_p} K(x) dx, \,\, si  \,\, p_{min} < p < 0 $$

Avec $ F_p, p >0 $, le domaine totale, union de celui delimite par $ C $ et celui entre $ C $ et la vraie parallele. De meme pour $ p <0 $.

2- La methode Bandle:

\bigskip

La m\'ethode de Bandle est celle des lignes de niveaux. Elle consiste a obtenir une in\'equation differentielle pour la symmetris\'ee de Schwarz d'une certaine fonction. Elle utilise aussi l'in\'egalit\'e isop\'erim\'etrique de Nehari qui, pour une foncton harmonique $ h $ definie sur un domaine $ D $ simplement connexe, est de la forme:

$$ 4\pi \int_D e^h dx \leq \left ( \int_{\partial D} e^{\frac{h}{2}} d\sigma \right )^2, $$

L'\'egalit\'e dans la m\'ethode de Bandle a lieu seulement pour une classe particuliere de fonctions $ u $.

\bigskip

\section{Quantisation et r\'esulats de Brezis-Merle: convergence faible dans $ L^1 $.}
~~\\
\begin{Theorem}On consid\`ere deux suites de fonctions $ (u_n, V_n) $ solutions de $ (E') $, avec,

$$ V_n \geq 0, \,\, ||V_n||_{L^{p}(\Omega)} \leq C_1, \,\,  ||e^{u_n}||_{L^{p'}(\Omega)} \leq C_2, $$

avec, $ 1 < p \leq + \infty $. Alors, on a, ou bien

$$ \sup_K |u_n|  \leq c=c(K, C_1, C_2, \Omega) $$

ou bien,

$$ u_n \to -\infty, \,\, {\rm sur \,\, tout \,\,  compact \,\,  de} \,\, {\Omega}, $$

ou bien,

Il existe un ensemble fini de points $ S $ tel que $ u_n \to -\infty $ uniform\'ement sur tout compact de ${\Omega}-S $ et $ u_n \to +\infty $ sur $ S $ et, au sens faible, on a,

$$ V_ne^{u_n} \to \sum_i\alpha_i \delta_{a_i}, \,\, {\rm avec } \,\, \alpha_i \geq  4\pi/p', \,\, {\rm et} \,\, a_i \in  S. $$
\end{Theorem} 

C'est la convergence faible dans $ L^1 $ et vers une mesure positive,  de:

$$  \int V_ne^{u_n} \rightharpoonup \mu. \,\, {\rm au \, sens \, faible } \,\,  L^1(\Omega) $$

qui permet d'avoir l'alternative. le principe de Harnack et le maximum dans $ W_0^{1,1}(B_R(x_0)) $  ($ x_0 $ un point "blow-up") obtenu par l'inegalite de Kato permettent de de conclure.

\smallskip

Comme corollaire de ce r\'esultat, on a:
\begin{Corollary}

$$  \sup_K u  \leq c=c(\inf_{\Omega} V, ||V||_{L^{\infty}(\Omega)}, \inf_{\Omega} u, K, \Omega). $$
\end{Corollary}

Voir [19]

\begin{Theorem}{\it(Li-Shafrir.)}On consid\`ere deux suites de fonctions $ (u_n, V_n) $ solutions de $ (E') $, avec,

$$  0 \leq V_n \leq C_1, \,\,  ||e^{u_n}||_{L^{1}(\Omega)} \leq C_2, $$

avec, $ V_n \to V $ dans $ C^0(\bar \Omega) $.  On suppose que le troisieme cas de l'alternative de Brezis-Merle est verifie (c'est a dire que la suite blow-up), alors:

$$ \int V_ne^{u_n} \rightharpoonup \sum_i\alpha_i \delta_{a_i}, \,\, {\rm avec } \,\, \alpha_i = 8 m_i \pi, \, m_i\in {\mathbb N}^*, \, {\rm et} \,\, a_i \in  S. $$

\end{Theorem} 

Voir [38]

\bigskip

\section{Quantization au bord.}

~~\\

Ici, on ajoute une condition au bord.

\begin{Theorem}(Bahoura) On consid\`ere deux suites de fonctions $ (u_n, V_n) $ solutions de $ (E') $, avec,

$$ u_n = 0  \,\, {\rm sur } \,\, \partial \Omega, $$

$$ V_n \geq 0, \,\, ||V_n||_{\infty} \leq C_1, \,\,  ||e^{u_n}||_{L^{1}(\Omega)} \leq C_2, $$

 alors, on a, 

\bigskip

Il existe un ensemble fini de points $ S \subset \partial \Omega $ tel que $ u_n \to u $ uniform\'ement sur tout compact de $ \bar \Omega-S $ et $ u_n \to +\infty $ sur $ S $ et, au sens faible, on a,

$$ V_ne^{u_n} \to Ve^u + \sum_i\alpha_i \delta_{a_i}, \,\, {\rm avec } \,\, \alpha_i \geq  4\pi, \,\, {\rm et} \,\, a_i \in  S. $$
\end{Theorem} 

Comme corollaire, on a, si $ V_n $ est supposee uniform\'ement Lipschitzienne, alors,

\begin{Corollary}(Bahoura)

$$ S=\{\emptyset \}, $$

et,

$$ ||u_n||_{\infty} \leq c=c(C_1, C_2,||\nabla V_n||_{\infty},  \Omega), $$

\end{Corollary}

Voir [12]

Esquisse de La preuve: 

Deux estimations uniformes sont necessaires pour obtenir ce resultat:

$$ \int_{\partial \Omega}\partial_{\nu} u_n d\sigma \leq C \,\, {\rm et } \,\,  \int_{\Omega}|\nabla u_n|^q dx \leq C,\,\, 1\leq q < 2. $$

D'une part on a convergence faible dans $ L^1 $ des deriv\'ees normales, d'autre part l'estim\'ee $ L^q $ permet de pr\'eciser la convergence.

$$  \int \partial_{\nu} u_n d\sigma \rightharpoonup \mu  \,\, {\rm dans } \,\,  L^1(\partial \Omega). $$
 
 On conclut grace a la formule de Pohozaev. (Pour une bonne formulation du probleme on suppose les elements de la suite dans $ W_0^{1,1}(\Omega) $).
 
 \smallskip 
 
\underbar {Remarque importante:}

\smallskip

Il existe une autre preuve de ce dernier resultat (compacit\'e lorsque $ V_n $ est lipschitzienne) par la methode "moving-plane". Elle consiste a prouver que les points crititques sont loin du bord, en se placant  suivant la normale, apres avoir utiliser la transformee de Kelvin. Voir les articles de Chen-Li et Ma-Wei [24] et [42].

\smallskip

On a le resultat de la compacit\'e globale de Chen-Li. Le resultat de Chen-Li utilise le fait qu'on a compacit\'e au voisinage du bord lorsque $ ||\nabla \log V||\leq A $, puis il etend ce resultat lorsque $ ||\nabla V||\leq A_1 $, pour cela il utilise l'extension des resultats de Brezis-Merle (Theorem 1 de Brezis Merle), qui reste vrai dans des domaines Lipschitzien d\`es qu'on a la regularit\'e des solutions dans $ W_0^{1,2} $, car le principe du maximum est valable dans ce cas (on utilise l'integration par parties, qui est vrai des que la regularit\'e du bord est Lipschitzien. Pour l'inegalit\'e de Sobolev aussi et la resolution d'un probleme variationnel dans $ L^2 $. Lipschitz suffit). (et aussi l'extension des fonctions harmoniques et la formule de la moyenne). Voir la preuve du corollaire de l'article de Chen-Li.

\smallskip

Remarque sur la preuve de Chen-Li : pour etendre la partie $ u_1 $ harmonique, il faut supposer le domaine analytique, pour pouvoir utiliser une transformation conforme qui reste invariante par le Laplacien. Principe de symetrisation de Schwarz. Donc, le resultat de compacit\'e reste vrai avec la regularit\'e « smooth » lorsque on suppose  $ ||\nabla \log V||\leq A $. Mais la regularit\'e du domaine doit etre suposs\'ee analytique lorsqu'on passe a $ ||\nabla V||\leq A_1 $.

\smallskip

Dans leur preuve Chen-Li, utilisent le fait que l'operateur est invariant par appplication de carte, ceci est possible si cette application est conforme, elle preserve le Laplacien. Puis, symetrise la fonction en symetrisant un probleme de Dirichlet, puis soustrayent les valeurs aux bords et ils obtiennent l'image de $ v_1 $. Alors $ u_1 $ est l'image de $ v_1 $ par l'application de carte. Maintenant pour construire $ v_1 $ ils utilisent une symetrisation d'un probleme de Dirichlet, qui requiert les solutions dans $ W^{2,p}\cap C^2(B_{\epsilon}) \cap C^1(\bar B_{\epsilon}), p >2 $ (la formule de representation de Green reste valable, dans ce cas, voir la preuve dans Gilbarg-Trudinger). Puis, ils utilisent la formule integrale de Poisson (qui necessite d'avoir l'operateur Laplacien).

\smallskip

a) Pour utiliser la formule de Poisson, on conserve le Laplacien : transformation conforme $\phi $.
\smallskip

b) Ils symetrisent $ uo\phi $ ils obtiennet une fonction $ u_v \in C^1(\bar B_{\epsilon}(0)) \cap W^{2,p} $.
\smallskip

c) Ils resolvent : $ -\Delta v_1=-\Delta u_v $ avec condition de Dirichlet sur $ B_{\epsilon}(0) $.
\smallskip

d) Ils utilisent la formule integrale de Poisson pour $ v_1-u_v \in W^{2,p}\cap C^2(B_{\epsilon}) \cap C^1(\bar B_{\epsilon}), p >2 $.
Sur le bord, il n'y a que les valaurs de $ u $.  

(Ce travail revient \`a symetriser une fonction harmonique qui n\'ecessite le theoreme de symetrisation de Schwarz
, qui necessite une application conforme, donc un domaine de depart $ \Omega $ analytique).

\bigskip

\underbar {\it Questions (Problemes ouverts de Brezis-Merle):}

\bigskip

On considere sur un domaine $  \Omega $ de $ {\mathbb R}^2 $:

$$ \begin{cases}

-\Delta u=V(x)e^u, \,\, {\rm dans } \,\, \Omega,\\

             u=0 \,\, {\rm sur } \,\, \partial \Omega, 
             
 \end{cases}  $$
 
On suppose:

$$  0 < a \leq V \leq b $$

1- A t on :

$$  \sup_{\Omega} u \leq c(a,b,\Omega) ? $$

2- Peut on avoir l'inegalite precedente si on suppose la condition de Brezis-Merle:

$$ \int_{\Omega}  e^u dx  \leq  C_1 ? $$

3- On suppose $ a=0 $ et $ V_n \to V $ dans $ C^0(\Omega) $. A t on la meme estimation uniforme que precedemment ?

\chapter{Le Cas de la dimension 2: In\'egalit\'es du type $ \sup+C\inf $.}

\bigskip

Sur un ouvert $ \Omega $ de $ {\mathbb R}^2 $, on considere l'equation suivante :

$$ -\Delta u=V(x)e^u, \qquad (E') $$

\section{R\'esultat de Shafrir.}

~~\\

\begin{Theorem} Si, on suppose $  0 < a \leq V(x) \leq b < + \infty $

$$ C \sup_K u + \inf_{\Omega} u \leq c=c(a,b, K, \Omega), $$

avec, $ C=C \left ( \dfrac{a}{b} \right ) $

\end{Theorem}

\bigskip

Pour prouver son resultat, Shafrir utilise l'in\'egalit\'e isop\'erimetrique d'Alexandrov-Bol-Fiala-Bandle et la fonction "blow-up" suivante, 

$$ G(r)=u(0)+ \dfrac{C_1}{2\pi r}\int _{\partial B_r} u d\sigma + 2(C_1+1)\log r, $$

Voir [46]

\bigskip

\section{R\'esultat de Brezis-Li-Shafrir.}

~~\\

La methode "moving-plane"  et un passage en coordonn\'ees polaires a permis a Brezis-Li-Shafrir de prouver que,

\bigskip

\begin{Theorem} Si, on suppose $  0 < a \leq V(x) \leq b < + \infty $ et que $ V $ est uniform\'ement Lipschitzienne,

$$ \sup_K u + \inf_{\Omega} u \leq c=c(a,b, ||\nabla V||_{\infty}, K, \Omega). $$
\end{Theorem}

Voir [18]

Notons que la methode moving-plane sous sa forme actuelle a ete introduite par James Serrin et Louis Nirenberg.

\bigskip

Dans la preuve on utilise le resultat de Classification par la methode moving-plane des solutions d'une EDP elliptique nonlineaire :

\begin{Theorem}(Chen-Li)([23]) Les solutions de:

$$ \begin{cases}

-\Delta u=Ke^u, \, u(0)=0, \, {\rm dans } \,\, {\mathbb R}^2 \\

            \int_{{\mathbb R}^2 } e^u <+ \infty. 
             
 \end{cases}  $$
 
Sont:

$$   u(x) = -2\log(1+\gamma |x|^2), \,\,\gamma= \sqrt {(K/8)} $$

\end{Theorem}

\bigskip

\section{R\'esultat de Chen-Lin.}

~~\\

La m\'ethode de symmetrisation, une in\'egalit\'e differentielle g\'eom\'etrique appliqu\'ees aux fonctions "courbures int\'egrales"  et l'in\'egalit\'e de Suzuki a permis a Chen-Lin de prouver que,

\begin{Theorem} Si, on suppose $  0 < a \leq V(x) \leq b < + \infty $

$$ \sqrt {\dfrac{a}{b}} \sup_K u + \inf_{\Omega} u \leq c=c(a, b, K, \Omega). $$

\bigskip

Si de plus, $ V $ est supposee uniform\'ement holderienne, on a, 

$$ \sup_K u + \inf_{\Omega} u \leq c=c(a,b, ||V||_{C^{\alpha}(\Omega)}, K, \Omega). $$

avec $ \alpha \in (0, 1] $.

\end{Theorem}

\bigskip

Notons que dans leur resultat Chen-Lin ont classifier les EDP associ\'ees sur $ {\mathbb R}^2 $ suivant les valeurs de leur courbure totale(ou int\'egrales).

\bigskip

Voir [26]

Esquisse de la preuve:

\bigskip

1) Symmetrisation

\bigskip

 Si $ u^* $ est la symmetrise de Schwarz de $ u $, alors, on considere les courbures integrales:

$$ F(r)=\int_{\Omega_{u^*(r)}} V(x)e^{u(x)} dx \,\, {\rm et } \,\, \bar K(r)=\dfrac{F'(r)}{2\pi r}e^{-u^*(r)}, $$

avec,

$$ \Omega_{u^*(r)}= \{x, u(x)>u^*(r)\}, \,\, {\rm qu'ils \,comparent \,a \,des \,boules.} $$

Alors, $ F(r) $ est solution de l'inequation differentielle suivante:

2) In\'egalit\'e differentielle g\'eom\'etrique et utilisation de l'in\'egalit\'e de Suzuki:

$$ \dfrac{rF'(r)}{ \bar K(r)} \geq \dfrac{1}{4\pi a} \left ( F(r) - 4\pi \left (1-\sqrt{\dfrac{a}{b}} \right ) \right ) \left ( 4\pi \left (1+\sqrt{\dfrac{a}{b}} \right )- F(r)\right ) $$

3) Cas ou la courbure int\'egrale est minimale:

\bigskip

Dans le cas ou la courbure est minimale, on a l'egalite suivante sur $ {\mathbb R}^2 $ :

$$\int_{{\mathbb R}^2}|\nabla u|^2 dx =\int_{{\mathbb R}^2}|\nabla u^*|^2 dx, $$

qui entraine que $ u $ est radiale.

Voir [6] et [7] des references supplementaires.

\bigskip

\section{Une nouvelle preuve dans le cas Holderien.}

~~\\

La m\'ethode moving-plane, nous permet de prouver, sans r\'esultat de classification d'EDP sur $ {\mathbb R}^2 $, l'in\'egalit\'e suivante:

\begin{Theorem}(Bahoura)([13])  $$ \sup_K u + \inf_{\Omega} u \leq c=c(a, b, ||V||_{C^{\alpha}(\Omega)}, K, \Omega). $$

avec $ \alpha \in (0, 1] $.
\end{Theorem}

\bigskip

Notons que la m\'ethode moving-plane sous sa forme actuelle a ete introduite par James Serrin et Louis Nirenberg.

\bigskip

\underbar {\it Question (Probleme ouvert): Brezis-Li-Shafrir\\}

\bigskip

Ici, on suppose que $ V $ une fonction continue  et positive sur $ \Omega $, peut on avoir l'estimation suivante:

$$ \sup_K u + \inf_{\Omega} u \leq c=c(||V||_{C^{0}(\Omega)}, K, \Omega) ? $$

\bigskip

\section{Le cas de la sph\`ere $ {\mathbb S}_2 $.}

~~\\
Sur la sph\`ere unit\'e de dimension 2, on consid\`ere l'\'equation suivante:

$$ \Delta u + K e^{2u} = 1. $$

\begin{Theorem}{\it (Chang-Gursky-Yang)}. Si $ K $ v\'erifie, $ 0 < m \leq K \leq M $, alors:

$$ \left |\int_{{\mathbb S}_2} |\nabla u|^2 + 2 \int_{{\mathbb S}_2} u \right | \leq C(m, M). $$
\end{Theorem}

\underbar {\bf Obstruction du type Kazdan-Warner: g\'en\'eralisation}

\bigskip

\begin{Theorem} {\it (Bourguignon-Ezin)}. Pour un champ de vecteurs conforme $ X $ sur une vari\'et\'e Riemannienne compacte $ (M,g) $, on a:

$$ \int_M X(R) dV_g =0 .$$

\end{Theorem}

\bigskip

\begin{Theorem} {\it ( Kazdan-Warner)}. On note $ \xi $ un \'el\'ement de l'espace vectoriel correspondant \`a la prem\`ere valeure propre et $ R $ la courbure scalaire, alors:

$$ {\rm Pour } \,\, {\mathbb S}_2, \,\, \int_{{\mathbb S}_2} \nabla^i \xi \nabla_i R e^u dV_g =0 , $$

\end{Theorem}

\begin{Theorem}{\it (Chang-Gursky-Yang)}. Si $ K $ v\'erifie, $ 0 < m \leq K \leq M $, $ \Delta K \not = 0 $, alors:

$$ |u|_{L^{\infty}}\leq c=c(m, M,|K|_{C^2}). $$
\end{Theorem}

Voir [22]

Les estimations a priori de Chang-Gursky-Yang ont pour but de resoudre le probleme de la courbure prescrite sur la sphere (Probleme de Nirenberg) par la methode de continuite. Pour des courbures prescrites proches des constantes.

\bigskip

on a:

\begin{Theorem}(Bahoura)([4-5]) Sur la sph\`ere $ {\mathbb S}_2 $, on consid\`ere l'\'equation suivante,

$$ -\Delta u+ 2=V(x)e^u,  $$%\qquad (E) $$ 

avec, 

$$ 0 \leq V(x)\leq b, $$ 

alors,

$$ \sup_{{\mathbb S}^2} u + \inf_{{\mathbb S}^2} u \geq c=c(b), $$

\end{Theorem}

Comme corollaire, on a:

\begin{Corollary}(Bahoura) ([4-5]) De plus, on peut estimer la constante $ c=c(b) $:

$$ -2+2\log 2-2\log b \leq c \leq 2\log 2-2\log b, $$

\end{Corollary}

\bigskip

Ces r\'esultats sont bas\'es sur les in\'egalit\'es de Sobolev et la formule de repr\'esentation de Green.

1- In\'egalit\'e d'Aubin: Sur $ {\mathbb S}^2 $, on a:

$$ \int_{{\mathbb S}^2} e^u d\sigma \leq C(\epsilon) e^{(\frac{1}{32\pi} + \epsilon) \int_{{\mathbb S}^2} |\nabla u|^2 d\sigma +\frac{1}{4\pi} \int_{{\mathbb S}^2} u d\sigma }, $$

pour toute $ u \in W^{1, 2}({\mathbb S}^2) $, telle que:

$$ \int_{{\mathbb S}^2} e^u x_j d\sigma = 0, $$

avec, les $ x_j $ les fonctions composantes de $ {\mathbb S}^2 $, qui sont les valeurs propres du laplacien sur  $ {\mathbb S}^2 $.

\bigskip

L'in\'egalit\'e d'Aubin est bas\'ee sur l'in\'egalit\'e de Moser, qui stipule:

$$ \int_{{\mathbb S}^2} e^u d\sigma \leq K  e^{\frac{1}{16\pi} \int_{{\mathbb S}^2} |\nabla u|^2 d\sigma +\frac{1}{4\pi} \int_{{\mathbb S}^2} u d\sigma } $$

pour toute $ u \in W^{1, 2}({\mathbb S}^2) $.

\bigskip

2-In\'egalit\'e d'Onofri:

$$ \int_{{\mathbb S}^2} e^u d\sigma \leq 4\pi e^{\frac{1}{16\pi} \int_{{\mathbb S}^2} |\nabla u|^2 d\sigma +\frac{1}{4\pi} \int_{{\mathbb S}^2} u d\sigma } $$

pour toute $ u \in W^{1, 2}({\mathbb S}^2) $. C'est l'in\'egalit\'e de Moser avec $ K=4\pi $.

\section{Sur la minoration de la somme $  \sup + C \inf $: cas des surfaces compactes sans bord.}
~~\\
Le resultat suivant est du a Cheeger et Gromov pour les surfaces de Riemann compactes, qu'on peut appliquer a une suite de metriques conformes.

\begin{Theorem} (Cheeger-Gromov). Considerons une suite de surfaces compactes sans bord $ (M, g_j) $ telles que leur courbures sectionelles $ K_j $ verifient:

$$ 0 < c_1 \leq K_j \leq c_2 < +\infty $$

alors cette suite contient une sous-suite convergente a des diffeomorphismes pres. De plus, le rayon d'injectivite $ i_M $  est tel que:

$$ i_M \geq c(c_1,c_2) >0. $$

\end{Theorem}

Notons que dans [52], on a un principe de concentration-compacite pour les suites de surfaces compactes

On peut etendre le resultat de la sphere a une surface compacte sans bord pour des valeurs particulieres de la courbure scalaire en utilisant le theoreme suivant:

\bigskip

\begin{Theorem}{\it (Fontana)}. Sur une vari\'et\'e Riemannienne compacte sans bord, on consid\`ere une fonction $ \phi \in H_1^2(M) $ telle que, $ \int_M \phi = 0 $, alors:

$$ \int_M e^{\phi} dV_g \leq C \exp \left ( \dfrac{1}{16 \pi} ||\nabla \phi ||_2^2 \right ). $$
\end{Theorem}

\bigskip

Le resultat suivant est une consequence d'un theoreme de YY.Li sur les surfaces compactes sans bord pour l'equation:

$$ -\Delta u+ R=V(x)e^u, \qquad (E) $$

avec, 

$$ 0 < a \leq V(x)\leq b,  \,\,\, ||\nabla V||_{L^{\infty}} \leq A, $$

\begin{Theorem} (YY. Li [37]) On a: $$   |\sup_{M} u + \inf_{M} u | \leq c=c(a, b, A, M, g), $$

\end{Theorem}

Nous obtenons une g\'en\'eralisation du r\'esultat precedent, sur une surface de Riemann compacte sans bord $ (M, g) $ de volume 1, en consid\'erant l'\'equation suivante,

$$ -\Delta u+ k=V(x)e^u,  $$

avec, 

$$ 0 \leq V(x)\leq b, $$ 

Alors,

\begin{Theorem} (Bahoura)$$  \dfrac{k-4\pi}{4 \pi} \sup_{M} u + \inf_{M} u \geq c=c(k, b, M, g), $$

\end{Theorem}

Voir [9]

Esquisse de la preuve:

\bigskip

1-Technique "blow-up":

\bigskip

2-Nature de la fonction de Green, en dimension 2 et principe du maximum:

\bigskip

\underbar{ \bf Fonctions de Green du Laplacien}

\bigskip

\underbar{ \it 1. Cas des vari\'et\'es Riemanniennes sans bord}

Soit $ M $ une vari\'et\'e Riemannienne compacte de dimension $ n \geq 2 $. On appelle fonction de Green du Laplacien, une fonction $ G $ d\'efinie sur $ M \times M $ et r\'eguli\`ere en dehors de la diagonale, telle que:

$$ -\Delta_g G(x,.) = \delta_x - \dfrac{1}{V_g} \,\, {\rm au \,\, sens \,\, des \,\, distributions}. $$

O\`u $ V_g $ est le volume de $ M $.

\smallskip

On peut prendre $ G(x,y) $ telle que:

$$ G \geq 0, \,\, \int G \equiv cte. $$

\bigskip

\underbar{ \it 2. Cas des vari\'et\'es Riemanniennes avec bord}

\bigskip

Soit $ M $ une vari\'et\'e Riemannienne compacte avec bord de dimension $ n \geq 2 $. On appelle fonction de Green du Laplacien, une fonction $ G $ d\'efinie sur $ M \times M $ et r\'eguli\`ere en dehors de la diagonale, telle que (au sens des distributions):

$$ -\Delta_g G(x,.) = \delta_x, \,\, {\rm sur} \,\, M, $$

$$ G(x,.) = 0, \,\, {\rm sur} \,\, \partial M. $$

\underbar {\it 3. Propri\'et\'es des fonctions de Green}

\bigskip

$$ G(x,y)=G(y,x), $$

$$ G(x,y) \leq c_0 \log | d_g(x,y)| \,\, {\rm si} \,\, n=2, \,\, G(x,y) \leq \dfrac{c_1}{{[d_g(x,y)]}^{n-2}} \,\, {\rm si} \,\, n \geq 3, $$

$$ |\nabla_x G(x,y)|\leq \dfrac{c_2}{{[d_g(x,y)]}^{n-1}}, $$

$$ |\nabla_{x,x} G(x,y)| \leq \dfrac{c_3}{{[d_g(x,y)]}^n}. $$

Pour la sphere unite de dimension 2, on peut prendre la fonction suivante:

$$ G(x,y) = \dfrac{\log 2}{4\pi}-\dfrac{\log (1-\cos d(x,y))}{4\pi}. $$

\underbar {\it Questions (Problemes ouverts):}

\bigskip

Sur un domaine $  \Omega $ de $ {\mathbb R}^2 $, l'\'equation se r\'eduit  :

$$ -\Delta u=V(x)e^u, \qquad (E') $$

On suppose:

$$ 0 \leq V \leq b $$

1- A quelle conditions (C) sur $ u  $ et $ V $ peut on avoir

$$  \sup_{\Omega} u + \inf_K u  \geq c=c(b, (C), K, \Omega) ? $$

2- Peut on avoir l'inegalite precedente si on suppose la condition de Brezis-Merle:

$$ \int_{\Omega}  e^u dx  \leq  C_1 ? $$

Les conditions $ (C) $ se reduisent a la precedente.

\bigskip

\chapter{Le cas des dimesions $ n \geq 3 $: In\'egalit\'es du type  $  \sup u \times  \inf u  \leq c $. }

%\underbar {\bf II) Cas de la dimension $ n \geq 3 $:} 

\bigskip

On consid\`ere ici deux types d'\'equations provenant du changement de m\'etrique conforme:

L'equation de Yamabe:

$$ -\Delta_g u+R_gu=n(n-2) u^{N-1},\,\, u >0, \,\, {\rm et } \,\, N=\dfrac{n+2}{n-2}. \qquad (E_1) $$

Et plus g\'en\'eralement l'\'equation de la courbure scalaire prescrite,

$$ -\Delta_g u+R_gu=Vu^{N-1},\,\, u >0, \,\, {\rm et } \,\, N=\dfrac{n+2}{n-2}. \qquad (E_2) $$

avec $ R_g $ la courbure scalaire et $ V $ la courbure scalaire prescrite. On peut remplacer $ R_g $ par une fonction $ h $ ou $ N $ par $ p <N $ ces nouvelles \'equations sont des variantes des premi\`eres.

\bigskip

\section{Les R\'esultats de Aubin-Schoen.} 

~~\\
\bigskip

\begin{Theorem}{\it (Aubin-Trudinger)}. Le probl\`eme de Yamabe a une solution quand l'invariant de Yamabe est n\'egatif.
\end{Theorem}
\bigskip

\begin{Theorem}{\it (Aubin-Schoen)}. Le probl\`eme de Yamab\'e poss\`ede une solution positive lorsque l'invariant de Yamabe est positif.
\end{Theorem}

Voir [2, 3, 44, 53]

\bigskip

\underbar {\bf Sur le Th\'eor\`eme de la masse positive }

\bigskip

Une vari\'et\'e Riemannienne $ (M,g) $ de dimension $ n $ et $ C^{\infty} $, est dite asymptotiquement plate d'ordre $ \tau > 0 $, s'il existe un compact $ K \subset M $ tel que $ M-K $ est diff\'eomorphe \`a $ {\mathbb R}^n - B_0 $ ( $ B_0 $ est la boule unit\'e de $ {\mathbb R}^n $ de centre $ 0 $ ), les composantes de la m\'etrique $ g $ sont:

$$ g_{ij} = \delta_{ij} + O({\rho}^{-\tau}), \,\, \partial_k g_{ij} = O({\rho}^{-\tau -1}), \,\, \partial_{kl} g_{ij} = O({\rho}^{-\tau -2}). $$

La masse $ m(g) $ d'une vari\'et\'e asymptotiquement plate $ (M,g) $ de dimension n, est la limite, si elle existe de la quantit\'e suivante:

$$ {\omega_{n-1}}^{-1} \int_{{\mathbb S}_{n-1}(\rho)} \sqrt {|g(\rho, \theta)|} g^{ij} (\partial_i g_{\rho j} - \partial_{\rho} g_{ij})(\rho,\theta) d\tau(\theta) , $$

avec, $ \rho \to + \infty $ et $ d\tau $ l\'el\'ement d'aire de $ {\mathbb S}_{n-1}(\rho) $.

\bigskip

\begin{Theorem}{\it (Schoen-Yau) $ (n=3) $([2])}. Si $ (M,g) $ est une vari\'et\'e asymptotiquement plate de dimension 3 et d'ordre $ \tau > (n-2)/2 $ avec une courbure scalaire positive ou nulle et dans $ L^1(M) $, alors, $ m(g) \geq 0 $ et $ m(g) = 0 $ si et seulement si $ (M,g) $ est isom\'etrique \`a l'espace euclidien.
\end{Theorem}

\bigskip

\underbar{\it Cons\'equence du Th\'eor\`eme de la masse positive en dimension 3}

\bigskip

Sur une vari\'et\'e Riemannienne compacte de dimension 3 de courbure scalaire $ R_g $, on consid\`ere $ k $ une fonction r\'eguli\`ere qu'on peut supposer, sans nuire \`a la g\'en\'eralit\'e, strictement positive. On suppose que la premi\`ere valeur propre $ \lambda_1 $ de $ -\Delta_g + k $ est positive. Pour $ y \in M $ et $ G_y $ la fonction de Green de $ -\Delta_g + k $ de p\^ole $ y $ :

$$ (-\Delta_g + k)(G_y) =\delta_y . $$

Alors, il existe un r\'eel $ A_{k,g}(y) $ tel que:

$$ G_y(x) = \dfrac{1}{3\omega_3 |x|} + A_{k,g}(y) + O (|x|^{\alpha}), \,\,\, |x| \to 0. $$

\bigskip

\begin{Theorem}{\it (Schoen-Yau)}. On a:

$$ A_{k,g} \geq A_{R_g/8,g} >0. $$
\end{Theorem}
\bigskip

\underbar { \bf L'Obstruction de Kazdan-Warner}

\bigskip

\begin{Theorem} {\it ( Kazdan-Warner)}. On note $ \xi $ un \'el\'ement de l'espace vectoriel correspondant \`a la prem\`ere valeure propre, alors:
$$ {\rm Pour } \,\, {\mathbb S}_n, \,\, n \geq 3, \,\, \int_{{\mathbb S}_n} \nabla^i \xi \nabla_i R u^{(n+2)/(n-2)} dV_g = 0. $$

o\`u $ u $ est solution de l'\'equation de la courbure scalaire prescrite avec $ R $ comme courbure scalaire, sur la sph\`ere unit\'e $ {\mathbb S}_n $, $ n\geq 2 $.
\end{Theorem}
\bigskip

\underbar {\bf Obstruction du type Kazdan-Warner: g\'en\'eralisation}

\bigskip

\begin{Theorem}{\it (Bourguignon-Ezin)}. Pour un champ de vecteurs conforme $ X $ sur une vari\'et\'e Riemannienne compacte $ (M,g) $, on a:

$$ \int_M X(R) dV_g =0 .$$
\end{Theorem}
\underbar{\bf L'Obstruction de Pohozaev}

Sur un ouvert born\'e r\'egulier $ \Omega $ de $ {\mathbb R}^n, n\geq 3 $, on consid\`ere l'\'equation suivante:

$$ -\Delta u= u^{(n+2)/(n-2)}, \,\, u >0, \,\, u=0 \,\, {\rm sur} \,\, \partial \Omega. $$

\begin{Theorem}{\it (Identit\'e de Pohozaev)}. Les solutions $ u > 0 $ du probl\`eme pr\'ec\'edent, v\'erifient:

$$ (1-n/2) \int_{\Omega} u^{2n/(n-2)} dx + n \int_{\Omega} u^{2n/(n-2)} dx = (1/2) \int_{\partial \Omega} \partial_{\nu} (||x||^2/2)(\partial_{\nu} u)^2 d\sigma. $$
\end{Theorem}
\begin{Theorem}{\it (Pohozaev)}. Si l'ouvert $ \Omega $ est \'etoil\'e, alors, le probl\`eme pr\'ec\'edent ne poss\`ede pas de solutions.
\end{Theorem}

\begin{Theorem}{\it (Kazdan-Warner)}. Le probl\`eme pr\'ec\'edent a au moins une solution pour des couronnes.
\end{Theorem}

\begin{Theorem}{\it (Coron, Bahri-Coron)}. Le probl\`eme pr\'ec\'edent a au moins une solution pour des domaines $ \Omega $ \`a trous, en particulier pour des couronnes.
\end{Theorem}

\underbar{Esquisse de la preuve:}

Soit $ b_k = inf J(u) $, apres blow-up et le theoreme de classification de Gidas-Cafferelli-Spruck (obtenu a partir de la symmetrisation de Gidas-Ni-Nirenberg), $b_k=k S^{k/2} $. Le blow-up fait que chaque fonction $ u $ est proche (procede de Struwe) des fonctions (concentrations) ou blow-ups $ \delta_j $ qui ressemblent aux masse de Dirac car ce sont les fonctions cocentrations. 
Donc $ \Sigma_+ $ peut etre considere comme parametr\'e par ces fonctions. Donc les ensembles de niveaux $ W_k $ sont comme $ \tilde \Omega^k \times \Delta_{k-1} $ les simplexes et $ \tilde \Omega $ est une variete compacte obtenue a partir de $ \Omega $ par la condition $ H_d(\Omega) \not =0 $ ($ \tilde \Omega $ est de dimension $ d $ et de classe fondamentale non nulle). Comme le bord de $ W_k $ c'est $ W_{k-1} $, on peut consider que $ W_{k-1} $ comme $ \tilde \Omega^k \times \partial \Delta_{k-1} $, donc:

$$ H_{kd+k-1}(W_k,W_{k-1})=H_{kd+k-1}(\tilde \Omega^k\times \Delta_{k-1}, \tilde \Omega^k \times \partial \Delta_{k-1})\equiv H_d(\tilde \Omega)\otimes H_{k-1}(\Delta_{k-1}, \partial \Delta_{k-1}).$$

Ce que font Bahri et Coron, c'est qu'ils prouve que $ (W_k,W_{k-1}) $ est retract\'e par deformation d'un ensemble intermediaire $ (F, W_{k-1}) $, ceci est visible car $ F\equiv J^{(k+1/2)S^{k/2}} $ qui est retracte par defomation de $ W_k $ par le procede de Struwe.($ u_k \to u_0+ u_1+...u_r $, $ u_j, j \geq 1 $ sont les fonctions concentrations, leur masse vaut $ S^{k/2} $ comme la masse totale ne doit pas depasser $ b_{k+1} $ car on est dans $ W_k $, on a $ r\leq k $. Si $ u_0\equiv 0 $, la masse totale doit etre egale a $ b_k $ ce n'est pas possible car, on est dans l'ensemble de niveau, strictement superieure a $ b_k $, dont $ u_0 \not = 0 $, donc d'apres Struwe, notre suite est de Palais-Smale et d'apres le lemme de deformation $ F $ est "presque " deformation par retract de $ W_k $. (voir le Kavian), (noter que l'hypothese de l'absurde, qu'il n'y a pas de solution, fait qu'il y aun flot, car $ \nabla J $, n'a pas de point critique et dont, on peut parler de deformation par retract, d'apres l'existence du flot.

Maintenant, Bahri et Coron, ont construit, l'ensemble $ F $ de telle maniere que l'homologie de $ (F, W_{k-1}) $  soit nulle, pour $ p=kd+k-1 \not = k-1 $, car $ F $ est dans la parmmetrisation $ V(n,\epsilon) $ de meme pour $ W_{k-1} $.

Donc, Bahri et Coron, en supposant qu'il n'y a pas de solution de cette equation non-lineaire, prouve que:

$$ H_{kd+k-1}(W_k,W_{k-1})=H_{kd+k-1}(F,W_{k-1})=0, $$

et,

$$ H_{kd+k-1}(W_k, W_{k-1})=H_d(\tilde \Omega)\otimes H_{k-1}(\Delta_{k-1}, \partial \Delta_{k-1}) \not =0 $$

si, il y a au moins un $ d $ tel que $ H_d(\Omega) \not = 0 $ qui implique qu'il existe une classe fondamentale non nulle pour $ H_d(\tilde \Omega) $.

\bigskip

\begin{Theorem}{\it (Br\'ezis-Nirenberg)}. Si on perturbe l'\'equation pr\'ec\'edente par le terme lin\'eaire $ \lambda u $ ( pour des $ \lambda >0 $ particuliers), alors le nouveau probl\`eme poss\`ede au moins une solution.
\end{Theorem}
\bigskip

Sur la sph\`ere unit\'e de dimension 3, on consid\`ere l'\'equation suivante:

$$ -8\Delta u + 6u= Ru^5, \,\,\, u > 0. $$

Par une methode de "blow-up "on a:

\begin{Theorem}{\it (Chang-Gursky-Yang)}. Si $ 0 < m \leq R \leq M $, alors, il existe une constante positive $ c=c[m, M, ||R||_{C^2({\mathbb S}_3)}] $, telle que:

$$ \int_{{\mathbb S}_3} |\nabla u|^2 \leq c. $$

\end{Theorem}

L'obstruction de Kazdan-Warner permet d'avoir

\begin{Theorem}{\it (Chang-Gursky-Yang)([22])}. Si $ K $ v\'erifie, $ 0 < m \leq R \leq M $, $ \Delta R \not = 0 $, alors:

$$ 0< 1/c \leq |u|\leq c=c(m, M,|R|_{C^2}). $$
\end{Theorem}

Rappelons le but des estimations a priori de Chang-Gursky-Yang est de prouver que l'equation de la courbure scalaire prescrite a des solutions lorsque les courbures prescrites sont proches des constantes. (Par la methode de continuite). (Probleme de Nirenberg). De plus, ils obtiennent une formule d'indice pour certaines courbures prescrites comme fonctions de Morse.

\underbar{\bf La M\'ethode "moving-plane" et applications}

\bigskip

La m\'ethode "moving-plane" consiste \`a rechercher, si possible, les points de sym\'etrie pour des E.D.P d\'efinies sur des domaines ayant des axes de sym\`etries, puis, de caract\'eriser ces solutions. On part de "l'infini", un point tr\'es loin, puis, on consid\`ere la fonction et son sym\'etris\'ee par rapport au plan contenant ce point, puis on ram\`ene, le plan jusqu'\'a l'annulation de la diff\'erence entre cette fonction et son sym\'etris\'ee, si c'est le cas, le plan limite est le plan de sym\'etrie.

\bigskip

\underbar{\it Un exemple}

\bigskip

Sur la boule unit\'e $ B $ de $ {\mathbb R}^n $, $ n \geq 3 $, on consid\`ere le probl\`eme suivant:

$$ -\Delta u = u^{(n+2)/(n-2) - \epsilon}, \,\, u > 0 \,\, u = 0 \,\, {\rm sur} \,\, \partial B. $$

O\`u $ \epsilon >0 $ est tr\'es petit.

\bigskip

\begin{Theorem}{\it (Gidas-Ni-Nirenberg)}. La solution $ u $ du probl\`eme pr\'ec\'edent est radiale et strictement d\'ecroissante.

\end{Theorem}
Voir [31]

En utilisant la transformee de Kelvin et le principe du maximum, on classifie les solutions de l'equation :

\begin{Theorem}{\it (Caffarelli-Gidas-Spruck)}. Les solutions de:

$$ -\Delta u=u^{(n+2)/(n-2)}, \, u(0)=1, u >0 \, {\rm dans } \,\, {\mathbb R}^n \\  $$
 
Sont:

$$   u(x) = (1+\gamma |x|^2)^{(2-n)/2}, \,\,\gamma >0 $$

\end{Theorem}

La preuve de Chen-Li (voir [23]) est plus courte et elle utilise les memes arguments que ceux de Caffarelli-Gidas-Spruck, c'est essentielement une transformation de Kelvin pour avoir un comportement asymptotique  et aussi l'utilisation du principe du minimum pour des fonctions regulieres singulieres en un point puis l'argument moving-plane bas\'e sur le principe du maximum.

Voir [21, 23]

\bigskip

\underbar{\it Applications de la m\'ethode "moving-plane"}

\bigskip

Si on remplace $ B $ par un ouvert r\'egulier quelconque, not\'e $ \Omega $, en maintenant la m\^eme \'equation, les solutions ne sont pas forc\'ement radiales.

\bigskip

\begin{Theorem}{\it (Han)}. Il existe un voisinage du bord $ \omega $ ne d\'ependant que de la g\'eom\`etrie du domaine $ \Omega $ et de la dimension $ n $, ainsi qu'une constante $ c=c(\Omega,n) >0 $ telle que:

$$ ||u||_{L^{\infty}(\omega)} \leq c. $$

O\`u $ u $ est la solution du probl\`eme pr\'ec\'edent.
\end{Theorem}

Voir [33]

\bigskip

Lorsque $ M= {\mathbb S}_n $, $  n \geq 3 $, T. Aubin a prouv\'e que:

\begin{Theorem}{\it (Aubin)}. Les solutions de l'\'equation de Yamabe sur la sph\`ere unit\'e de dimension $ n $, $ {\mathbb S}_n $, sont donn\'ees par la formule suivante:

$$ \phi_{\beta,P}(Q)=\left [ (\beta^2-1)/(\beta-\cos r)^2 \right ]^{(n-2)/4}, \,\, \beta >1, \,\, r=d(P,Q), \,\, P,Q \in {\mathbb S}_n. $$
\end{Theorem}
\begin{Theorem} On a:
$$  \sup_{{\mathbb S}_n} \phi_{\beta,P}(Q)\times \inf_{{\mathbb S}_n} \phi_{\beta,P}(Q) = 1. $$

\end{Theorem}

Quant \`a  Schoen-Pacard-Mazzeo-Korevaar, ont prouv\'e dans le cas, $ M= \Omega $ un ouvert de $ {\mathbb R}^n $:
\begin{Theorem} 

$$ \sup_K  u  \times \inf_{\Omega} u \leq c(n, K, \Omega) $$

\end{Theorem}

\begin{Theorem}{\it (Schoen)}. Pour une vari\'et\'e Riemannienne compacte sans bord et conform\'ement plate, l'ensemble des solutions de l'\'equation de Yamabe est compact.
\end{Theorem}

Voir [2-45]

\bigskip

\section{R\'esultat de Chen-Lin et YY.Li.}

~~\\

YY.Li prouva dans le cas de la sph\`ere que:

\begin{Theorem} 

 si, on suppose les courbures $ V $ v\'erifient:

$$ 0 < a \leq V(x) \leq b < + \infty $$

$$ ||\nabla^{\alpha} V|| \leq C_{\alpha} ||\nabla V||^{\beta(\alpha)}, \,\,\, {\rm avec} \,\,\, \alpha \leq n-2. $$

alors,

$$  ||u||_{H^1} \leq C(a, b, C_{{\alpha}}, n) $$

et,

$$ \sup_{{\mathbb S}_n} u \times \inf_{{\mathbb S}_n} u \leq c= c(a, b, C_{{\alpha}}) $$
\end{Theorem} 

Li utilisa la notion de blow-up isol\'e et blow-up isol\'es simples ainsi qu'un proc\'ed\'e d'iteraton de ces notions (blow-up de blow-up), des estimations autour des blow-up isol\'es simples et une formule de Pohozaev pour les fonctions avec singularit\'es

La methode suivante est celle utilisee par YY.Li ainsi que YY.Li et M. Zhu pour prouver la compacite des solutions de l'equation de Yamabe. C'est aussi, peut etre la methode Schoen pour prouver la compacite des solutions de l'equation de Yamabe dans le cas conformement plat, notons que les inegalites de type Harnack permettent d'aboutir a ce resultat toujours en utilisant le theoreme de la masse positive.

Rappelons que le but de YY. Li est de resoudre les Probleme de Nirenberg sur la sphere par la methode de continuite en prouvant, grace aux estimations a priori que les solutions existent pour des courbures scalaires proches des constantes. Ils prouvent  que les suites  de solutions sont bornees  et qu'il n'y a pas de blow-up isoles simples. Il existe des solutions pour le probleme de Nirenberg pour des courbures prescrites verifiant des conditions de platitudes.

D'autre part il donne une formule d'indice  pour des courbures prescrites comme fonctions de Morse.

Nous allons donner la deinition d'un point bloow-up isol\'e et point blow-up isol\'e simple.

 \begin{Definition} On dit qu'un point $ x_0 $ est un point blow-up isol\'e de la suite de fonctions $ (u_k) $  definies sur un domaine $ \Omega $, s'il existe une suite de points $ (x_k)_k $ et un r\'eel $ r >0 $ telle que :
 
 a) $ x_k \to x_0 $, $ x_k $ est maximum local de $ u_k $
 
 b) $ u_k(x_k) \to + \infty $. 
 
 c) $ 0 < r < d(x_0, \partial \Omega) $ et il existe $ c >0 $ tel que:
 
$$ u_k(x) \leq c |x-x_k|^{-\frac{n-2}{2}} \,\,\, y \in B_r(x_k). $$

 \end{Definition} 

Soit $ x_0 $ un point isol\'e de $ (u_k)_k $, on definit:

$$ \bar u_k(s)=\dfrac{1}{|\partial B_s|} \int_{\partial B_s(x_k)} u_k, \,\,\,  s>0, $$

et,

$$ \bar w_k(s)=r^{\frac{n-2}{2}} \bar u_k(s), \,\,\,  s>0. $$

\begin{Definition}On dit qu'un point $ x_0 $ est un point blow-up isol\'e simple de la suite de fonctions $ (u_k) $  definies sur un domaine $ \Omega $, s'il est isol\'e et pour un $ \rho >0 $ independant de $ k $, $ \bar w_k $ a precisement un point critique dans $ (0, \rho) $ pour $ k $ assez grand.

 \end{Definition} 

1-Le procede de selection de Schoen:

Soit $ V $ $ C^1({\mathbb S}_n) $ et $ v $ une fonction solution de:

$$ -\Delta_g v-\dfrac{(n-2)}{4(n-1)}R_gv=\dfrac{(n-2)}{4(n-1)}Vv^{(n+2)/(n-2)} $$

avec $ R_g $ la courbure sclalire de la sphere.

On suppose:

$$ 0 < a \leq V \leq b < + \infty $$

$$ ||\nabla V||_{\infty} \leq A. $$

Alors, pour tout $ 1> \epsilon >0 $ et $ R >1 $, il existe deux constantes positives $ C_0=C_0( \epsilon, R, n, a, b, A) $ et $ C_1=C_1( \epsilon, R, n, a, b, A)>1 $ tels que pour $ v $ solution de l'equation precedente avec,

$$  \max_{{\mathbb S}_n} v >C_0, $$

il existe un entier $ k=k(v) $ et un ensemble de points,

$$ S(v)=\{ q_1,\ldots, q_k \} \subset {\mathbb S}_n $$

tels que:

a)-Les $ q_j $ sont des maximums locaux de $ v $ pour $ 1\leq j \leq k $, et dans un systme de coordonnees geodesiques normales:

$$||v(0)^{-1}v(v(0)^{-(n-2)/2}y)-\delta_j (y)||_{C^2(B_{2R(0))}} < \epsilon, $$

avec,

$$ \{B_{Rv(q_j)^{-(n-2)/2}}(q_j)  \}_{1 \leq j \leq k } \,\,\, {\rm disjointes}, $$

et, 

$$ \delta_j (y)=(1+k_j|y|^2)^{(2-n)/2}, $$

est l'unique solution de:

$$ \begin{cases}

-\Delta \delta_j=\dfrac{(n-2)}{4(n-1)}K(q_j)\delta_j^{(n+2)/(n-2)}, \,\, {\rm dans } \,\, {\mathbb R}^n,\\
             \delta_j>0\,\, {\rm dans } \,\, {\mathbb R}^n,\\
              \delta_j(0)=1,
              k_j=\dfrac{1}{4n(n-1)}K(q_j).                      
             
 \end{cases}  $$

b)-"blow-up" isol\'es

$$ v(q) \leq C_1 (d(q, S(v)))^{-(n-2)/2}. \,\,{\rm pour \,\, tout } \,\, q \in {\mathbb S}_n. $$

3-Les blow-up isol\'es sont des blow-up isol\'es simples et sont en nombre fini, par compacite de la sphere. Pour prouver que les blow-up isoles sont isoles simples, l'auteur procede par l'absurde. Apres avoir fait un "blow-up", il existe des fonctions $ w_i $ telles que:

$$ w_i(0) \to + \infty $$

$$ w_i(0)w_i \to a |y|^{2-n}+b=h, \,\, a, b >0 \,\, {\rm dans} \,\, C^2({\mathbb R}^n-\{0 \}), $$

et,

$ 0 $ est un blow-up isole de $ w_i $

Grace aux estimations a priori decoulant du fait que $ 0 $ est blow-up isole simple (en raisonnant par l'absurde, par blow-up, les nouveaux points deviennent simples isoles), il prouve que:

$$  \limsup_i (w_i(0)^2 \times B(\sigma, w_i, \nabla w_i)) \geq 0,  $$

Mais la formule de Pohozaev donne:

$$  \limsup_i (w_i(0)^2 \times B(\sigma, w_i, \nabla w_i) )=B(\sigma, h, \nabla h) <0.  $$

Contradiction. Ici, $ B(\sigma, u, \nabla u) $ est la quantite ne contenant que $ u $ et ses derivees dans la deuxieme partie de la formule de Pohozaev. 

4- Par le meme raisonnement, il prouve qu'il n'y a qu'un nombre fini de blow-up isoles et la distance entre ces blow-ups est minoree par une constante positive $ \delta=\delta(\epsilon, R, n, a, b, A) >0 $. En, effet, en supposant que la distance entre les points tend vers $ 0 $, on exhibe une suite de fonctions notee $ w_i $ et telle que:

$$ w_i(0) \to + \infty $$

$$ w_i(0)w_i \to a_1|y|^{2-n}+a_2 |y-q|^{2-n}+ b(y)=h, \,\, a_1, a_2, b >0 \,\, {\rm dans} \,\, C^2({\mathbb R}^n-(S-\{0, q \})), $$

avec, $ b $ harmonique dans $ {\mathbb R}^n-\cup S $ et $ S $ un ensemble au plus denombrable.

et,

$ 0 $ est un blow-up isole de $ w_i $

Grace aux estimations a priori decoulant du fait que $ 0 $ est blow-up isole simple (en raisonnant par l'absurde, par blow-up, les nouveaux points deviennent simples isoles) il prouve que:

$$  \limsup_i (w_i(0)^2 \times B(\sigma, w_i, \nabla w_i)) \geq 0,  $$

Mais la formule de Pohozaev donne:

$$  \limsup_i (w_i(0)^2 \times B(\sigma, w_i, \nabla w_i) )=B(\sigma, h, \nabla h) <0.  $$

Contradiction.

5- Contradiction grace aux estimations uniformes et la formule de Pohozaev.
  
  \bigskip  
  
 \underbar {Remarques importantes:} 
 
 \bigskip
  
 1) Dans l'identite de Pohozaev, l'operateur conforme implique l'inexistence de termes du type:
 
 $$  \int_{B_{\sigma}(0)} (x.\nabla R_g+R_g)u^2 dx,  \,\,\, {\rm et } \,\,\,   \int_{\partial B_{\sigma}(0)} R_gu^2 d\sigma     $$ 
  
   ( $ R_g $, la courbure scalaire) qui compliqueraient les estimations et pousserait a regarder le degre d'annulation du tenseur de Weyl, dans le cas ou la variete n'est pas la sphere.
  
  \bigskip  
  
  2) Les conditions de platitudes, permettent, grace a la formule de Taylor, d'augmenter la presence du terme $ u_i(0)^{-2/(n-2)} $ autant de fois , que la derivation, ce qui permet de faire apparaitre le terme $ u_i(0)^{-2}   $ et pouvoir appliquer correctement la formule de Pohozaev.

Voir [36]
 
\bigskip
 \underbar{\bf 2.Cas des vari\'et\'es Riemanniennes quelconques}

\bigskip

\underbar{\bf Une cons\'equence du th\'eor\`eme de la masse positive}

\bigskip

\underbar{ \it Compacit\'e des solutions de l'\'equation de Yamabe}

\bigskip

\begin{Theorem}{\it (Li-Zhu)}. Sur toute vari\'et\'e Riemannienne compacte de dimension 3, non conform\'ement diff\'eomorphe \'a la sph\`ere unit\'e de dimension 3, et pour toute fonction positve $ K \in C^2(M) $, $ K \geq a > 0 $, il existe une constante strictement positive
$ C = C(a, ||K||_{C^2(M)}, M, g) $ telle que:

$$  0 < \dfrac{1}{C} \leq u \leq C, $$

o\`u $ u $ est solution de l'\'equation de la courbure scalaire prescrite  relativement \`a $ K $ ($ K $ est la courbure scalaire prescrite).
\end{Theorem}
\bigskip

Sur une vari\'et\'e Riemannienne compacte de dimension 3 et $ k $ une fonction r\'eguli\`ere, on consid\`ere l'\'equation suivante:

$$ -\Delta_g u + k u = u^5, \,\,\, u >0. $$

\begin{Theorem}{\it (Li-Zhu)}. Il existe une constante positive $ c = c( ||k||_{C^1(M)}, M, g) $ telle que:

$$ || u ||_{H^1(M)}  \leq c, $$

avec $ u $ solution de l'\'equation pr\'ec\'edente.
\end{Theorem}
\bigskip

Sur une vari\'et\'e Riemannienne quelconque $ (M,g) $ compacte, de dimension $ n\geq 3 $ et de courbure scalaire $ R_g $, on consid\`ere l'\'equation de Yamabe:

$$ -\Delta_g u + R_g u = u^{(n+2)/(n-2)}, \,\,\, u >0 . $$

\begin{Theorem}{\it (Druet)}. Les solutions de l'\'equation de Yamabe forment un ensemble compact dans $ C^0(M) $, en diemensions 4 et 5.
\end{Theorem}
\bigskip

\begin{Theorem}{\it (Marques)}. Les solutions de l'\'equation de Yamabe forment un ensemble compact dans $ C^0(M) $, en dimensions $ 4, 5, 6, 7 $.
\end{Theorem}
\bigskip

\begin{Theorem}{\it (Li-Zhang)}. Les solutions de l'\'equation de Yamabe forment un ensemble compact dans $ C^0(M) $, en dimensions $ 4, 5, 6, 7,8 $.
\end{Theorem}
\bigskip

II)- En utilisant la m\'ethode moving-plane, Chen et Lin ont prouv\'e le m\^eme type d'estimations sur un ouvert $ \Omega $ de $ {\mathbb R}^n $:

\begin{Theorem}  Si,

$$ 0 < a \leq V(x) \leq b < + \infty $$

$$ ||\nabla^{\alpha} V|| \leq C_{\alpha} ||\nabla V||^{\beta(\alpha)}, \,\,\, {\rm avec} \,\,\, \alpha \leq n-2. $$

alors,

$$ \sup_K u \times \inf_{{\Omega}} u \leq c= c(a, b, C_{{\alpha}}, n,  K, {\Omega}) $$
\end{Theorem}

Voir [25]

Esquisse de la preuve:

\bigskip

La methode utiis\'ee est celle des d\'eplacements de plans. Le but est de prouver que dans une voisinage d'un certain point critique la deriv\'ee d'une certaine fonction est strictement positive, ce qui contredirait le fait qu'il existe des  maximas locaux  au voisinage du point consid\'er\'e au d\'epart. Nous donnons ici une esquisse de la preuve. La condition de platitude est utilis\'ee pour estimer les deriv\'ees successives au point blow-up en fonction " de la fonction blow-up" et avoir des estim\'ees du type:

$$  |\nabla V(x_k)| \leq C u_i(x_k)^{-2/(n-2)} $$

par exemple, cette estimation est \'a l'ordre 1. Et ces dernieres estimations impliquent la positivit\'e de la fonction auxiliere, qui corrige la perturbation dans l'\'equation. Ici, on a note $ (x_k)_k $ la suite blow-up.

\bigskip

1-Technique "blow-up":

Etape 1: On veut etablir une inegalite du type:

$$ \epsilon^{n-2} \max_{B(0, \epsilon)} u \times \min_{B(0, 4 \epsilon) } u \leq c=c(a, b, A, M, g). $$

Pour cela, on raisonne par l'absurde et on suppose que:

$$  \max_{B(0, \epsilon_k)} u_k \times \min_{B(0, 4 \epsilon_k) } u_k \geq k {\epsilon_k }^{2-n}.$$

Etape 2:

Notre hypothese de l'absurde se trouve mise en defaut si on arrive  prouver que:

$$ \forall \,\, \epsilon >0, \,\, \exists \,\, \eta = \eta(\epsilon)>0, \,\, \min_{|y|\leq r} v_k(y) \leq (1+\epsilon) U(r) \,\,\,  \forall \,\, 0\leq r \leq \eta l_k, \,\,\, l_k=\dfrac {1}{2}R_k M_k^{2/(n-2)}$$  

On raisonne alors par l'absurde et on suppose qu'il existe $ \epsilon >0 $ et une suite de r\'eels $ r_k \to + \infty $ tels que:

$$  \min_{|y|\leq r_k} v_k(y) \geq (1+\epsilon) U(r_k),   \,\, 0\leq r \leq \eta l_k, \,\,\, l_k=\dfrac {1}{2}R_k M_k^{2/(n-2)}$$

\bigskip

2-Methode "moving-plane":

Pour $ x\in {\mathbb R}^2 $ et $ \lambda \in {\mathbb R}^*_- $, on pose,

$$ v_k^1(y):= \dfrac {1}{|y|^{n-2}}v_k \left (e+\dfrac {y}{|y|^2}\right  ), \,\,\, \bar U_0(y)=\dfrac {1}{|y|^{n-2}}U \left (e+\dfrac {y}{|y|^2}\right  ) $$

c'est la transformation de Kelvin de $ v_k $ pour la boule unite. On pose$ e=(1,0,...0) $.

On cherche a utiliser le principe du maximum et le lemme de Hopf, la methode moving-plane, pour prouver que dans un voisinage du point $ e^*=(-1/2,0,...0) $, la derivee est strictement positive, ce qui contredirait le fait que, puisqu'il y a convergence uniforme, $ v_k^1 $ n'aurait pas de maximum local au voisinage de ce point , ce qui est absurde. Les deux auteurs utilisent des estimees asymptotiques et comparent la fonction obtenue par le procede de symetrisation a une fonction auxiliere, la fonction auxiliere est obtenue a partir de la fonction de Green, les estimees asymptotiques sont utilies pour prouver la positivite de la fonction auxiliere. 

\bigskip

La comparaison de la fonction obtenue par le procede de symmetrisation a la fonction auxiliere, la positivite de la fonction auxiliere et les conditions aux bord  constituent des conditions suffisantes permettant d'appliquer la methode moving-plane, le principe du maximum combin\'e au Lemme de Hopf.

\bigskip

\underbar {\it L'identit\'e de Pohozaev :}

\bigskip

Soit $ \Omega $ un domaine de $ {\mathbb R}^n, n \geq 3 $ de frontiere reguliere. On considere l'equation suivante:

$$ -\Delta u=f(x, u), \,\, x \in \Omega, $$

avec $ f $ continue en $ x $ et $ u $. On pose:

$$ F(x, u) =\int_0^u f(x,s) ds. $$

Alors:

$$ \int_{\Omega} \left ( nF(x,u)-\dfrac{n-2}{2} uf(x,u)+\sum_{i=1}^n x_i F_{x_i}(x,u) \right )= $$

$$ = \int_{\partial \Omega} \left (  \left ( \sum_{i=1}^n x_i \nu_i (F(x,u)-\dfrac{1}{2}|\nabla u|^2 \right )+ \partial_{\nu} u ( \sum_{i=1}^n x_i u_i)+\dfrac{n-2}{2}u|\partial_{\nu} u| \right ) d\sigma. $$

Dans le cas ou notre equation se reduit a:

$$ -\Delta u=c(n)K(x) |u|^{p-1}u, \,\, x \in \Omega, $$

L'identite de Phozaev devient:

$$ c(n)\int_{B_{\sigma}} \left ( \dfrac{1}{p+1} \left ( \sum_{i=1}^n x_i \partial_i K \right ) |u|^{p+1}+\left ( \dfrac{n}{p+1}-\dfrac{n-2}{2} \right ) K|u|^{p+1} \right )-\dfrac{\sigma c(n)}{p+1}\int_{\partial B_{\sigma}} K |u|^{p+1}= $$

$$ = \int_{\partial B_{\sigma}}  \left ( \dfrac{n-2}{2} u \partial_{\nu} u- \dfrac{\sigma}{2}|\nabla u|^2 +\sigma (\partial_{\nu} u)^2 \right ) d\sigma. $$

On pose alors,

$$ B(\sigma, x, u, \nabla u) = \dfrac{n-2}{2} u \partial_{\nu} u- \dfrac{\sigma}{2}|\nabla u|^2 +\sigma (\partial_{\nu} u)^2. $$

et, 

$$ B(\sigma, u, \nabla u) =\int_{\partial B_{\sigma}} B(\sigma, x, u, \nabla u), $$

 \bigskip

\section{Estimation a priori pour de petites variations de la courbure scalaire prescrite.}

~~\\
\underbar {\it Estimation asymptotique : 1. op\'erateurs d'ordre 2, cas de la boule unit\'e}

\bigskip

Sur la boule unit\'e de $ {\mathbb R}^n, n\geq 3 $, $ B_1 $, on consid\`ere le probl\`eme suivant ( $ \epsilon > 0 $):

$$ -\Delta u = n(n-2) u^{(n+2)/(n-2) - \epsilon} \,\,\, {\rm dans } \,\,\, B_1, $$

$$ u > 0 \,\,\, {\rm dans} \,\,\, B_1 \,\,\, u = 0 \,\,\, {\rm sur} \,\,\, \partial B_1. $$

D'apr\'es le travail de Gidas-Ni-Nirenberg, les solutions $ u_{\epsilon} $ de l'\'equation pr\'ec\'edente sont radiales.

\bigskip

\begin{Theorem}{\it (Atkinson-Peletier)}. Les fonctions $ u_{\epsilon} $ sont telles que:

$$ \lim_{\epsilon \to 0 } \epsilon u_{\epsilon}^2 (0) = \dfrac{32}{\pi}, $$

$$ \lim_{\epsilon \to 0} {\epsilon}^{-1/2} u_{\epsilon}(x) = \dfrac{1}{4} \sqrt (\pi/2) \left (1/|x| -1 \right ). $$
\end{Theorem}

Voir [1]

En dimension 3, on consid\`ere le probl\`eme suivant sur la boule unit\'e $ B_1 $:

$$ -\Delta u - \lambda u = 3 u^{5 - \epsilon} \,\,\, {\rm dans } \,\,\, B_1, $$

$$ u > 0 \,\,\, {\rm dans} \,\,\, B_1 \,\,\, u = 0 \,\,\, {\rm sur} \,\,\, \partial B_1. $$

\begin{Theorem}{\it (Br\'ezis-Peletier)}. La fonction $ u_{\epsilon} $ solution du probl\`eme pr\'ec\'edent est telle que:

$$ (a) \,\,\,  {\rm si } \,\,\, 0 \leq \lambda \leq {\pi }^2/4, \,\,\, {\rm alors } \,\,\, \lim_{\epsilon \to 0} \epsilon  u_{\epsilon}^2(0) =\dfrac{ 32 \sqrt \lambda}{\pi \tan \sqrt \lambda}, $$

$$ (a) \,\,\,  {\rm si } \,\,\, 0 \leq \lambda  <  {\pi }^2/4, \,\,\, {\rm alors \,\, pour \,\, x \not = 0, } \,\,\, \lim_{\epsilon \to 0} \epsilon^{-1/2}u_{\epsilon}(x) = \left ( \dfrac {\pi^3  \tan \sqrt \lambda}{2 \sqrt \lambda} \right )^{1/2} G_{\lambda } ( {\rm x } ), $$

avec, $ G_{\lambda } $ la fonction de Green de l'op\'erateur $ \Delta - \lambda $ avec condition de Dirichlet. 
\end{Theorem}
\bigskip

Notation : $ \dfrac{ \sqrt \lambda}{ \tan \sqrt \lambda} = 1 $ si  $ \lambda = 0 $.

\bigskip

\underbar {\it Estimation asymptotique : 1. op\'erateurs d'ordre 2, cas d'un ouvert r\'egulier quelconque}

\bigskip

Soit $ \Omega $ un ouvert r\'egulier born\'e de $ {\mathbb R}^n, n\geq 3 $. On consid\`ere le probl\`eme suivant ( $ \epsilon > 0 $):

$$ -\Delta u = n(n-2) u^{(n+2)/(n-2) - \epsilon} \,\,\, {\rm dans } \,\,\, \Omega, $$

$$ u > 0 \,\,\, {\rm dans} \,\,\, \Omega \,\,\, u = 0 \,\,\, {\rm sur} \,\,\, \partial \Omega. $$

\begin{Theorem}{\it (Han)}. Pour $ u_{\epsilon} $ solution du probl\`eme pr\'ec\'edent, on suppose que:

$$ \dfrac{ \int_{ \Omega } |\nabla u_{\epsilon} |^2 }{ ||u_{\epsilon} ||_{L^{2n/(n-2) - \epsilon }}^2 } = S_n + o(1) \,\,\, {\rm pour } \,\,\, \epsilon \to 0, $$

avec, $ S_n $ est la meilleure constante de Sobolev dans $ {\mathbb R}^n $:

$$ S_n = \pi n(n-2) \left [ \dfrac{\Gamma (n/2)}{\Gamma(n)} \right ]. $$

Alors, apr\'es passage au sous-suites, pour $ \epsilon \to 0 $ :

\bigskip

Il existe $ x_0 $ de $ \Omega $ tel qu'on ait:

$$ u_{\epsilon} \to 0 \,\,\, {\rm sur} \,\,\, C^1(\Omega - \{ x_0 \} ), $$

$$ |\nabla u_{\epsilon} |^2 \to n (n-2) \left [\dfrac{S_n}{n(n-2)} \right ]^{n/2} \delta_{x_0}, $$

au sens des distributions.

$$ \lim_{\epsilon \to 0 } \epsilon ||u_{\epsilon}||_{L^{\infty}}^2 = 2 \sigma_n^2 \left [\dfrac{n(n-2)}{S_n} \right ]^{n/2} |g|, $$

$$ \dfrac{u_{\epsilon}(x)}{\sqrt \epsilon} \to \left [ \dfrac{n(n-2)}{S_n} \right ]^{n/4} \dfrac { (n-2) G(x,x_0) }{ \sqrt { 2 |g|}}, $$

avec, $ g = g(x,x) $ la valeur de la partie r\'eguli\`ere de la fonction de Green du Laplacien avec condition de Dirichlet, en $ x_0 $. Le r\'eel $ \sigma_n $ est l'aire de la sph\`ere unit\'e de $ {\mathbb R}^n $.
\end{Theorem}
\bigskip

On suppose ici, $ u_{\epsilon} $ solution de:

$$ -\Delta u_{\epsilon} =V_{\epsilon} {u_{\epsilon}}^{N-1-\epsilon},\,\, u_{\epsilon} >0, \,\, {\rm et } \,\, N=\dfrac{n+2}{n-2}. \qquad (E_{\epsilon}) $$

On a:

\begin{Theorem}(Bahoura) Si 

$$ 0 < a \leq V_{\epsilon}(x) \leq b < + \infty $$

$$ ||\nabla V_{\epsilon}||_{\infty} \leq  A,  $$

alors,

$$ {\epsilon}^{(n-2)/2} (\sup_K u_{\epsilon})^{1/4}  \times \inf_{{\Omega}} u_{\epsilon}  \leq c= c(a, b, A, n,  K, {\Omega}) $$

\end{Theorem}

Voir [6]

\begin{Theorem}(Bahoura) Si 

$$ 0 < a \leq V_{\epsilon}(x) \leq b < + \infty $$

$$ ||\nabla V_{\epsilon}||_{\infty} \leq k\epsilon\,\, {\rm avec} \,\,\, k >0, $$

alors,

$$ (\sup_K u_{\epsilon})^{4/5}  \times \inf_{{\Omega}} u_{\epsilon}  \leq c= c(a, b, k, n,  K, {\Omega}) $$

\end{Theorem}

Voir [6]

Esquiisse de la preuve:

\bigskip

1-Technique "blow-up": proc\'ed\'e de Schoen, esimtations elliptiques et theoreme de classification de Caffarelli-Gidas-Spruck.

\bigskip

2-Methode "moving-plane": ou plutot "moving-sphere" c'est essentielement  un passage en coordonnees polaires, application du principe du maximum de Hopf, Lemme de Hopf et proprietes des fonctions classifiees par Caffarelli-Gidas-Spruck.

\bigskip

\section{Cas d'une perturbation nonlin\'eaire de l'\'equation.}

~~\\

Si on suppose $ u $ solution de :

$$ -\Delta u =V u^{N-1} + W u^{\alpha},\,\, u >0, \,\, {\rm et } \,\, N=\dfrac{n+2}{n-2}, \,\, \dfrac{n}{n-2} \leq \alpha < \dfrac{n+2}{n-2}, $$

\begin{Theorem} (Bahoura) Si

$$ 0 < a \leq V(x) \leq b \,\, {\rm et } \,\, 0 < c \leq W(x) \leq d, $$

$$ ||\nabla V||_{\infty} \leq A\,\, {\rm et} \,\, ||\nabla W||_{\infty} \leq B, $$

Alors,

$$ \sup_K u \times \inf_{{\Omega}} u  \leq c= c(a, b, c, d, A, B, n, K, {\Omega}) $$

\end{Theorem}

Voir [6]

\bigskip

\underbar {\it Question :}

\bigskip

Que se passe t il si on prend $ \alpha = 1 $ ?

\bigskip

\section{Le cas Holderien en dimension 3.}

~~\\

Sur un ouvert de $ \Omega $ de $ {\mathbb R}^3 $, les solutions de :

$$ -\Delta u =V u^5 \,\, {\rm et } \,\, u >0 $$

\begin{Theorem} (Bahoura) Si

$$ 0 < a \leq V(x) \leq b \,\, {\rm et } \,\, |V(x)-V(y)| \leq A|x-y|^s \,\, {\rm et} \,\,  1/2 < s \leq 1 $$

verifient,

$$ (\sup_K u)^{2s-1} \times \inf_{{\Omega}} u  \leq c= c(a, b, A, s, K, {\Omega}) $$

\end{Theorem}

Voir [9]

La preuve de ce resultat s'inspire de la technique proposee par YY.Li et L. Zhang.

\bigskip

\bigskip

\underbar {\it Question :}

\bigskip

Peut on etendre ce resultat aux varietes riemanniennes de dimension 3 ?

\bigskip

\section{Les vari\'et\'es Riemanniennes et l'\'equation de Yamabe, le cas g\'en\'eral.}

\subsection{Influence de la courbure scalaire.}

~~\\

Sur une vari\'et\'es Riemanniennes $ (M,g) $ de dimension $ n\geq 3 $ ( pas n\'ec\'essairment compacte), on consid\`ere l'\'equation suivante (de type Yamabe)

$$ -\Delta u-h(x) u=n(n-2)u^{N-1},\,\, N=2n/(n-2). $$

\begin{Theorem}(Bahoura) On note $ S_g $ la courbure scalaire de $ M $. et on suppose pour $ m >0 $ donn\'e,

 $$ 0 < m \leq h+\dfrac{(n-2)}{4(n-1)}S_g \leq 1/m. $$
 
Alors, pour tout compact $ K $ de $ M $, il existe une constante positive $ c=c(K,M,m,n,g) $ telle que:

$$ \sup_K u \times \inf_M u \leq c. $$

\end{Theorem}

Voir [7]

\subsection{Resultat de Li-Zhang en dimension 3 et 4.}

~~\\

YY.Li et L. Zhang, ont prouv\'e que les solutions de l'\'equation de Yamabe en dimension 3 et 4:

\begin{Theorem} Si,

 $$ -\Delta u+\dfrac{(n-2)}{4(n-1)}S_g u=n(n-2)u^{N-1},\,\, N=2n/(n-2), \,\, n=3, 4 $$

alors on a l'estimee a priori suivante,

$$ \sup_K u \times \inf_M u \leq c. $$

\end{Theorem}

Voir [39]

\bigskip

\subsection{Le cas des dimensions 5 et 6.}

~~\\

En dimension 5 et 6, nous montrons que si $ u >0 $ est solution de 

$$ -\Delta u+\dfrac{(n-2)}{4(n-1)}S_gu=n(n-2)u^{N-1},\,\, u >0 \,\, {\rm et} \,\, N=\dfrac{2n}{n-2}, \,\, n=5, 6 $$ 

alors,

 \begin{Theorem} (Bahoura) Si $ n=5 $, alors pour tout compact $ K $ de $ M $, il existe une constante positive $ c=c(K,M,g) $ telle que:

$$ (\sup_K u)^{1/7} \times \inf_M u \leq c. $$

Si $ n=6 $ alors, pour tout $ m>0 $, pour tout compact $ K $ de $ M $, il existe une constante positive $ c=c(K,M,m,g) $ telle que:

$$ \sup_K u \leq c, \,\,\, {\rm si} \, \,\, \min_M u \geq m. $$

\end{Theorem}

Voir [10]

\bigskip

\subsection{L'equation de la courbure scalaire prescrite en dimension 3 et 4.}

~~\\

En dimension 3, nous montrons que si $ u>0 $ est solution de:

$$ -\Delta u +\dfrac{1}{8} S_g u= V u^5, $$

Alors:

\begin{Theorem}(Bahoura). Si

$$ 0 < a \leq V \leq b < +\infty,\,\, {\rm et} \,\, ||\nabla V||_{\infty} \leq A, $$

alors,

pour tout compact $ K $, il existe une constante positive $ c=c(a,b,A, K, M, g) $ telle que:

$$ (\sup_K u)^{1/3} \times \inf_M u \leq c. $$

\end{Theorem}.

\bigskip

En dimension 4, nous montrons que si $ u >0 $ est solution de 

$$ -\Delta u+\dfrac{1}{6}S_gu=Vu^3,\,\, $$ 

\begin{Theorem} (Bahoura) Si

$$ 0 < a \leq V(x) \leq b, \,\,  ||\nabla V||_{\infty} \leq A\,\, {\rm et} \,\, A \to 0, $$

alors, 

pour tout $ m>0 $, pour tout compact $ K $ de $ M $, il existe une constante positive $ c=c(a, b, A, K,M,m,g) $ telle que:

$$ \sup_K u \leq c, \,\,\, {\rm si} \, \,\, \min_M u \geq m. $$

\end{Theorem}

Voir [14] et [6] pour le cas euclidien.

\bigskip

Esquisse de la preuve pour le cas Riemannien:

\smallskip

1)-Elements de geometrie: pour tout point de la variete, il existe un voisinage geodesiquement convexe (Theoreme de Whitehead). L'exponentielle tranforme une boules ouvertes en boules ouvertes, les boules fermees en boules fermees et les spheres en spheres. Tout se passe, via l'application exponentielle, comme dans le cas plat. Puis on applique le procede 'blow-up" de Schoen. Les estimations elliptiques et le Theoreme de classification de Caffarelli-Gidas-Spruck nous permet d'avoir le compotement des fonctions "blow-up". On tire des proprietes ineterssantes pour les fonctions blow-up.

2)-Coordonnes geodesiques polaires et Laplacien en ces coordonnes. Puis on applique la methode "moving-sphere" qui est essentielement, le principe du maximum de Hopf, Le lemme de Hopf.

3)-Cette methode qu'on a etendu au cas Riemannien a ete utilisee par Brezis-Li-Shafrir en dimension 2 et par Korevaar-Mazzeo-Pacard-Schoen dans le cas plat en dimension plus grande ou egale a 3. Voir [18] et [34] pour le cas euclidien.

\bigskip

\chapter{Sur la minoration du produit $  \sup \times \inf $ et influence de la fonction de Green.}

\bigskip

\section{Sur une vari\'et\'e compacte.}

~~\\

Soit $ (M,g) $ une vari\'et\'e Riemannienne  compacte de dimension $ n\geq 3 $ et $ \epsilon >0 $. On consid\`ere les solutions de :

\begin{Theorem} (Bahoura) Si,

$$ -\Delta u_{\epsilon}+\epsilon u_{\epsilon}=V_{\epsilon} {u_{\epsilon}}^{N-1}, \,\, N=\dfrac{2n}{n-2}$$

avec,

$$  0 < a \leq V_{\epsilon}(x) \leq b < +\infty, $$ 

alors, 

pour tout $ a,b,m>0 $, il existe une constante positive $ c=c(a,b,m,M,g) $ telle que pour tout $ \epsilon >0 $ et toute solution positive $ u_{\epsilon} $ de:

$$ \max_M u \geq m >0 \,\, \Rightarrow \,\, \epsilon \max_M u_{\epsilon} \times \min_M u_{\epsilon} \geq c. $$
\end{Theorem}

Voir [7]

Dans la preuve de ce resultat, on prouve l'estimation suivante pour  la fonction de Green $ G_{\epsilon} $ de l'operateur  $-\Delta +\epsilon $ et le procede d'iteration de Moser avec les inegalites de Sobolev.

\begin{Theorem} (Bahoura) Pour  $ x_i \to x $ et $ y_i \to y $, il existe une constante positive $ c= c(x, y, g, M) $ telle que:

$$ G_{\epsilon_i} (x_i,y_i) \geq \dfrac{c}{\epsilon_i}, $$

\end{Theorem}

\begin{Theorem} (Hebey-Vaugon). Sur une variete riemannienne compacte sans bord, il existe une constante optimale $ A $ et une autre constante $ B $ telles que:

$$ ||u||_{L^{2^*}}  \leq A ||u||_{H^1(M)}+ B ||u||_{L^2(M)} \,\,  \forall \,\, u \in H^1(M).     $$

\end{Theorem}

Rappelons que sur une variete compacte $ (M, g) $ la construction de la fonction de Green se fait grace au noyau de Green usuel et le produit de convolution. Les estimes de Giraud donne la regularite des parties  formant cette fonction de Green:

$$ G(x,y)=H(x,y)+\sum_{i=1}^k \int_M \Gamma_i(x,z)H(z,y)dv_g(z) + u(x,y), $$

avec,

$$ H(x,y)=\dfrac{\eta}{(n-2)\omega_{n-1} d_g(x,y)^{n-2}}, \,\, {\rm \eta \,\, fonction\,\, cutoff },\,\, $$

et, 

$$ \Gamma_1(x,y) =\Delta_{g,y} H(x,y)-\epsilon H(x,y), $$

et,

$$  \Gamma_{i+1}(x,y) = \int_M \Gamma_i(x,z) \Gamma_1(z,y)dv_g(z), $$

et $ u $ est telle que:

$$ -\Delta_{g,y} u_x +\epsilon u_x = \Gamma_{k+1, x} = \Gamma_{i+1}(x,y) $$

avec, $ k=E(n/2)+1 $, $ E (s) $ la partie entiere de $ s $.

\bigskip

Rappelons les estimees de Giraud:

$$ |\Gamma_{i}(x,y)|\leq   \dfrac{c}{ d_g(x,y)^{n-2i}}\,\, {\rm si},\,\, 2i < n, $$

$$ |\Gamma_{i}(x,y)|\leq   c(1+  |\log d_g(x,y)|) \,\, {\rm si},\,\, 2i = n, $$

$$  |\Gamma_{i}(x,y)|\leq  c \,\, {\rm si},\,\, 2i > n, $$

et la derniere ligne, la fonction est continue.

\bigskip

\section{Sur un ouvert de $ {\mathbb R}^n $ avec condition au bord.}

~~\\

Sur la boule unit\'e $ B $ de $ {\mathbb R}^n $, on consid\`ere l'\'equation suivante:

\begin{Theorem} (Bahoura) Si,

 $$ -\Delta u =V u^{N-1},\,\, u >0, \,\, {\rm et }, u=0 \,\, {\rm sur} \,\, \partial B, \,\, N=\dfrac{n+2}{n-2},  $$

et,

$$ 0 \leq V(x) \leq b, $$

alors, pour tout compact $ K $ de $ B $ on a,

$$  (\sup_B u)^7 \times \inf_K u \geq c=c(b,n,K) > 0. $$

\end{Theorem}

Voir [9]

\bigskip

Esquiisse de la preuve:

1-Technique "blow-up":

\bigskip

2-Nature de la fonction de Green, en dimension n et principe du maximum:

\bigskip

\section{Sur un ouvert de $ {\mathbb R}^n $ sans condition au bord.}

~~\\

Sur la boule unit\'e $ B $ de $ {\mathbb R}^n $, on consid\`ere l'\'equation suivante:

\begin{Theorem} (Bahoura) $$ -\Delta u =V u^{N-1},\,\, u >0, \,\, {\rm et },\,\, N=\dfrac{n+2}{n-2},  $$

et,

$$ 0 <a \leq V(x) \leq b < + \infty , $$

On suppose qu'il existe un r\'eel positif $ m $ et point $ x_0 \in B $, tels que:

$$ u(x_0) \geq m >0 $$

alors, pour tout compact $ K $ de $ B $ on a,

$$  \sup_B u \times \inf_K u \geq c=c(a, b, m, x_0, n,K) > 0. $$

\end{Theorem}

Voir [2] des references supplementaires.

\bigskip

Esquiisse de la preuve:

La preuve est basee sur le procede d'iteration de Moser et une minoration de la fonction de Green sur les compacts.

\bigskip

\section{Application des estimations du type $  \sup \times \inf $.}

\subsection{Estimations d'\'energie et convergence vers la fonction de Green.}

~~\\

\begin{Theorem} (Bahoura)([8]). Sur $ \Omega $ un ouvert born\'e de $ {\mathbb R}^n $, on consid\`ere les solutions du probl\`eme suivant:

$$ -\Delta u_{\epsilon_i}=n(n-2){u_{\epsilon_i}}^{N-1-\epsilon_i},\,\,\, u_{\epsilon_i}>0 \,\, {\rm dans}\,\, \Omega \,\, {\rm et} \,\, u_{\epsilon_i}=0 \,\, {\rm sur} \,\, \partial \Omega \qquad (E), $$

\bigskip

Nous avons,

\bigskip

 Il existe $ c_1=c_1(n, \Omega)>0, c_2=c_2(n, \Omega)>0 $ telles que:

$$ c_2 \leq ||u_{\epsilon_i}||_{H^1(\Omega)} \leq c_1. $$

2) Si $ \Omega $ est \'etoil\'e, alors, il existe une sous-suite $ (u_{\epsilon_j}) $ pour laquelle, il existe un $ m \in {\mathbb N}^* $ et un nombre fini de points de concentrations $ x_1,x_2,...,x_m\in \Omega $ tels que:

$$ i) \lim_{\epsilon_j \to 0} u_{\epsilon_j} = 0 \,\,\, {\rm dans } \,\, C_{loc}^2(\bar \Omega-\{x_1,\ldots x_m\}), $$

$ \forall \,\, k\in \{1,\ldots, m\}, \,\, \exists \,\, (x_{j,k}) $ avec , $ x_{j,k} \to x_k $ et $ u_{\epsilon_j}(x_{j,k})\to+\infty $.

$$ ii) \lim_{\epsilon_j \to 0}  u_{\epsilon_j}^{N-\epsilon_j} =\sum_{i=1}^m \mu_i\delta_{x_i} \,\,\, {\rm avec } \,\, \mu_i\geq \dfrac{\omega_n}{2^n}. $$

Ici la convergence est au sens des distributions.

\bigskip

 iii)  Pour tout compact $ K $ de $ \Omega-\{x_1,\ldots,x_m \} $, il existe une constante positive $ c=c(K,\Omega,n)>0 $ telle que:

$$ \sup_{\Omega} u_{\epsilon_j} \times \sup_K u_{\epsilon_j} \leq c. $$

 iv)  Il existe un voisinage $ \omega $ du bord  $ \partial \Omega $ et une constante positive $ \bar c=\bar c(\omega,\Omega,n) $ tels que:

$$ \sup_{\Omega} u_{\epsilon_j} \times \sup_{\omega} u_{\epsilon_j} \leq \bar c. $$
 
 v)  il existe deux constantes positives, $ \beta_1 $ et $ \beta_2 $ telles que:

$$ \beta_1 \leq \epsilon_j \left ( \sup_{\Omega} u_{\epsilon_j} \right )^2 \leq \beta_2, $$

plus pr\'ecis\'ement, il existe une fonction $ g \in C^2(\partial \Omega) $, telle que,

$$  \epsilon_j \left ( \sup_{\Omega} u_{\epsilon_j} \right )^2 \to \dfrac{ c_n \int_{\partial \Omega} <x|\nu(x)>[ \partial_{\nu} g(\sigma)]^2 d\sigma }{\sum_{k=1}^m \mu_k}. $$

  vi) il existe $ m $ r\'eels positifs $ \gamma_1, \ldots, \gamma_m $, $ \gamma_k\geq n(n-2)\dfrac{\omega_n}{2^n} $, $ k\in\{1,\ldots,m \} $, tels que:
  
  $$ \sup_{\Omega} u_{\epsilon_j} \times u_{\epsilon_j}(x) \to \sum_{k=1}^m \gamma_k G(x_k,x) \,\, {\rm dans } \,\, C_{loc}^2(\bar \Omega-\{x_1,\ldots,x_m\}), $$
 
  o\`u $ G $ est la fonction de Green du laplacien avec condition de Dirichlet. On peut prendre, $ g=\sum_{k=1}^m \gamma_k G(x_k,.) $ dans le v).
 
 \end{Theorem}

\bigskip

\subsection{Equations dont les solutions sont r\'eduites aux constantes.}

~~\\

\underbar {\bf Quelques Propri\'et\'es de l'espace Euclidien.}

\bigskip

Soit $ \Omega  $ un ouvert de $ {\mathbb R}^n $ et $ \lambda_1 $ la premiere valeur propre du laplacien avec condition de Dirichlet. On a le resultat suivant:

\begin{Theorem}{\it(Zhu)}. Pour $ \alpha_i \to \lambda_1 $ et $ n \geq 3 $, les solutions de l'\'equation suivante:

$$ -\Delta u_i + \alpha_i u_i = u_i^{(n+2)/(n-2)}, \,\, {\rm dans}\,\, \Omega \,\, {\rm et} \,\, u_i=0 \,\, {\rm sur} \,\, \partial \Omega$$

sont telles qu'\`a partir d'un certain rang $ i_0 $, $ (u_i)_i \equiv cte $

\end{Theorem}

\begin{Theorem}{\it(Zhu)}. Sur un ouvert convexe $ \Omega  $ de $ {\mathbb R}^3 $, les solutions de l'\'equation suivante:

$$ \begin{cases}
  -\Delta u_i + \epsilon_i u_i = u_i^{(n+2)/(n-2)} \,\, \text{dans} \,\,  \Omega, \\
   \,\, \partial_{\nu} u_i =0                                    \qquad \qquad \qquad \qquad  \text{sur}   \quad \partial \Omega.           
 \end{cases}  $$

sont telles qu'\`a partir d'un certain rang $ i_0 $, $ u_i\equiv cte $

\end{Theorem}

Voir [55] et [56]

On a le resultat suivant qui donne l'existence de branches de solutions non-constantes du probleme de Lin-Ni.

\begin{Theorem}{\it(Adimurthi-Yadava)}. Sur une boule $ \Omega =B_R(0)  $ de $ {\mathbb R}^n $, $ (n=4, 5, 6) $ il existe des solutions radiales non constantes de l'\'equation suivante:

$$ \begin{cases}
  -\Delta u_i + \epsilon_i u_i = u_i^{(n+2)/(n-2)} \,\, \text{dans} \,\,  \Omega, \\
   \,\, \partial_{\nu} u_i =0                                    \qquad \qquad \qquad \qquad  \text{sur}   \quad \partial \Omega.           
 \end{cases}  $$

\end{Theorem}

On a l'unicite de solutions du probleme de Lin-Ni pour des fonctions radiales.

\begin{Theorem}{\it(Adimurthi-Yadava)}. Sur une boule $ \Omega =B_R(0)  $ de $ {\mathbb R}^n $, $ (n \geq 7) $, les solutions radiales de l'equation suivante:

$$ \begin{cases}
  -\Delta u_i + \epsilon_i u_i = u_i^{(n+2)/(n-2)} \,\, \text{dans} \,\,  \Omega, \\
   \,\, \partial_{\nu} u_i =0                                    \qquad \qquad \qquad \qquad  \text{sur}   \quad \partial \Omega.           
 \end{cases}  $$

sont telles qu'\`a partir d'un certain rang $ i_0 $, $ u_i\equiv cte $

\end{Theorem}

Voir [1] des references supplementaires.

\bigskip

Si on suppose $ \Omega $ ouvert particulier de $ {\mathbb R}^n $, $ n \geq 4 $, dont le bord a une courbure moyenne $ H $ quelconque (en particulier on peut prendre $ H $ positive), alors on a:

\begin{Theorem}{\it(Wang-Wei-Yan)}. Sur $ \Omega $ de $ {\mathbb R}^n $ avec $ n \geq 4 $,  il existe une suite de solutions  non constantes de l'equation suivante:

$$ \begin{cases}
  -\Delta u_i + \epsilon_i u_i = u_i^{(n+2)/(n-2)} \,\, \text{dans} \,\,  \Omega, \\
   \,\, \partial_{\nu} u_i =0                                    \qquad \qquad \qquad \qquad  \text{sur}   \quad \partial \Omega.           
 \end{cases}  $$

telles que:

$$ \int_{\Omega} u_i^{2n/(n-2)} \to +\infty.. $$

\end{Theorem}

 Voir [26] des references supplementaires.

\smallskip

Si on suppose $ \Omega $ ouvert de $ {\mathbb R}^n $, $ n=3 $ ou $ n \geq 7 $, dont le bord a une courbure moyenne $ H $ strictement positive partout ($ H >0 $), c'est dire convexe dans un certain sens, alors on a:

\begin{Theorem}{\it(Druet-Robert-Wei)}. Sur $ \Omega $ de $ {\mathbb R}^n $ avec $ n= 3 $ ou $ n \geq 7 $, les solutions de l'equation suivante:

$$ \begin{cases}
  -\Delta u_i + \epsilon_i u_i = u_i^{(n+2)/(n-2)} \,\, \text{dans} \,\,  \Omega, \\
   \,\, \partial_{\nu} u_i =0                                    \qquad \qquad \qquad \qquad  \text{sur}   \quad \partial \Omega.           
 \end{cases}  $$

avec,

$$ \int_{\Omega} u_i^{2n/(n-2)} \leq \Lambda, \,\, \forall \,\, i. $$

sont telles qu'\`a partir d'un certain rang $ i_0= i_0(n, \Omega, \Lambda) $, $ u_i\equiv cte $

\end{Theorem}

Voir [10] des references supplementaires.

\bigskip

\underbar {\bf Quelques Propri\'et\'es des vari\'et\'es Riemanniennes.}

\bigskip

Soit $ (M,g) $ une vari\'et\'e Riemannienne.

\begin{Theorem}{\it (Bochner-Lichnerowicz-Weitzenb\"ock)}. Pour toute fonction $ u $ de classe $ C^3 $, on a la formule suivante:

$$ -\dfrac{1}{2} \Delta_g(|\nabla u|^2) = -|Hess (u)|^2 -<\nabla (\Delta_g u),\nabla u > - Ricci(\nabla u, \nabla u). $$
\end{Theorem}
On suppose que $ M $ est compacte sans bord de dimension $ \geq 2 $ et que le tenseur de Ricci v\'erifie l'in\'egalit\'e suivante:

$$ Ricci_{jk} \geq k g_{jk}, $$

avec $ k >0 $.

\bigskip

\begin{Theorem}{\it(Lichnerowicz)}. La premi\`ere valeur propore du Laplacien $ \lambda_1 $ v\'erifie:

$$ \lambda_1 \geq kn/(n-1). $$
\end{Theorem}

\begin{Theorem}{\it(Gidas-Spruck, Bidaut-V\'eron-V\'eron)}. Pour $ \epsilon >0 $, assez petit et $ n \geq 3 $, les solutions de l'\'equation suivante:

$$ -\Delta_g u + \epsilon u = u^{(n+2)/(n-2)}, $$

sont telles qu'\`a partir d'un certain rang $ i_0 $, $ u_i\equiv [\epsilon_i n(n-2)]^{(n-2)/4} $

\end{Theorem}

Voir [32] de la bibliographie et [4] des references supplementaires.

Consid\'erons une vari\'et\'e riemannienne  compacte $ (M,g) $  de courbure scalaire $ R_g $ ne verifiant pas necessairement la condition precedente sur le tenseur de Ricci. On a le resultat suivant de Brezis-Li qui donnent une preuve d'un resultat d'unicite sur la sphere unite par une methode de symetrie: 

\begin{Theorem} (Brezis-Li) Si $ n=3 $  et $ R_g $ quelconque, ou $ (M,g)= ({\mathbb S}_n, \delta) $, la suite de fonctions $ u_i >0 $ solutions de:

$$ -\Delta u_i+\epsilon_i u_i=n(n-2) u_i^{N-1}, $$

est telle qu'\`a partir d'un certain rang $ i_0 $, $ u_i\equiv [\epsilon_i n(n-2)]^{(n-2)/4} $
\end{Theorem}

Voir [17]

\bigskip

On a le r\'esultat plus g\'eneral suivant:

\begin{Theorem} (Bahoura [7]) Sur une vari\'et\'e Riemannienne compacte $ (M,g) $ de courbure scalaire $ R_g >0 $ partout, la suite de fonctions $ u_i >0 $ solutions de:

$$ -\Delta u_i+\epsilon_i u_i=n(n-2) u_i^{N-1}, $$

est telle qu'\`a partir d'un certain rang $ i_0 $, $ u_i\equiv [\epsilon_i n(n-2)]^{(n-2)/4} $
\end{Theorem}

\bigskip 

On a le r\'esultat suivant plus g\'eneral pour un systeme d'equations du type Yamabe:

$$ \begin{cases}

-\Delta_g u_i^k+ \sum_{j=1}^p A_{jk}(\epsilon) u_i^j=|U_i|^{4/(n-2)}u_i^k  \,\, {\rm dans }  \,\, M, \\

            u_i^k >0. 
             
 \end{cases}  $$

avec,

$$ |U_i|^2=\sum_{j=1}^p (u_i^j)^2 \,\, {\rm et }  \,\, ||A_{jk}(\epsilon)||\leq \epsilon. $$

\begin{Theorem} (Hebey) Sur une vari\'et\'e Riemannienne compacte $ (M,g) $ de dimension 3 et de courbure scalaire quelconque ou de dimension $ n\geq 4 $ et de courbure scalaire $ R_g >0 $ partout, la suite de fonctions $ u_i^k >0 $ solutions du systeme d'equations du type Yamabe precedent est telle qu'\`a partir d'un certain rang $ i_0 $, $ u_i^k\equiv cte.$
\end{Theorem}

Voir [14] des references supplementaires.

\section{Quelques estimations suppl\'ementaires.}

\subsection {La dimension 3 et r\'esultat du type unicit\'e.}

~~\\

En dimension 3, nous montrons que si $ u >0 $ est solution de 

\begin{Theorem} (YY.Li-L.Zhang) $$ -\Delta u+h(x)u=Vu^5,\,\, u >0, $$ 

avec, 

$$ 0 < a \leq V(x) \leq b, \,\,  ||\nabla V||_{\infty} \leq A, $$

alors, 

pour tout compact $ K $ de $ M $, il existe une constante positive $ c=c(a, b, A, h_0=||h||_{\infty}, K,M,g) $ telle que:

$$ \sup_K u \times \inf_M u \leq c. $$

\end{Theorem}

Voir [3] des r\'ef\'erences suppl\'ementaires. En particulier si $ h=\dfrac{1}{8}S_g $, L. Zhang obtient  une inegalit\'e optimale du type $ \sup \times \inf $ pour les solutions de l'equation de la courbure prescrite en dimension 3.

\bigskip

\chapter{} 

\section{Appendice.}

~~\\

Preuve de l'inegalite du type Harnack de Siu-Tian.

\bigskip

\underbar{\bf Une propri\'et\'e des vari\'et\'es K\"ahleriennes: In\'egalit\'e de Harnack.}

\bigskip

Soit $ (M,g) $ une vari\'et\'e K\"ahlerienne compacte de dimension complexe $ m $. Pour $ t\in ]0,1] $, on consid\`ere une fonction $ \phi_t $ solution de :

$$ \log M(\phi_t)=-t\phi +f. $$

avec, $ M ( \phi_t ) = det [ (g + \nabla^2 \phi_t)o g^{-1}]. $

\bigskip

{\bf Th\'eor\`eme}{\it (Tian)}. Il existe une constante $ C(t) $ telle que pour toute fonction admissible $ \psi $ qui satisfait:

$$ \int_M e^{f-t \psi} dV_g = V_g. $$

La solution $ \phi_t $ v\'erifie:

$$ \sup_M (\psi-\phi_t) \leq m \sup_M(\phi_t-\psi) + C(t). $$

De plus, si la m\'etrique initiale est Einstein-K\"ahler, alors pour toute fonction admssible $ \psi $ telle que:

$$ \int_M e^{-\psi} dV_g = V_g. $$

On a:

$$ \sup_M \psi + m \inf_M \psi \leq C. $$

La derni\`ere in\'egalit\'e est appell\'ee {\it in\'egalit\'e de type Harnack de Tian}.

\bigskip

\underbar{\it  Quelques \'el\'ements de la preuve:}

\bigskip

\underbar{\it  Un r\'esultat plus g\'en\'eral:}

\bigskip

Si on cherche une m\'etrique d'Einstein-K\"ahler sur une vari\'et\'e K\"ahlerienne compacte $ (M,g) $ avec comme premi\`ere forme de Chern $ C_1(M) > 0 $, on a besoin de r\'esoudre les \'equations de Monge-Amp\`ere suivantes:

$$ (*)_t  \,\,\, \begin{cases}

(\omega_g + \partial \bar \partial \phi )^n = e^{f-t \phi}\omega_g^n, \\

\omega_g + \partial \bar \partial \phi > 0 \,\,\, \text { sur} \,\,\, M. \

\end{cases}  $$

avec, $ \omega_g $ la forme de K\"ahler associ\'ee \'a la m\'etrique $ g  $, $ \omega_g^n = \omega_g \wedge \ldots \wedge \omega_g $ est la forme volume, $ 0 \leq t \leq 1 $, $ \partial \bar \partial f = Ricci(g) - \omega_g $, $ \int_M e^f  \omega_g^n =  \int_M \omega_g^n = Vol_g(M) $ et $ n = dim(M) $. 

\bigskip

\underbar{\it  Th\'eor\`eme :}

\bigskip

Soit $ (M,g) $ une vari\'et\'e K\"ahlerienne compacte avec $ C_1(M) > 0 $, $ n = dim(M) $. Alors, pour $ \psi  \in C^2(M,R) $ avec $ \omega_g + \partial \bar \partial \psi \geq 0 $ et $ \int_M e^{f-t \psi}\omega_g^n  = Vol_g(M) $, la solution $ \phi $ de  $ (*)_t $ satisfait l'in\'egalit\'e de Harnack du type:

$$ - \dfrac{1}{Vol_g(M)} \int_M (\phi- \psi ) (\omega_g + \partial \bar \partial \phi )^n \leq n \sup_M (\phi -\psi). \qquad (0.1) $$

De plus, il existe une constante $ C(t) $ d\'ependant seulement de $ t $ tel que pour $ t > 0 $, la solution $ \phi $ de $ (*)_t $ satisfait,

$$ - \inf_M (\phi - \psi) \leq n \sup_M  (\phi - \psi ) + C(t). \qquad (0.2) $$

\underbar{\it  Corollaire  :}

\bigskip

Il existe une constante universelle $ C $ telle que pour toute vari\'et\'e  d'Einstein-K\"ahler compacte $(M, g) $ avec $ C_1(M) > 0 $, i.e., $ Ricci(g) = \omega_g $, et pour toute fonction  $ \psi $ de classe $ C^2 $ avec $ \omega_g +  \partial \bar \partial \psi \geq 0 $,  $ \int_M e^{- \psi} \omega_g^n = Vol_g(M) $, on a l'in\'egalit\'e suivante,

$$ \sup_M \psi \leq - n \inf_M \psi + C. \qquad (0.3) $$

\underbar{\it Preuve du Th\'eor\`eme}

\bigskip

On suppose que $ \psi \in C^{\infty}(M,{\mathbb R}), \omega_g +\partial \bar \partial \psi > 0 $. Alors, on peut d\'efinir une nouvelle m\'etrique K\"ahlerienne $ \tilde g $ telle que la forme de K\"ahler associ\'ee est $ \omega_g + \partial \bar \partial \psi $. On consid\`ere $ \tilde f = f + \log (\omega_g^n / \omega_{\tilde g}^n) -t \psi $. Alors,

$$ \int_M e^{\tilde f} \omega_{\tilde g}^n = \int_M e^{f-t\psi} \omega_g^n = Vol_g(M). $$

On r\'ecrit $ (*)_t $ en termes de $ \omega_{\tilde g} $ et $ \tilde f $ comme,

$$ (1.2)  \,\,\, \begin{cases}

((\omega_g + \partial \bar \partial  (\phi - \psi))^n = e^{f-t (\phi -\psi)}\omega_g^n, \\

\omega_g + \partial \bar \partial (\phi - \psi) > 0 \,\,\, \text { sur} \,\,\, M. \\

\end{cases}  $$

Premi\`erement, on va prouver que les \'equations suivantes ont des solutions pour $ 0 \leq s \leq t $ :

$$ (1.3)_s  \,\,\, \begin{cases}

((\omega_g + \partial \bar \partial  \theta)^n = e^{f- s \theta}\omega_g^n, \\

\omega_g + \partial \bar \partial \theta  > 0 \,\,\, \text { sur} \,\,\, M. \\

\end{cases}  $$

On utilise la m\'ethode de continuit\'e. Soit $ S = \{ s\in [0,t] {\rm tel \,\, que} \,\, (1.3)_{s'} \,\, {\rm a \,\, une \,\, solution \,\, pour } \,\, s' \in [s,t] \} $. Comme $ (1.2) $ a une solution $ \phi, t\in S $ et $ S $ est non vide. Il est suffisant de montrer que $ S $ est \'a la fois ouvert et ferm\'e. Pour l'unicit\'e, on peut estimer la premi\`ere valeur propre de la m\'etrique $ g_{\theta} $ associ\'ee \'a la forme de K\"ahler $ \omega_g + \partial \bar \partial \theta $ pour la solution $ \theta $ de $ (1.3)_s $.

\bigskip

{\it Lemme 1.1}. La premi\`ere valeur propre non nulle $ \lambda_1(g_{\theta}) $ est plus grande que $ s $.

\bigskip

{\it Preuve}. D'apr\'es l'in\'egalit\'e bien connue de Bochner, il est suffisant de montrer que $ Ricci(g_{\theta}) $ est strictement minor\'e par $ s $. De $ (1.3)_s $, on a:

$$ Ricci(g_{\theta}) = Ricci(\tilde g) - \partial \bar \partial \tilde f + s \partial \bar \partial \theta $$

$$ = Ricci(\tilde g) - \partial \partial f - \partial \partial \log \left (\dfrac{\omega_g^n}{\omega_{\tilde g}^n} \right ) + t\partial \bar \partial \psi + s\partial \bar \partial \theta $$

$$ = Ricci(g) - \partial \bar \partial f + t\partial \bar \partial f + t\partial \bar \partial \psi +s\partial \bar \partial \theta $$

$$ = (1-t)\omega_g + (1-s)\omega_{\tilde g} +S(\omega_{\tilde g} + \partial \bar \partial \theta) > s \omega_{g_{\theta}}. $$

La premi\`ere variation de $ (1.3)_s $ en $ \theta $ est $ -\Delta_s u = -su $, o\`u $ \Delta_s $ est le Laplacien de la m\'etrique $ g_{\theta} $. Le lemme 1.1 implique que l'op\'erateur lin\'eaire $ -\Delta_s -s $ de $ (1.3)_s $ est inversible, le fait que S soit ouvert est une cons\'equence du Th\'eor\`eme des fonctions implicites.

\bigskip

Pour prouver que S est ferm\'e, on utilise la th\'eorie des \'equations elliptiques et les estimations de Yau des d\'eriv\'ees d'ordre sup\'erieures des solutions complexes des \'equations de Monge-Amp\`ere du type $ (1.3)_s $, il suffit d'estimer les solutions de $ (1.3)_s $ en normes $ C^0 $.

\bigskip

On suppose que $ (1.3)_s $ a une solution pour $ s \in (s_0,t] $ et $ \theta_s $ est la solution. On utilise la preuve du fait que S soit ouvert, on peut alors conclure que  $ \{ \theta \}_{s \in \{(s_0,t]\}} $ est une famille $ C^{\infty} $ dans $ C^{\infty}(M,{\mathbb R}) $, i.e., $ \theta_s $ varie de mani\`ere $ C^{\infty} $ en $ s $.

\bigskip

On consid\`ere deux quantit\'es introduites par T. Aubin,

$$ I(\theta_s) = \dfrac{1}{Vol_g(M)} \int_M \theta_s[\omega_{\tilde g}^n - (\omega_{\tilde g} + \partial \bar \partial_s)^n], \,\,\, J(\theta_s) =\int_0^1 \dfrac{I(x\theta_s)}{x} dx. $$

{\it Lemme 1.2}. (i) $ (n+1)J(\theta_s)/n \leq I(\theta_s) \leq (n+1) J(\theta_s) $,\\ (ii) $ \dfrac{d[I(\theta_s) - J(\theta_s)]}{ds} = -\dfrac{1}{Vol_g(M)} \int_M \theta_s(\Delta \dot \theta_s) \omega_g^n $.\\

Avec $ \dot \theta_s = {\dfrac{d\theta_s}{ds'}}_{|s'=s} $, $ g_s $ est la m\'etrique de K\"ahler associ\'ee \`a $ \omega_{\tilde g} + \partial \bar \partial \theta_s $ et $ \Delta_s $ son Laplacien.

\bigskip

En corollaire, on a le lemme suivant connu d\'eja par Bando et Mabuchi.

\bigskip

{\it Lemme 1.3}. $ I(\theta_s) - J(\theta_s) $ est d\'ecroissante.

\bigskip

{\it Preuve}. On diff\'erencie $ (1.3)_s $ en $ s $:

$$ -\Delta \dot {\theta_s} = -s \dot {\theta_s}-\theta_s. \qquad (1.4)_s $$

On remplace $ (1.4)_s $ dans le membre de droite dans la formule du lemme 1.2 (ii), on obtient,

$$ \dfrac{d}{ds} [I(\theta_s)-J(\theta_s)] = \dfrac{1}{Vol_g(M)} \int_M (-\Delta_s \dot {\theta_s} +s \dot {\theta_s} )-\Delta_s \dot {\theta_s} \omega_{g_s}^n. \qquad (1.5) $$

On \'ecrit $ \dot {\theta_s} $ dans le d\'eveloppement de la fonction propre, I.e.,

$$ \dot {\theta_s} = \sum_{i=0}^{+ \infty} a_i u_i, \qquad (1.6)_s $$

o\`u $ \Delta_s u_i = \lambda_i u_i, 0 = \lambda_0 < \lambda_1 \leq \lambda_2 \ldots $.

\bigskip

On utilise le lemme 1.1., $ \lambda > s $. Alors,

$$ \dfrac{d}{ds} [I(\theta_s) - J(\theta_s)] = \dfrac{1}{Vol_g(M)} \int_M \left [ \sum_{i=0}^{+ \infty} a_i(\lambda_i-s)u_i \right ] \left (\sum_{i=0}^{+ \infty} a_j \lambda_j u_j \right ) \omega_{g_s}^n $$

$$ = \dfrac{1}{Vol_g(M)} \sum_{i=0}^{+ \infty} |a_i|^2 \int_M (\lambda_i - s) \lambda_i |u_i|^2 \omega_{g_s}^n \geq 0, $$

et lemme est prouv\'e.

\bigskip

Dans ce qui suit, on note $ C $ une constante ind\'ependante de $ s $.

\bigskip

{\it Lemme 1.4}. Il existe une constante $ C > 0 $ telle que pour toute solution $ \theta_s $ de $ (1.3)_s $, $ 0 < s \leq t $, on ait, $ \sup_M |\theta_s| \leq C. $

\bigskip

{\it Preuve}. On d\'efinit un invariant holomorphe $ \alpha(M) $ sur la vari\'et\'e K\"ahlerienne compacte $ M $ avec $ C_1(M) > 0 $. Pour tout $ \lambda < \alpha(M) $, il existe une constante $ C_{\lambda} $ qui d\'epend de la m\'etrique $ \tilde g $, telle que,

$$ \int_M e^{-\lambda (u-\sup_M u)} dV_g \leq C_{\lambda} \,\,\, {\rm pour} \,\,\, u \in C^2(M, {\mathbb R}), \,\,\, \omega_{\tilde g} + \partial \bar \partial u \geq 0. \,\,\, (1.7) $$

Dans le cas o\`u $ s\in (0,\alpha(M)/(n+2)] $, $ \int_M e^{-(n+1)s(\theta_s - \sup_M \theta_s)} dV_{\tilde g} \leq C $ pour une constante $ C $. Pour $ p > 0 $,

$$ \int_M e^{-p(\theta_s - \sup_M \theta_s)} (e^{\tilde f-s\theta_s} -1) dV_{\tilde g} = \int_M e^{-p(\theta_s - \sup_M \theta_s)}(\omega_{g_s}^n -\omega_{\tilde g}^n) $$

$$ = \int_M e^{-p(\theta_s - \sup_M \theta_s)} \partial \bar \partial (\theta_s-\sup_M \theta_s)(\omega_{g_s}^{n-1} + \omega_{g_s}^{n-2} \wedge \omega_{\tilde g} + \ldots + \omega_{\tilde g}^{n-1}) $$

$$ = \dfrac{4}{p} \int_M \partial [e^{-p(\theta_s-\sup_M \theta_s)/2}] \wedge \bar \partial [e^{-p(\theta_s -\sup_M \theta_s)/2}] $$

$$ \wedge (\omega_{g_s}^{n-1} + \omega_{g_s}^{n-2} \wedge \omega_{\tilde g} + \ldots + \omega_{\tilde g}^{n-1} ) $$

$$ \geq \dfrac{4}{p} \int_M |\tilde \nabla [e^{-p(\theta_s -\sup_M \theta_s)/2}]|^2 dV_{\tilde g} $$

$$ \geq \dfrac{4c}{p} \left [ \int_M e^{-np(\theta_s-\sup_M \theta_s)/(n-1)} dV_{\tilde g} \right ]^{(n-1)/n} $$

$$ -\dfrac{4C}{p} \int_M e^{-p(\theta_s-\sup_M \theta_s)} dV_{\tilde g}, $$

o\`u $ c $ est la constante de Sobolev, ne d\'ependant que de $ (M,\tilde g) $. En utilisant l'in\'egalit\'e de H\"older dans le membre gauche de ce qui pr\'ec\`ede, on a,

$$ \left [ \int_M e^{-np(\theta_s-\sup_M \theta_s)/(n-1)} dV_{\tilde g} \right ]^{(n-1)/n} $$

$$ \leq C p \left [ \int_M e^{-(n+1)p(\theta_s -\sup_M \theta_s)/n} dV_{\tilde g} \right ]^{n/(n+1)} $$

$$ \times  \left [ \left [\int_M e^{-(n+1)s\theta_s} dV_{\tilde g} \right ]^{1/(n+1)} + 1 \right ] . \qquad (1.8) $$

On a, $ \sup_M \theta_s \geq 0 $, car,

$$ Vol_{\tilde g}(M) = \int_M \omega_{\tilde g}^n = \int_M \omega_{g_s}^n = \int_M e^{\tilde f- s\theta_s} dV_{\tilde g} \geq e^{-s \sup_M \theta_s} \int_M e^{\tilde f} dV_{\tilde g} $$

$$ = e^{-s\sup_M \theta_s} Vol_{\tilde g}(M). $$

Alors,

$$ \left [ e^{-(n+1)\theta_s} dV_{\tilde g} \right ]^{n+1} \leq \left [\int_M e^{-(n+1)s(\theta_s-\sup_M \theta_s)} dV_{\tilde g} \right ]^{n+1} $$

$$ \leq C^{n+1}. \qquad (1.9) $$

On remplace $ (1.9) $ dans $ (1.8) $, on a,

$$ \left |e^{-(\theta_s-\sup_M \theta_s)} \right |_{np/(n-1)} \leq C^{1/p}p^{1/p} \left |e^{-(\theta_s - \sup_M \theta_s)} \right|_{(n+1)p/n}. \qquad (1.10) $$

On met $ p = n s $ et $ p_{m+1} = p_m n^2/(n^2-1) $. Alors,

$$ \left | e^{-(\theta_s - \sup_M \theta_s)} \right |_{p_{m+1}} \leq C^{[n/(n-1)](1/p_{m+1})} \left ( \dfrac{n-1}{n}p_{m+1} \right )^{[n/(n-1)](1/p_{m+1}) } $$

$$ \times \left |e^{-(\theta_s-\sup_M \theta_s)} \right |_{p_m} $$

$$ (C p_m)^{1/p_m}\left |e^{-(\theta_s-\sup_M \theta_s)} \right |_{p_m} $$

$$ \ldots $$

$$ C^{1/ns} \sum_{m=0}^{+ \infty} \left (\dfrac{n^2-1}{n^2} \right )^m \left |e^{-(\theta_s-\sup_M \theta_s)} \right |_{p_0} $$

$$ \times \exp \left [\dfrac{1}{n p} \sum_{m=1}^{+ \infty} \left (\dfrac{n^2-1}{n^2} \right )^m \left ( m \log \dfrac{n^2}{n^2-1} + \log(ns) \right) \right ] $$

$$ \leq C, $$

et il s'en suit que,

$$ - \inf_M (\theta_s -\sup_M \theta_s) = \log \left [ \lim_{m \to + \infty} \left |e^{-(\theta_s-\sup_M \theta_s)} \right |_{p_m} \right ] \leq C, $$

i.e., pour $ s \in (0, \alpha(M)/(n+2)] $, $ \sup_M |\theta_s| \leq C $.

\bigskip

Dans le cas o\`u $ s\geq \alpha(M)/(n+2) $, $ Ricci(g_s) \geq s \geq \alpha(M)/(n+2) $. Alors, en utilisant l'identit\'e de Bochner et, les r\'esultats de Crkoe et de Peter Li, on a, en usant de l'in\'egalit\'e de Sobolev et l'in\'egalit\'e de Poincar\'e avec leurs constantes uniform\'ement born\'ees sur $ (M,g_s) $. Comme $ \Delta_s \theta_s \geq -n $, l'it\'eration de Moser implique que,

$$ -\inf_M \theta_s \leq C \int_M (-\theta_s) \omega_{g_s}^n + C, \qquad (1.11) $$

D'autre part, en utilisant la formule de Green sur $ (M,\tilde g) $, il s'en suit que,

$$ \sup_M \theta_s \leq \int_M \theta_s \omega_{\tilde g}^n + C. \qquad (1.12) $$

En utilisant le lemme 1.3 et le point (i) du lemme 1.2,

$$ I(\theta_s) \leq (n+1) [I(\theta_s) -J( \theta_s)] \leq (n+1) [I(\theta_t) -J(\theta_t)] \leq C.  \qquad (1.13) $$

Comme $ \int_{\{\theta_s > 0 \} } \theta_s e^{\tilde f -s\theta_s} dV_{\tilde g} $ et $ \int_{\{\theta_s < 0 \} } (-\theta_s) dV_{\tilde g} $ sont born\'ees par une constante $ C $ ind\'ependante de $ s $, le lemme provient de $ (1.11)-(1.13) $ et de la d\'efinition de $ I(\theta_s) $.

\bigskip

Le fait que S soit ferm\'e provient du lemme pr\'ec\'edent. Alors, $ (1.3)_s $ a une solution pour $ 0 \leq s \leq t $. Il s'en suit qu'il existe une famille r\'eguli\`ere de $ \{ \theta_s \}_{ s \geq 0} $ telle que $ \theta_t =\phi - \psi $. En utilisant le lemme 1.2 (ii) et $ (1.4)_s $, on a,

$$ \dfrac{d}{ds} [I(\theta_s)-J(\theta_s)] = \dfrac{1}{Vol_{\tilde g}(M)} \int_M \theta_s (-s \dot {\theta_s} - \theta_s) \omega_{g}^n, $$

$$ = \dfrac{1}{Vol_{\tilde g}(M)}\dfrac{d}{ds}\left ( \int_M \theta_s  \omega_{g_s}^n \right ) + \dfrac{1}{Vol_{\tilde g}(M)} \int_M \dot {\theta_s}  \omega_{g_s}^n. $$

En diff\'erentiant $ Vol_{\tilde g}(M) =  \int_M e^{\tilde f - s\theta_s} dV_{\tilde g} $, on aura,

$$ \int_M (-s \dot {\theta_s} - \theta_s) e^{\tilde f - s\theta_s} dV_{\tilde g} = 0.  \qquad (1.14) $$

Alors,

$$ \dfrac{d}{ds} [I(\theta_s)-J(\theta_s)]  = \dfrac{1}{ s Vol_{\tilde g}(M)} \dfrac{d}{ds}\left ( \int_M \theta_s e^{\tilde f - s \theta_s} \omega_{g_s}^n \right ); $$

I.e.,

$$  \dfrac{d}{ds} [ s [I(\theta_s)-J(\theta_s)]  ] - [I(\theta_s)-J(\theta_s)]   $$

$$  = \dfrac{d}{ds} \left ( \dfrac{s}{Vol_{\tilde g}(M)} \int_M - \theta_s e^{\tilde f - s \theta_s} \omega_{g_s}^n \right ) . \qquad (1.15) $$

Notons que,

$$ I(\theta_s)-J(\theta_s)  \geq \dfrac{1}{n+1}I(\theta_s) = \dfrac{1}{(n+1)Vol_{\tilde g}(M)} \int_M  \theta_s (\omega_{\tilde g}^n  - \omega_{g_s}^n )  $$

$$ =  \dfrac{1}{(n+1)Vol_{\tilde g}(M)} \int_M [ \partial \theta_s  \wedge \partial \bar \theta_s $$

$$  \wedge (\omega_{\tilde g}^{n-1} + \omega_{\tilde g}^{n-2}  \wedge \omega_{g_s} +  \ldots + \omega_{ g_s}^{n-1} ) ] \geq 0. $$
 
Une cons\'equence de $ (1.15) $ et le point (ii) du lemme 1.2 est,
 
$$ \dfrac{1}{Vol_{\tilde g}(M)} \int_M  (- \theta_t )  ( \omega_{g_s}^n ) \leq  I(\theta_t)-J(\theta_t) \leq   \dfrac{n}{n+1} I(\theta_t)  $$

$$ =  \dfrac{n}{(n+1)Vol_{\tilde g}(M)} \int_M  \theta_t (\omega_{\tilde g}^n  - \omega_{g_t}^n ) ,\qquad (1.16) $$

i.e.,

$$ - \dfrac{1}{Vol_{\tilde g}(M)} \int_M (\phi - \psi ) ( \omega_g + \partial \bar \partial \phi )^n  \leq \dfrac{n}{Vol_{\tilde g}(M)} \int_M (\phi - \psi ) \omega_{\tilde g}^n $$

$$ \leq \sup_M (\phi - \psi ), $$

c'est justement  $ (0.1) $. L'in\'egalit\'e $ (0.2) $ provient de l'it\'eration de Moser et du fait que $ Ricci(g_t) \geq t > 0 $. Le Th\'eor\`eme 1 est prouv\'e. (Voir le lemme 1.4 pour plus de d\'etails).

\bigskip

\underbar {\it Preuve du Corollaire : }

\bigskip

Du point $ (0.2) $ dan le Th\'eor\`eme 1, pour toute $ \psi \in C^2(M, {\mathbb R})$  avec  $ \omega_g + \partial \bar \partial \psi  \geq 0 $ et  $ \int_M e^{ - \psi} \omega_g^n = Vol_g(M) $, on a,

$$ -\inf_M ( \phi - \psi ) \leq n \sup_M ( \phi - \psi ) + C(1), \qquad (2.1) $$

o\`u $ \psi $  est la solution de $ (*)_1 $ et $ C(1) $ est une constante universelle. Notons qu'ici $ f \equiv  0 $, car $ g $ est une m\'etrique d'Einstein-K\"ahler.  Ce qui implique que $ \phi \equiv 0 $ est une solution de $ (*)_1 $. Pour $ \phi = 0 $, $ (2.1) $ devient,

$$ -\inf_M ( - \psi ) \leq n \sup_M ( - \psi ) + C(1).  $$

Comme $  \inf_M ( - \psi ) =  - \sup_M (  \psi ) $ et $  -sup_M ( - \psi ) =  -\inf_M ( \psi ) $, le corollaire est prouv\'e.

\bigskip

\underbar{\bf Principe de concentration-compacite de Brezis-Merle.}

\bigskip

Soit $ \Omega $ un ouvert de $ {\mathbb R}^2 $ et $ V $ une fonction telle que:

$$ 0 \leq V \leq b, $$

pour un r\'eel positif $ b $ donn\'e.

\bigskip

On consid\`ere l'\'equation suivante:

$$ -\Delta u = V e^u. $$

{\bf Th\'eor\`eme.}{\it (Brezis-Merle)}. Soient $ (u_n) $ et $ (V_n) $ deux suites de fonctions solutions de l'\'equation pr\'ec\'edente. On suppose qu'il existe deux constantes $ C_0 $ et $ C_1 $ telles que:

$$ 0 \leq V_n \leq C_0, $$

$$ \int_{{\Omega}} e^{u_n} \leq C_1. $$

Il existe une sous-suite not\'ee $ (u_j) $ telle que:

\bigskip

1. La suite $ (u_j) $ est born\'ee dans $ L_{loc}^{\infty} $, ou,

\bigskip

2. $ u_j \to - \infty $ uniform\'ement sur tout compact de $ \Omega $, ou,

\bigskip

3. Il existe un ensemble fini de points blow-up $ S = \{ a_1, \ldots, a_m \} $ tel que, pour tout $ k $ il existe une suite de points $ x_{j,k} $, $ x_{j,k} \to a_k $, $ u_j(x_{j,k}) \to + \infty $ et $ u_j \to - \infty $ sur tout compact de $ \Omega - S $. De plus, $ V_j e^{u_j} \to \sum_{k=1}^m \alpha_k \delta_{a_k} $ au sens des distributions, avec, $ \alpha_k \geq 4 \pi $.

\bigskip

\underbar {\bf Preuve}

\bigskip

Avant de prouver le Th\'eor\`eme, nous devons donner quelques d\'efinitions.
\bigskip

L'ensemble des points blow-up:  

\bigskip

S = $ \{  $ x $ \in  \Omega $, tel qu'il existe une suite $ (x_n) $ dans  $  \Omega $ telle que  $ x_n  \to x $ et $ u_n (x_n)  \to + \infty \} $.
   
\bigskip 
Comme  la suite $ (V_n e^{u_n})_n $ est bornee dans $ L^1 $ on peut en extraire une sous-suite qui converge au sens des mesures vers une mesure positive ou nulle $ \mu $.

On dit qu'un point $ x_0 $ est un point r\'egulier de $ \mu $ s'il existe une fonction $ \psi \in C_c(\Omega) $, $ 0 \leq \psi \leq 1 $, avec $ \psi =1 $ dans un voisinage de $ x_0 $, tel que :

$$ \int \psi d\mu < 4 \pi. \qquad (1)  $$

 L'ensemble des points non r\'eguliers: on le note $ \Sigma $.

Une cons\'equence du fait qu'un point $ x_0 $ est r\'egulier, est :

$$ \exists \,\, R_0 > 0 \,\, tel \,\, qu'on  \,\, puisse  \,\, borner  \,\, (u_n^+) \,\, dans  \,\, L^{\infty}(B_{R_0}(x_0)).\qquad (2)  $$ 

On commence la preuve du Th\'eor\`eme:

\bigskip

\underbar { Etape 1}: S = $ \Sigma $.

\bigskip

Il ais\'e de voir que S  $ \subset \Sigma $. Inversement, soit $ x_0 \in \Sigma $. Alors, on a:

$$   \forall \,\,\,  R >0,\,\,\,  \lim ||u_n^+ ||_{ L^{\infty }(B_R(x_0))} = + \infty. \qquad (3) $$

Sinon, il existerait $ R_0 > 0 $ et une sous-suite telle que:

$$  ||u_n^+ ||_{ L^{\infty }(B_{ R_0}(x_0))} \leq C. $$

En particulier,

$$  ||e^{u_{n_k}} ||_{ L^{\infty }(B_{ R_0}(x_0))} \leq C, \,\,\, {\rm et  }$$

$$ \int_{B_{ R_0}(x_0)}     V_{n_k}e^{u_{n_k}} \leq C. $$

Ce qui implique $ (1) $ pour une fonction $ \psi $ particuli\`ere et $ x_0 $ serait un point r\'eguli\`er, contradiction. Aisni, on a \'etablit $ (3) $. On choisit $ R_ 0 > 0  $ assez petit pour que $ B_{ R_0}(x_0) $ ne contient pas un autre point de $ \Sigma $. Soit $ x_n \in B_R(x_0)  $ tel que :

$$ u_n^+(x_n) = \max_{B_R(x_0)} u_n^+ \to + \infty. $$

On a $ x_n \to x_0 $. Sinon, il existerait une sous-suite $ x_{n_k} \to \bar x \not = x_0 $ et $ \bar x  \not \in \Sigma $, i.e. $ \bar x $ est un point r\'egulier. Ceci n'est pas possible si on utilise $ (2) $. Ce qui prouve que $ x_0 \in S $ et l'\'etape 1 est prouv\'ee.

\bigskip

\underbar { Etape 2}: S = $ \emptyset $ implique que  $ 1) $ ou $ 2) $ est vraie.

\bigskip

En utilisant $ (2) $, $ (u_n^+) $ est born\'ee dans $ L_{loc}^{\infty}(\Omega) $ et alors $ f_n = V_n e^{u_n} $ est born\'ee dans $ L_{loc}^{p}(\Omega) $. Ceci implique que $ \mu \in L^1(\Omega) \cap  L_{loc}^{p}(\Omega) $. Soit $ v_n $ la solution de :

$$ -\Delta v_n = f_n \,\, {\rm dans } \,\,\, \Omega,\,\,\,{\rm et} \,\,\, v_n = 0 \,\,\,{\rm sur } \,\,\, \partial  \Omega. $$

Clairement, $ v_n \to v $ uniform\'ement sur tout compact de $ \Omega $ et $ v $ est solution de :

$$ -\Delta v = \mu \,\, {\rm dans } \,\,\, \Omega,\,\,\,{\rm et} \,\,\, v = 0 \,\,\,{\rm sur } \,\,\, \partial  \Omega. $$

Soit $ w_n = u_n - v_n $ tel que $ \Delta w_n = 0 $ sur $ \Omega $ et $ w_n^+ $ est born\'ee dans $  L_{loc}^{\infty}(\Omega) $. En utilisant le principe de Harnack, on trouve que:

$$ Il\,\, existe \,\, une \, sous-suite \,\, (w_{n_k}) \,\, qu'on \,\, borne \,\, dans \,\, L_{loc}^{\infty}(\Omega). $$

ou bien

$$ (w_n) \,\, converge \,\, vers  \,\, -\infty  \,\, dans \,\, L_{loc}^{\infty}(\Omega). $$

Ces deux derniers cas correspondent aux cas $ 1) $ et $ 2) $.

\bigskip

\underbar { Etape 3}: S $\not = \emptyset $. implique que  $ 3) $ est vraie.

\bigskip

En utilisant $ (2) $, $ (u_n^+) $ est born\'ee dans $ L_{loc}^{\infty}(\Omega -  S) $ et alors $ f_n $ est born\'ee dans $ L_{loc}^{\infty}(\Omega -  S) $. Ce qui implique  que $ \mu $ est born\'ee sur $ \Omega $ et $ \mu \in L_{loc}^{p}(\Omega -  S) $. Comme dans l'\'etape 2, on d\'efinit $ v_n $, $ v $ et $ w_n $. Alors, $ v_n \to v $ uniform\'ement sur tout compact de $ \Omega - S $. Comme ci-dessus, le principe de Haranck donne :

$$ Il\,\, existe \,\, une \, sous-suite \,\, (w_{n_k}) \,\, qu'on \,\, borne \,\, dans \,\, L_{loc}^{\infty}(\Omega - S). $$

ou bien

$$ (w_n) \,\, converge \,\, vers  \,\, -\infty  \,\, dans \,\, L_{loc}^{\infty}(\Omega- S). $$

Icic, on prouve que le premier cas n'est pas possible. On fixe un point $ x_0 \in S $ et $ R > 0 $ assez petit tel que $ x_0 $ soit le seul point de $ S $ dans $ \bar B_R(x_0) $. Supposons que le premier point soit vrai, alors, il existe une sous-suite $ (u_{n_k}) $ born\'ee dans $ L^{\infty}(\partial B_R(x_0)) $ par $ C $. Soit  $ (z_{n_k}) $ la solution de :

$$ -\Delta z_{n_k} = f_{n_k} \,\,\, {\rm dans } \,\,\, B_R(x_0),\,\,\,{\rm et} \,\,\, z_{n_k} = -C \,\,\,{\rm sur } \,\,\, \partial  B_R(x_0). $$

Par le principe du maximum, $ u_{n_k} \geq z_{n_k} $ dans $ B_R(x_0) $.

\bigskip

En particulier,

$$ \int  e^{z_{n_k}}  \leq  \int  e^{u_{n_k}}  \leq C. \,\,\,\,\, (*)$$

D'autre part, $ z_{n_k} \to z $ p.p. ( sur tout compact de  $ B_R(x_0) -  \{ x_0\} $) avec $ z $ solution de :

$$ -\Delta z = \mu \,\,\, {\rm dans } \,\,\, B_R(x_0),\,\,\,{\rm et} \,\,\, z = -C \,\,\,{\rm sur } \,\,\,  \partial B_R(x_0). $$

Finallement, comme $ x_0 \in S $ n'est pas un point r\'egulier,  on a $ \mu(\{ x_0\}) \geq 4 \pi $, ce qui implique que,  $ \mu \geq 4 \pi \delta_{x_0} $ et alors, par le principe du maximum dans $ W^{1,1}_0(B_R(x_0)) $ (obtenu par l'inegalite de Kato)

$$ z(x) \geq  2 \log \dfrac{1}{| x - x_0 |} + O(1) \, \, {\rm si } \, \, x  \to x_0. $$

Donc, 

$$ e^{z} \geq  \dfrac{C}{| x - x_0 |^2}, \, \, C > 0. $$

Par cons\'equent:

$$  \int_{B_R(x_0)}  e^z = \infty. $$

D'autre part, on utilise $ (*) $ et lemme de Fatou pour obtenir:

$$ \int  e^z   \leq C. $$

Ce qui est contradictoire. 

On suppose:

$$ 0 < a \leq V \leq b <+\infty $$

{\bf Corollaire}{\it (Brezis-Merle)}. Si $ m > -\infty $, alors:

$$ \sup_K u \leq c=c(a, b, m, K, \Omega), \,\, {\rm si} \,\, u \geq m, $$

o\`u $ K $ est compact de $ \Omega $.

\bigskip

\underbar{\bf Inegalites de type Harnack sur un ouvert du plan euclidien.}

\bigskip

Preuves des inegalites du type Harnack de Shafrir et Brezis-Li-Sahfrir:

\bigskip

{\bf Th\'eor\`eme}{\it (Shafrir)}. Il existe une constante $ C=C \left ( \dfrac{a}{b} \right ) $ telle que:

$$ C \sup_K u + \inf_{\Omega} u \leq c(a, b, K, \Omega), $$

o\`u $ K $ est un compact de $ \Omega $.

\bigskip

\underbar {\bf Preuve}

\bigskip

\underbar {\bf Lemme}

\bigskip

Par soucis de simplification,  on suppose que $ \Omega = B_1 $. Pour, $ a, b >0 $, il existe $ \alpha_0 > 4 \pi $ et une constante positive $ C_0 $ telle que : 

$$ u(0) \leq C_0 \,\,\, {\rm  si } \,\,\, -\Delta u = V e^u, \,\,\,  a \leq  V \leq b \,\,\,{\rm  et } \,\,\, \int_ {B_1} V e^u \leq \alpha_0. $$

\underbar {\bf Preuve du lemme}

\bigskip

\underbar {\it Etape 1: In\'egalit\'e isop\'erim\'etrique d'Alexandrov.}

\bigskip

On consid\`ere le disque $ B \subset  {\mathbb R}^2 $ et une fonction $ u  \in C^2(B) $ v\'erifiant  $  \Delta u = V e^u $ dans $ B $. On pose $ d\sigma^2 = e^u (dx_1^2 + dx_2^2) $, on consid\`ere  $ (B, d\sigma) $, comme une surface. Son aire est donn\'ee par $ M =  \int_B e^u dx $ et la longueur du bord est donn\'ee par, $ L= \int_{\partial B} d\sigma =  \int_{\partial B} e^{ u / 2} ds $. La courbure scalaire est $ K=V/2 $. Pour un r\'eel $ K_0 $, on note :

$$  \omega_{K_0}^+(B) =  \int_{\{ x \in \beta, K(x) > K_0 \}} (K-K_0) e^u  dx. $$

On obtient l'in\'egalit\'e isop\'erim\'etrique suivante, dite d'Alexandrov :

$$ L^2 \geq (2 \alpha-K_0 M) M. \qquad (1) $$

\underbar {\it Etape 2: Application et preuve du lemme.}

\bigskip

Sans perdre en g\'en\'eralit\'e, on suppose $ V  $ r\'eguli\`ere, le cas g\'en\'eral sera trait\'e par approximation. Pour $ r \in (0,1) $, on note $ a(r) = \int_ {B_r} V e^u dx $. On choisit un $ K_0 $ optimal pour  $ (1) $ avec $ B = B_1 $. C'est donn\'e par un $ K_m $ pour lequel:

$$ \int_ {\{ K > K_m \}}  e^u \leq a(1)/2  \,\,\, {\rm  et } \,\,\,\int_ {\{ K < K_m \}}  e^u \leq a(1)/2. $$

On \'ecrit :

$$ \int_ {\{ K > K_m \}}  e^u = a(1)/2 - A  \,\,\, {\rm  et } \,\,\,\int_ {\{ K < K_m \}}  e^u = a(1)/2-B, $$

avec, $ A $ et $ B $ positifs ou nuls. Pour $ K_m $, on a:

$$ 2 \alpha-K_m a(1) = 4 \pi - \int_ {\{ K \geq K_m \} }  K e^u dx + 2 K_m B. $$

Ensuite, on veut prouver que cette quantit\'e est strictement positive si $  \int_ {B_1}  K e^u dx $ est assez proche de $ 2 \pi  $.

$$ \int_ {\{ K  < K_m \}}  K e^u dx \geq \dfrac{a}{2} \int_ {\{ K  < K_m \}}  e^u dx = \dfrac{a}{2} \left ( \dfrac{ a(1) }{2}-B \right ) \geq   \dfrac{a}{2b} \int_ {B_1}  K e^u dx -  \dfrac{a B}{2}. $$

Alors,

$$ \int_ {\{ K   \geq K_m\}}  K e^u dx \leq  \left ( 1 - \dfrac{a}{2b} \right )  \int_ {B_1}  K e^u dx + \dfrac{a B}{2}. $$

On fixe $ \alpha_0 \in ( 4 \pi, 4 \pi / (1-(a/2b))) $, avec, $ 2 \alpha - K_m a(1) \geq 4 \pi - 2 (1-(a/2b)) \int_{B_1}  K e^u dx \geq \gamma_0 > 0 $ quand $  \int_{B_1}  V e^u dx  \leq  \alpha_0 $.

\bigskip

On peut appliquer $ (1) $ pour $ B_1 $. On consid\`ere la fonction absolue continue suivante:

$$ f(r) =  4 \pi - 2 \int_{ \{ K   \geq K_m  \} \cap B_r }  (K - K_m)e^u dx - K_m \int_{ B_r} e^u dx. $$

Clairement, $ f(r) $ est strictement d\'ecroissante, alors, $ f(r) > f(1) \geq \gamma_0 $ pour tout $ r \in (0,1) $, on applique $ (1) $ pour $ B = B_r $.

On \'ecrit:

$$ a(r) = \int_0^r \int_{\partial B_{ \tilde r } } e^u ds d\tilde r. $$

On utilise l'in\'egalit\'e de Cauchy-Schwarz et $ (1) $:

$$ a'(r) =  \int_{\partial B_r } e^u ds \geq \dfrac{1}{2\pi r} \left ( \int_{\partial B_r } e^{u/2} ds  \right )^2  \geq  \dfrac{1}{2\pi r} f(r) a(r). $$

En fixant  un $ r_0 > 0 $ assez petit puis en int\'egrant, on obtient:

$$   \int_0^1 \dfrac{a'(r)}{f(r) a(r)} dr  \geq -  \dfrac{1}{2\pi}  \log r_0. $$

Une int\'egration par partie donne:

$$  \int_0^1 \dfrac{a'(r)}{f(r) a(r)} dr  =   \dfrac{ \log a(1)}{f(1)}- \dfrac{  \log a(r_0)}{f(r_0)} + \int_0^1 \dfrac{ \log a(r) f'(r)}{f^2(r)} dr. $$  

Ensuite, on prend $ r_0  \to 0 $. Comme $ a(r_0)  $ est \'equivalent  \'a $ \pi r_0^2 e^{u(0)} $ et, $ f(r_0) $ est \'equivalent  \'a  $ 4 \pi - O(r_0^2) $, on conclut que :

$$  \lim_{ r_0  \to 0}    \left ( \dfrac{ \log a(r_0)}{f(r_0)}- \dfrac{ \log r_0}{2 \pi}   \right ) = \dfrac{u(0) + \log \pi }{4 \pi}. $$

et,

$ \int_0^1 \dfrac{ \log a(r) f'(r)}{f^2(r)} dr $ converge. Donc,

$$ u(0) \leq   \dfrac{ 4 \pi }{ f(1) }\dfrac{ \log a(1) }{ \pi } + \int_0^1 \dfrac{ \log a(r) f'(r)}{f^2(r)} dr. $$

Le deuxieme membre est born\'e par une constante $ C(a,b) $ car, $ -f'(r) \leq (3b-2a)a'(r) $. D'o\`u le lemme.

\bigskip

\underbar {\it  Preuve du Th\'eor\`eme.}

\bigskip

Clairement, il est suffisant de prouver le th\'eor\`eme pour des boules. On pose $ w(x) = u(rx)+ 2\log r $, cette fonction est solution de notre \'equation sur $ B_1 $  avec la m\^eme condition sur $ V $. Il est suffisant de prouver:

$$ u(0) + C_1 \inf_{B_1}  u  \leq C_2, $$

pour une solution $ u $ de notre \'equation avec la m\^eme condition sur $ V $. On peut avoir mieux, si, on remplace $  \inf_{B_1} $ par $ (1/2 \pi) \int_ {\partial B_1} u ds $. Soit $ r  \in (0,1) $ et $ G $ la fonction  suivante:

$$ G(r) = u(0) + \dfrac{ C_1}{2\pi r } \int_ {\partial B_r} u ds + 2 (C_1 +1)  \log r . $$

Avec, $ C_1 $ une constante qu'on d\'eterminera plustard. On d\'erive $ G $ et on cherche son maximum.
 
$$ G'(r) \geq 0  \Leftrightarrow  \dfrac{C_1}{2 \pi } \left ( \int_0^{2 \pi}  u(r e^{i \theta}) d\theta \right ) + \dfrac{ 2 ( C_1 + 1)}{ r } \geq 0  \Leftrightarrow   \dfrac{ C_1}{2 \pi }  \int_{ {\partial B}_r}  \partial_r u ds +  2 ( C_1 + 1)  \geq 0. $$

Mais,

$$ \int_{ {\partial B}_r}  \partial_r u ds = \int_{ B_r}  \Delta u dx = - \int_{ B_r}  V e^u dx. $$

Alors,

$$ G'(r) \geq 0 \Leftrightarrow   \int_{ B_r}  V e^u dx  \leq \dfrac{4 (C_1 + 1)}{C_1}. $$

Si,  $  \int_{ B_r}  V e^u dx  > \dfrac{4 (C_1 + 1)}{C_1} $ on prend $ r_0  $ tel que $  \int_{ B_{r_0}}  V e^u dx  = \dfrac{4 (C_1 + 1)}{C_1} $, sinon, on prend $ r_0 = 1 $. Dans chaque cas, $ G(1) \leq G(r_0) $. On choisit $ C_1 $ assez grand pour avoir $  \dfrac{4 (C_1 + 1)}{C_1}  \leq  \alpha_0 $ pour le $  \alpha_0 $ comme dans le lemme. Ensuite, on utilise la super-harmonicit\'e de $ u $ et le lemme pour avoir:

$$ G(r_0) = u(0) + \dfrac{ C_1}{2 \pi r_0}  \int_{ {\partial B}_{r_0}}   u ds + 2 (C_1 + 1)  \log r_0  \leq (C_1+1) (u(0) + 2  \log r_0 )  \leq (C_1+1)C_0.  $$

On pose $ C_2 = C_0(C_1+1) $, on obtient:

$$ u(0) + \dfrac{ C_1}{2 \pi r_0}  \int_{ \partial B_1}   u ds = G(1) \leq C_2. $$

et d'o\`u le r\'esultat.

\bigskip

On suppose que $ V $ est lipschitzienne avec:

$$ ||\nabla V||_{L^{\infty}(\Omega)} \leq A. $$

{\bf Th\'eor\`eme.}{\it (Brezis-Li-Shafrir)}. Il existe une constante $ c=c(a, b, A, K, \Omega) $ telle que:

$$ \sup_K u + \inf_{\Omega} u \leq c, $$

o\`u $ K $ est un compact de $ \Omega $.

\bigskip

\underbar{\bf Preuve}

\bigskip

Soit $ A = ||\nabla V||_{L^{\infty}} $. La preuve est divis\'ee en 5 \'etapes.

\bigskip

\underbar{\it Etape 1:} R\'eduction \`a $ \Omega = B_2 $ ( la boule centr\'ee en 0 et de rayon 2 et

$$ u(0) + \inf_{B_2} u \leq C(a, b, A). \qquad (9) $$

Le cas g\'en\'eral provient de $ (9) $. En effet, supposons que $ (9) $ est vraie et soit $ v $ la solution de:

$$ - \Delta v = V e^v \qquad {\rm sur } \qquad B_R. $$ 

Alors,

$$ u(x) = v \left ( \dfrac{R}{2} x \right ) + 2 \log (R/2) $$

v\'erifie,

$$ -\Delta u = V \left ( \dfrac{R}{2} x \right ) e^v \qquad {\rm sur} \qquad B_2 $$

et alors, par $ (9) $,

$$ v(0) + \inf_{B_R} v \leq C(a, b, RA/2) -4 \log (R/2) = C(a, b, A, R). \qquad (10). $$

Le Th\'eor\`eme  est donn\'e par $ (10) $.

\bigskip

Dans ce qui suit, on raisonne par l'absurde et on suppose que $ (9) $ n\'est pas vraie. Plus pr\'ecis\'ement, on suppose qu'il existe une suite $ (u_n) $ de solutions de,

$$ - \Delta u_n =V_n e^{u_n} \qquad {\rm sur} \qquad (11) $$

avec,

$$ a \leq V_n \leq b, \qquad ||\nabla V_n||_{L^{\infty}} \leq A \qquad (12) $$

telle que,

$$ u_n(0) + \inf_{B_2} u_n \to + \infty \qquad (13) $$

Apr\'es passage au sous-suites, on peut supposer que $ V_n \to V $ uniform\'ement sur $ \bar B_2 $ avec $ V(0) = K \geq a >0 $. Soit,

$$ \delta_n = e^{-u_n(0)/2} \qquad (14) $$

\underbar{\it Etape 2:}

$$ \delta_n \to 0 \qquad (15) $$

et,

$$ \limsup \int_{B_{R \delta_n}} V_n e^{u_n} \leq 8 \pi \,\,\, {\rm pour \,\, tout } \,\,\, R > 0 \qquad (16) $$

{\it Preuve.} 

On a,

$$ u_n(0) + \inf_{B_2} u_n \leq 2 u_n(0) $$

et alors, par $ (13) $, $ u_n(0) \to + \infty $. On introduit la fonction,

$$ G(r) = u_n(0) + \dfrac{1}{2\pi r} \int_{\partial B_r} u_n ds + 4\log r, \qquad 0 < r \leq 2. $$

Comme,

$$ G'(r) = \dfrac{1}{2\pi r} \int_{\partial B_r} \partial_r u_n ds + \dfrac{4}{r} $$

et,

$$ \int_{\partial B_r} \partial_r u_n ds = -\int_{B_r} \Delta u_n dx = -\int_{B_r} V_n e^{u_n}, $$

on conclut que,

$$ G'(r) \geq 0 \Leftrightarrow \int_{B_r} V_n e^{u_n} \leq 8 \pi \qquad (17) $$

et,

$$ G'(r) =0 \Leftrightarrow \int_{B_r} V_n e^{u_n} = 8 \pi \qquad (18) $$

La fonction $ G(r) $ atteint son maximum sur $ [0,2] $ en un point $ 0 < \mu_n \leq 2 $. Si $ \mu_n < 2 $, on a:

$$ \int_{B_{\mu_n}} V_n e^{u_n} = 8 \pi. $$

Sinon, $ \mu_n = 2 $ et on a,

$$ \int_{B_{\mu_n}} V_n e^{u_n} \leq 8 \pi. $$

Donc, dans tous les cas,

$$ \int_{B_{\mu_n}} V_n e^{u_n} \leq 8 \pi. $$

Comme $ u_n $ est super-harmonique, on a,

$$ 2[u_n(0)+ 2 \log \mu_n] \geq u_n(0) + \dfrac{1}{2\pi \mu_n} \int_{\partial B_{\mu_n}} u_n ds + 4 \log \mu_n $$

$$ = G(\mu_n) \geq G(2) = u_n(0) + \dfrac{1}{4\pi} \int_{\partial B_2} u_n ds + 4 \log 2 $$

$$ \geq u_n(0) + \inf_{\partial B_2} u_n + 4 \log 2 \geq u_n(0) + \inf_{B_2} u_n + 4 \log 2. $$

En utilisant $ (13) $, on conclut que,

$$ u_n(0) + 2 \log \mu_n \to + \infty. $$

i.e.,

$$ \log(\mu_n/\delta_n) \to +\infty $$

et donc, $ \mu_n /\delta_n \to + \infty $. Pour $ R >0 $ et $ n $ assez grand, $ R \delta_n \leq \mu_n $ et donc,

$$ \int_{B_{R\delta_n}} V_n e^{u_n} \leq \int_{B_{\mu_n}} V_n e^{u_n} \leq 8 \pi. $$

\underbar{Etape 3}. Il existe une suite $ x_n \to 0 $ et $ R_n > 0 $ tel que (pour une sous-suite),

$$ |x_n| < R_n \leq 1, $$

$$ x_n \,\,{\rm est\,\, un \,\,maximum \,\, de \,\, } u_n \,\, {\rm sur} \,\, B_{R_n}(x_n) \qquad (21) $$

$$ R_n e^{u_n(x_n)/2} \to + \infty \qquad (22) $$

et,

$$ \limsup \int_{B_{R_n}(x_n)} V_n e^{u_n} \leq 8 \pi. \qquad (23) $$

{\it Preuve.} Soit,

$$ v_n(x) = u_n(\delta_n x) + 2 \log \delta_n \qquad {\rm pour} \,\, |x|\leq 1/\delta_n. $$

On consid\`ere les restrictions des $ (v_n) $ \`a $ B_1 $, elles v\'erifient,

$$ -\Delta v_n = V_n(\delta_nx) e^{v_n} \qquad {\rm sur} \,\,\, (24) $$

De $ (16) $ ( avec $ R = 1 $) et $ (12) $ on en d\'eduit que,

$$ \limsup \int_{B_1} e^{v_n} \leq \dfrac{8\pi}{a}. \qquad (25) $$

On est maintenant en situation d'appliquer la technique {\it blow-up} du Th\'eor\`eme pr\'ec\'edent de Br\'ezis-Merle. Il y a trois possibilit\'es:

\bigskip

Cas 1. $ (v_n) $ est born\'e dans $ L_{loc}^{\infty}(B_1) $.

\bigskip

Cas 2. $ v_n \to - \infty $ uniform\'ement sur tout compact de $ B_1 $.

\bigskip

Cas 3. Il existe un ensemble non vide $ S $ de $ B_1 $ de points {\it blow-up} tel que $ v_n \to - \infty $ uniform\'ement sur tout compact de $ B_1 - S $ et pour chaque point $ a \in S $ il existe une suite $ (a_n) $ telle que $ a_n \to a $ et $ v_n(a_n) \to + \infty $.

\bigskip

Comme $ v_n(0) = 0 $, le cas 2 est exclu. On examine les cas 1 et 3 s\'epar\'ement.

\bigskip

Cas 1. On consid\`ere $ (v_n) $ restrinte \`a $ B_R $ pour un $ R >1 $ fix\'e. Pour $ n $ assez grand, $ (v_n) $ v\'erifie $ (24) $ et $ (25) $ (et $ B_1 $ est remplac\'ee par $ B_R $). On applique le Th\'eor\`eme de Br\'ezis-Merle dans $ B_R $ et on voit que $ (v_n) $ est born\'ee dans $ W_{loc}^{2,p}(B_R) $ pour chaque $ p < + \infty $. Alors, en passant au sous-suites, on peut supposer que $ (v_n) $ converge dans $ C_{loc}^1({\mathbb R}^2) $ vers une fonction $ v $ satisfesant,

$$ v \in L_{loc}^{\infty}({\mathbb R}^2), $$

$$ -\Delta v = K e^v \,\,{\rm sur} \,\, {\mathbb R}^2 \,\, (K= \lim V_n(0)), $$

$$ \int_{{\mathbb R}^2} e^v \leq \dfrac{8\pi}{a}, $$

et,

$$ v(0) = 0. $$

Il s'en suit que $ v $ est de la forme,

$$ v(x) = \log \left \{ \dfrac{8\lambda^2 /K}{(1 + \lambda^2 |x-y_0|^2)^2} \right \} $$

pour un point $ y_0 \in {\mathbb R}^2 $ et un $ \lambda >0 $. Pour $ \rho > |y_0| $ le maximum de $ v_n $ sur $ \bar B_{ \rho } $ est atteint en un $ y_n $. Clairement $ y_n \to y_0 $ et $ v_n \to v $ uniform\'ement sur $ B_{ \rho } $. En particulier, pour $ k $ entier assez grand, on a $ (v_{n_k}) $ et $ (y_{n_k}) $ telles que,

$$ \max_{B_k} v_{n_k} \,\, v_{n_k} \,\, {\rm est \,\, atteint \,\, en} \,\, y_{n_k} $$

et,

$$ y_{n_k} \to y_0 $$

Comme $ \mu_n/\delta_n \to + \infty $, on peut supposer que,

$$ k \delta_{n_k} \leq \dfrac{1}{2} \mu_{n_k}.  $$

Soit,

$$ x_{n_k} = \delta_{n_k} y_{n_k} \qquad {\rm et} \qquad R_{n_k} = (k - |y_{n_k}|) \delta_{n_k}. $$

Il est ais\'e de voir que les sous-suites correspondantes v\'erifient $ (20)-(23) $.

\bigskip

Cas 3. Clairement $ 0 _in S $ (sinon, on peut avoir $ v_n(0) \to - \infty $, mais $ v_n(0) = 0 $). On peut choisir $ r_0 \in (0,1) $ tel que $ (v_n) $ a un autre point blow-up dans $ B_{r_0} $ exc\'ept\'e \`a l'origine. Pour chaque $ n $ soit $ y_n $ le maximum de $ v_n $ sur $ B_{r_0} $. Alors, d'apr\'es l'assertion "blow-up", $ v_n(y_n) \to + \infty $ et $ y_n \to 0 $. Soit $ x_n = \delta_n y_n $ et $ R_n = \dfrac{1}{2} r_0 \delta_n $. Il est ais\'e de voir que les propri\'et\'es $ (20)-(23) $ sont satisfaites.

\bigskip

On pose, sur $ B_1 $,

$$ \bar u_n = u_n(x + x_n) \qquad {\rm et} \qquad \bar V_n(x) = V_n(x + x_n), $$

alors, on a,

$$ -\Delta \bar u_n = \bar V_n e^{\bar u_n} \qquad {\rm sur} \qquad B_1 \qquad (26) $$

$$ 0 \,\, {\rm est \,\, un \,\, point \,\, maximum \,\, de}, \,\, \bar u_n \,\, {\rm sur },\,\, B_{R_n}, \,\, {\rm avec }\,\, 0 < R_n \leq 1, \,\,\, (27) $$

$$ R_n e^{\bar u_n(0) /2} \to + \infty, \qquad (28) $$

$$ \limsup \int_{B_{R_n}} \bar V_n e^{\bar u_n} \leq 8 \pi, \qquad (29) $$

et,

$$ \bar u_n(0) + \inf_{B_1} \bar u_n \to + \infty. \qquad (30) $$

\underbar{Etape 4}. Soit,

$$ \eta_n = e^{\bar u_n(0)/2}, \,\,\, {\rm alors,} \,\,\, \eta_n \to 0, $$

$$ \bar v_n(x) = \bar u_n(\eta_n x) + 2 \log \eta_n, \,\,\, {\rm pour,} \,\,\, | x | \leq 1/\eta_n, $$

et,

$$ \bar w_n(x) = \bar v_n(x) + 2 \log | x |, \,\,\, {\rm pour,} \,\,\, | x | \leq 1/ \eta_n. $$

Clairement, $ \bar v_n $ v\'erifie,

$$ -  \Delta  \bar v_n = V_n( \eta_n x) e^ {\bar v_n} \,\,\, {\rm pour,} \,\,\, | x | \leq 1/\eta_n, $$

$$ \bar v_n(0) = 0, $$

et pour chaque $ R $,

$$ \max_{B_R}   \bar v_n \,\,\, {\rm est \,\, atteint \,\, en } \,\, 0 \,\, {\rm pour } \,\, n  \,\,  {\rm assez  \,\, grand }, $$

$$ \limsup \int_{B_R} e^{\bar u_n} \leq 8 \pi /a, \,\,\, {\rm en  \,\, utilisant \,\,  (28)  \,\, et \,\,  (29) }. $$

On peut utiliser le Th\'eor\`eme de Br\'ezis-Merle pour conclure que $ \bar v_n $ est born\'ee dans $ L_{loc}^{\infty}({\mathbb R}^2) $. En utilisant les estimations elliptiques, on trouve que $ \bar v_n $ est aussi born\'ee dans $ C_{loc}^{1, \alpha}({\mathbb R}^2) $. Alors, pour une sous-suite, $ \bar v_n $ converge dans $ L_{loc}^{\infty}({\mathbb R}^2) $ vers une fonction $  \bar v $ satisfaisant,

$$ - \Delta \bar v = K e^{\bar v} \,\,\, {\rm sur } \,\,\, {\mathbb R}^2, $$

$$ \bar v(0) = 0, $$

et,

$$ \int_{{\mathbb R}^2} e^{\bar v}  < +  \infty. $$

Rn utilisant un r\'esultat connu, $  \bar v $ est donn\'ee par,

$$  \bar v(x) =  \log \left \{ \dfrac{1}{(1+\gamma^2 |x|^2)^2} \right \}, $$

avec, $ \gamma (K/8)^{1/2} $. Il s'en suit que:

$$ \bar w_n - \bar w \to 0  \,\,\, {\rm dans} \,\,\,  C_{loc}^2({\mathbb R}^2), \qquad (31) $$

avec,

$$  \bar v(x) =  \log \left \{ \dfrac{|x|^2}{(1+\gamma^2 |x|^2)^2} \right \}. $$

Dans la suite, on travaille en coordon\'ees polaires $ (r, \theta ) $ et on pose $ t = \log r $. Soit, $ t < 0 $ et $ \theta \in (0,2\pi) $,

$$ \tilde w_n (t, \theta) = \bar u_n(e^t \cos \theta, e^t \sin \theta) + 2 t. \qquad (32) $$

Clairement $ \tilde w_n $ v\'erifie,

$$ -  \Delta  \tilde w_n = \tilde V_n(t, \theta) e^ {\tilde v_n} \,\,\, {\rm dans} \,\,\, Q, $$ 

avec,

$$ Q = \{ (t, \theta); t \leq 0,\,\, {\rm et } \,\, 0 \leq  \theta \leq 2 \pi \}, $$

$$ \Delta = \partial_{tt} + \partial_{\theta \theta}, $$

et,

$$ \tilde V_n(t, \theta) =  \bar V_n(e^t \cos \theta, e^t \sin \theta). $$

On introduit, pour $ s \in {\mathbb R}^2 $,

$$   \tilde w(s) = \log \left \{ \dfrac{e^{2s}}{(1+\gamma^2 e^{2s})^2} \right \} = 2s - 2  \log (1+\gamma^2 e^{2s}). $$

Notons que $ \tilde w $ atteint son maximum en $ s = - \log  \gamma $, $ \tilde w'(s) > 0 $ pour $ s < -  \log  \gamma  $, et $ \tilde w $ est symm\'etrique par rapport \'a $ s = - \log \gamma $. On utilise le fait que,

$$ \tilde w(s) \leq 2s \,\,\, \forall  s \in {\mathbb R}^2. \qquad (33) $$

Clairement, on a,

$$ \tilde w_n (s + \log \eta_n, \theta) - \tilde w(s) = \bar w_n(e^s \cos \theta, e^s \sin \theta) - \bar w(e^s \cos \theta, e^s \sin \theta). $$

Dans les nouvelles variables $ (31) $ implique que, pour tout $ \alpha \in {\mathbb R} $, quand $ n \to +\infty $,

$$ || \tilde w_n (s + \log \eta_n, \theta) - \tilde w(s) ||_{L^{\infty} \{ s \leq \alpha,\theta \in (0, 2\pi) \}} \to 0. $$

En particulier, on peut choisir $ n_0 $ suffisemment grand pour que, pour $ n \geq n_0 $, on a,

$$ | \tilde w_n (t, \theta)- \tilde w(t-\log \eta_n) || \leq 1\,\,\, {\rm si} \,\,\, t \leq 4 - \log  \gamma +  \log  \eta_n,\,\, 0\leq \theta \leq 2\pi, \,\,\, (34) $$

et,

$$  \tilde w_n (-\log \gamma + \log \eta_n, \theta) > \tilde w_n(-\log \gamma + \log \eta_n + 4, \theta) \,\,\, {\rm si} \,\,\,0\leq \theta \leq 2\pi.  \,\,\, (35) $$

Finalement, on introduit,

$$ \hat w_n(t, \theta) = \tilde w_n(t, \theta)- \dfrac{A}{a} e^t \,\,\, {\rm dans} \,\,\, Q. \qquad (36) $$

On clame que,

$$ \partial_t \left \{  \tilde V_n(t, \theta) e^{Ae^t/a}e^{\xi} + \dfrac{A}{a} e^t \right \} \geq 0 \,\,\, \forall  (t,\theta) \in Q, \,\,
\forall \,\,  \xi \in {\mathbb R}. \,\,\, (37) $$

Ceci provient du fait que:

$$ \tilde V_n  \geq a  \,\,\, {\rm et} \,\,\, |\partial_t \tilde V_n(t,\theta) | \leq Ae^t. $$

\underbar{Etape 5}. {\it (Conclusion via la m\'ethode de r\'eflexion )}. 

\bigskip

On utilise la m\'ethode "moving-plane " introduite par Alexandrov, d\'evelopp\'ee par Gidas-Ni-Nirenberg et utilis\'ee au paravant par Schoen.

\bigskip

Pour $ \lambda < 0 $ et $ \lambda \leq t \leq 0 $, on pose,

$$ t^{\lambda} = 2 \lambda - t $$

et,

$$ \hat w_n^{\lambda} (t, \theta) = \hat w_n(t^{\lambda}, \theta). $$

We have,

$$ - \Delta (\hat w_n^{\lambda} -\hat w_n) = \hat V_n^{\lambda}(t, \theta)e^{\hat w_n^{\lambda}} - \hat V_n(t, \theta)e^{\hat w_n} + \dfrac{A}{a} (e^{t^{\lambda}} -e^t), \qquad (38) $$

avec, $ \hat V_n (t, \theta) = \tilde V_n^{\lambda}(t, \theta)e^{Ae^t/a}  $ et $ \hat V_n^{\lambda}(t, \theta) = \hat V_n (t^{\lambda}, \theta) $.

\bigskip

Pour un $ \lambda $ assez n\'egatif (d\'ependant de $ n $), on a,

$$ \hat w_n^{\lambda} (t, \theta) - \hat w_n (t, \theta)  < 0 \,\,\, {\rm pour} \,\,\, \lambda < t \leq 0, \,\,  0 \leq  \theta 2\pi.  \qquad (39) $$

Pour prouver $ (39) $ on utilise le fait que pour $ n $ fix\'e, on a, en utilisant $ (32) $ et $ (36) $,

$$  \hat w_n (t, \theta) = 2t + a_n + O_n(e^t) \,\,\, {\rm quand } \,\,\, t \to - \infty $$

et,

$$ \partial_t  \hat w_n (t, \theta) = 2 + O_n(e^t)  \,\,\, {\rm quand } \,\,\, t \to - \infty. $$

On d\'efinit,

$$  \lambda_n = \sup \{ \lambda < 0;  \hat w_n^{\lambda} (t, \theta) - \hat w_n (t, \theta)(t, \theta) < 0 \,\, {\rm pour }\,\, \lambda < t \leq 0, \,\, 0 \leq  \theta \leq 2 \pi \}. $$

On clame que,

$$ \lambda_n \leq - \log \gamma +  \log \eta_n + 2. \qquad (40) $$

En effet, si on choisit $  \lambda_n = - \log \gamma +  \log \eta_n + 2 $ et $ t =  \lambda_n \leq - \log \gamma +  \log \eta_n + 4 $ alors $ t^{\lambda}  = - \log \gamma +  \log \eta_n $ et, en utilisant $ (35) $, 
$ \hat w_n^{\lambda} (t, \theta)  >  \hat w_n (t, \theta) $, $ \forall \,\, \theta \in (0,2\pi) $.

\bigskip

D'autre part, si on use de $ (38) $, $ (37) $ et la d\'efinition de $ \lambda_n $,

$$ - \Delta (\hat w_n^{\lambda} (t, \theta) - \hat w_n (t, \theta)) \leq 0 \,\,\, {\rm pour } \,\,\, \lambda \leq t \leq 0, \,\,\lambda \leq \lambda_n  \,\, {\rm et} \,\,0 \leq  \theta \leq 2 \pi. $$

Maintenant, on clame que,

$$ \min_{0 \leq \theta \leq 2\pi} \hat w_n^{\lambda} (0, \theta)  \leq \max_{0 \leq \theta \leq 2\pi}   \hat w_n (2\lambda_n, \theta). \qquad (41) $$

Supposons que ce ne soit pas le cas,

$$  \max_{0 \leq \theta \leq 2\pi}   \hat w_n (2\lambda_n, \theta)  <  \min_{0 \leq \theta \leq 2\pi} \hat w_n^{\lambda} (0, \theta); $$

alors, en utilisant le principe du maximum,

$$ \hat w_n^{\lambda} (t, \theta) - \hat w_n (t, \theta)  < 0 \,\,\, {\rm pour } \,\,\, \lambda_n < t < 0, \,\,0 \leq  \theta \leq 2 \pi. $$

et le lemme de Hopf donne,

$$ \partial_t (\hat w_n^{\lambda} (t, \theta) - \hat w_n (t, \theta))_{ t  = \lambda_n} < 0, \,\,\, 0 \leq  \theta \leq 2 \pi. $$

Ceci contredit la d\'efinition de $ \lambda_n $.

\bigskip

En utilisant $ (34) $ puis $ (33) $, on a,

$$ \max_{0 \leq \theta \leq 2\pi}   \hat w_n (2\lambda_n, \theta)  \leq   \tilde w (2\lambda_n- \log \eta_n) + 1 \leq 4 \lambda_n - 2 \log \eta_n + 1. $$

Donc, si on use de $ (40) $, on obtient,

$$ \max_{0 \leq \theta \leq 2\pi}   \hat w_n (2\lambda_n, \theta)  \leq \leq 2 \log \eta_n + C( \gamma ). $$

En combinant $ (41) $ et $ (42) $, on peut voir que,

$$ \min_{0 \leq \theta \leq 2\pi}   \hat w_n (0, \theta)  \leq 2 \log \eta_n + C( \gamma ). $$

Si, on regarde la d\'efinition de $  \hat w_n $, on a,

$$ \min_{ \partial B_1}   \bar  u_n   \leq 2 \log \eta_n + C(a, A, \gamma ).  \qquad (43) $$

Comme $ \eta_n = e^{- \bar u_n(0)/2} $, et donc, si on utilise $ (43) $,

$$ \bar u_n(0) + \min_{ \partial B_1} \bar u_n \leq C(a, A, \gamma), $$

ce qui contredit $ (30) $.

\bigskip

On suppose de plus que $ V $ uniformement $ \alpha $-holderienne de constante $ A $ alors:

\bigskip

{\bf Th\'eor\`eme}{\it (Chen-Lin)}(sans preuve). Il existe une constante $ c=c(a, b, A, \alpha, K, \Omega) $ telle que:

$$ \sup_K u + \inf_{\Omega} u \leq c, $$

o\`u $ K $ est un compact de $ \Omega $.

\bigskip

\underbar{\bf Inegalites de type Harnack sur une variete Riemannienne quelconque de dimension 3 et 4.}

\bigskip

Preuve de l'inegalite du type Harnack de Li-Zhang en dimensions 3 et 4:

Sur une vari\'et\'e Riemannienne quelconque $ (M,g) $ (non n\'c\'essairement compacte), de dimension $ n\geq 3 $ et de courbure scalaire $ R_g $, on consid\`ere l'\'equation de Yamabe:

$$ -\Delta_g u + R_g u = u^{(n+2)/(n-2)}, \,\,\, u >0 . $$

{\bf Th\'eor\`eme}{\it (Li-Zhang)}. En dimensions $ 3, 4 $, pour tout compact $ K $ de $ M $, on a:

$$ \sup_K u \times \inf_M u \leq c = c(K, M, g, n = 3 \,\, {\rm ou}\,\, n= 4). $$

\underbar{\it Preuve:}

\bigskip

\underbar{\it Cas n=3:}

\bigskip

On raisonne par l'absurde, il existe une constante $ \bar a > 0 $ pour laquelle il y a une suite de m\'etriques $ \{g_k\} $, une suite $ (\epsilon_k) $ avec $ \epsilon_k \to 0^+ $ et une suite $ (u_k) $ solutions de notre \'equation avec,

$$ \max_{\bar B(0,\epsilon_k)} u_k \times \min_{\bar B(0, 4\epsilon_k)} u_k > k{\epsilon_k}^{2-n}, \qquad (1) $$

o\`u $ B(0,\epsilon_k) $ est la boule g\'eod\'esique relative \`a $ g_k $.

\bigskip

D'autre part, il ais\'e de voir qu'il existe $ \bar \epsilon = \bar \epsilon (n,\bar a) > 0 $ telle que le principe du maximum est vrai pour l'op\'erateur $ L_g = \Delta_g -c(n) R_g $ sur $ ]0, \bar \epsilon [ $, avec $ c(n) =\dfrac{n-2}{4(n-1)} $ et $ R_g $ la courbure scalaire.

\bigskip

Comme, $ L_g u_k \leq 0 $, on a,

$$ \min_{\bar B(0,r)} u_k = \min_{\partial B(0,r)} u_k, \,\,\, \forall \,\, 0 < r \leq \bar \epsilon. \qquad (2) $$

{\bf Remarque: Dans leur preuve Li-Zhang utilisent la version generale du principe du maximum. Pour voir ca, on considere une fonction positive $ u $ verifiant $ Lu \leq 0 $. On choisit une fonction $ u_0 >0 $ telle que  $ u= u_0 v $ verifie une  EDP sans terme lineaire pour l'operateur:

$$ L u= u_0 Lv+ (termes \, d'ordre \, 1) \partial v +  v (L u_0) \leq 0, $$

Il suffit de resoudre  (voir Gilbarg-Trudinger, solution reguliere):

$$ Lu_0 = 0, \,\,dans \,\, \Omega \, u_0 =1\,\, sur\,\,  \partial \Omega .$$

On chosit l'ouvert de depart de telle maniere que $ L $ satisfasse le principe du maximum d'Alexandrov (forme faible). Alors par ce principe du maximum, $ u_0 \leq 0 $ et par le principe du maximum de Gidas-Ni-Nirenberg $ -u_0 < 0 $.

Donc, $ u_0 >0 $ et $ v $ verifie une EDP sans terme lineaire donc:

$$ \forall \,\, K \subset \Omega, \, \min_{\partial K} v= \min_{K} v,  $$

Comme $ u_0 >0 $ sur $ \bar \Omega $,  et $ u \geq 0 $, on obtient $ \exists \,\, c_1= c_1(\Omega) >0, c_2=c_2(\Omega) >0  $ telles que:

$$ \forall \, K \subset \Omega, \, \, c_1 \min_{\partial K}  u \leq \min_K u \leq c_2 \min_{\partial K} u. $$     

Cela est suffisant dans le preuve de Li-Zhang.

Via l'exponentielle $ \exp_{x_k}(y) $ on ramene les boules $ B(x_k, r_0) $ vers la boule de l'espace euclidien $ B(0, r_0) $. Alors $ L $ depend de $ k $ et devient $ L_k $ mais les coeficients  $ a^k_{ij}, b^k_j, c_k $ sont uniformement bornes et uniformement corecif. Il suffit de consider le coeficient lineaire $ c $ continue.

On aura,

$$ L_ku^k_0 = 0, \,\,dans \,\, \Omega=B(0, r_0), \, u^k_0 =1\,\, sur\,\,  \partial \Omega =\partial B(0, r_0).$$

Le reel positif $ r_0 $ est choisit au depart ne dependant que des parametres exterieurs des coeficients de $ L_k $ de telle maniere que $ L_k $ satisfasse le principe du maximum d'Alexandrov.

Par le principe du maximum d'Alexandrov (inegalite d'Alexandrov-Bakelman-Pucci) $ 0 < \sup u^k_0 \leq 1+ O(c) $, (ici aussi, le $ r_0>0 $ est choisit de telle maniere que cette inegalite soit verifiee, $ r_0 $ depend des parametres exterieurs des coefficients, on considere $ v^k_0= u^k_0-1 $, alors $ (a^k_{ij} \partial_{ij}+ b^k_j \partial_j -c_k^-) v^k_0 = (L_k +c_k^+) v^k_0 \geq -c_k-c_k^+ (v^k_0)^+ $ avec condition au bord nulle, on obtient $ v^k_0 \leq (v^k_0)^+ \leq  O(c_k) $, donc, $ 0 < u^k_0 \leq 1+ O(c_k) $), puis $ 1-u^k_0 $, on utiise les estimations elliptiques pour avoir  la convergence dans $ C^2 $ de $ u^k_0 $, puis on applique les principes du maximum d'Alexandrov et de Gidas-Ni-Nirenberg a l'operateur limite $ L_0 $ et la fonction limite $ u_0 $ pour avoir $ u_0 >0 $. Finalement cela revient a le faire pour un operateur au lieu d'une suite d'operateurs.

}

Pour un $ \bar x_k \in \bar B(0,\epsilon_k) $, $ u_k(\bar x_k) =\max_{\bar B(0,\epsilon_k)} u_k $, et, d'apr\'es ce qui pr\'ec\`ede,

$$ u_k(\bar x_k) {\epsilon_k}^{(n-2)/2} \to + \infty. $$

On peut trouver $ x_k \in B(\bar x_k,\epsilon_k /2) $ et $ \sigma_k \in (0,\epsilon_k/4) $ satisfaisant,

$$ u_k(x_k)^{2/(n-2)} \sigma_k \to + \infty, \qquad (3) $$

$$ u_k(x_k) \geq u_k(\bar x_k), \qquad (4) $$

et,

$$ u_k(x) \leq C_1 u_k(x_k), \,\,\, \forall \,\, x\in B(x_k,\sigma_k), \qquad (5) $$

o\`u $ C_1 $ est une constante universelle.

\bigskip

Il s'en suit de $ (4) $, $ (2) $ et $ (1) $ que,

$$ u_k(x_k) \times \min_{\partial B(x_k, 2\epsilon_k)} u_k \times {\epsilon_k}^{2-n} \geq u_k(\bar x_k) \times \min_{\bar B(0, 4\epsilon_k)} u_k \times {\epsilon_k}^{2-n} \geq k \to + \infty.  \qquad (6) $$

On note $ \{ z^1, \ldots, z^n \} $ les coordon\'ees g\'eod\'esiques normales centr\'ees en $ x_k $. Dans ces coordon\'ees, $ g = g_{ij}(z) dz^i dz^j $,

$$ g_{ij}(z) =\delta_{ij} + O(r^2),\,\,\, g:=det[g_{ij}(z)] =1+ O(r^2),\,\,\, R_g(z) = O(1), \,\,\, (7) $$

avec, $ r = |z| $. Alors,

$$ \Delta_g =\dfrac{1}{\sqrt g} \partial_i(\sqrt g g^{ij}\partial_ju) = \Delta u + b_i\partial_i u + d_{ij}\partial_{ij} u, $$

o\`u,

$$ b_j = O(r),\,\,\, d_{ij} = O(r^2). \qquad (8) $$

Ici, on note, $ \partial_i =\dfrac{\partial}{\partial z^i} $ et $ \partial_{ij} =\dfrac{\partial}{\partial z^i \partial z^j} $.

\bigskip

L'\'equation de $ u_k $ peut s\'ecrire,

$$ L_g u_k + u_k^{(n+2)/(n-2)} = \Delta u_k + b_i \partial_i u_k + d_{ij} \partial_{ij} u_k - c(n)R_g u_k + u_k^{(n+2)/(n-2)} = 0 \,\,\, B(0, 3\epsilon_k) \,\, (9) $$

On pose,

$$ v_k(y) =M_k^{-1} u_k \left (M_k^{-2/(n-2)} y \right ) \,\,\, {\rm pour} \,\,\, |y| \leq 3 \epsilon_k M_k^{2/(n-2)}, $$

o\`u $ M_k = u_k(0) $. En utilisant $ (3) $ et le fait que $ \epsilon_k \geq 4 \sigma_k $,

$$ \lim_{k \to + \infty} \epsilon_k M_k^{2/(n-2)} = \lim_{k \to + \infty} \sigma_k M_k^{2/(n-2)} = + \infty, \qquad (10) $$

En utilisant $ (9) $, $ (6) $ et $ (7) $,

$$ (11) \,\,\, \begin{cases}

\Delta v_k + \bar b_i \partial_i v_k + \bar d_{ij} \partial_{ij} v_k -\bar cv_k + v_k^{(n+2)/(n-2)} = 0 \,\,\, \text{pour}\,\, -y| < 3\epsilon_k M_k^{2/(n-2)}, \\

v_k(0)=1, \\

v_k(y) \leq C_1, \,\,\, \text{pour} \,\,\, |y| \leq \sigma_kM_k^{2/(n-2)}, \\

\lim_{k \to +\infty} \min_{|y|=2\epsilon_kM_k^{2/(n-2)}} \left (v_k(y)|y|^{n-2} \right ) = + \infty, 

\end{cases} $$

o\`u $ C_1 $ est une constante universelle dans $ (6) $,

$$ \bar b_i(y) = M_k^{-2/(n-2)} b_i(M_k^{-2/(n-2)}y),\,\,\, \bar d_{ij}(y) = d_{ij}(M_k^{-2/(n-2)}y), \qquad (12) $$

et,

$$ \bar c(y) =c(n)R\left (M_k^{-2/(n-2)} y \right ) M_k^{-4/(n-2)}. \qquad (13) $$

Ici, par soucis de simplifications, on omet la d\'ependance de $ \bar b_i, \bar d_{ij} $ et $ \bar c $ en $ k $.

\bigskip

Pour $ |y| \leq 3 \epsilon_k M_k^{2/(n-2)} $, on a, en utilisant $ (8) $,

$$ |\bar b_i(y)|\leq CM_k^{-4/(n-2)}|y|, \,\, |\bar d_{ij}(y)| \leq C M_k^{-4/(n-2)} |y|^2, \,\, |\bar c(y)|\leq CM_k^{-4/(n-2)}, \qquad (14) $$

o\`u $ C $ est une constante qui ne d\'epend que de $ n $ et $ \bar a $.

\bigskip

Il s'ensuit de $ (33) $, $ (34) $ et $ (37) $, en utilisant les estimations elliptiques standards, que, apr\'es passage au sous-suites, $ v_k $ converge en norme $ C^2 $ sur tout compact de $ {\mathbb R}^n $ vers une fonction positive $ U $ solution de,

$$ (15) \,\,\, \begin{cases}

\Delta U + U^{(n+2)/(n-2)} = 0 \,\, \text{ sur} \,\, {\mathbb R}^n, \\

U(0) =1, \,\, 0 < U \leq C_1. 

\end{cases} $$

Par soucis de simplifications, on notera $ v_k $ toute sous-suite de $ v_k $.

\bigskip

Pour tout $ a > 0 $, il existe une constante $ c(a) > 0 $ et $ \bar k_a > 0 $, ind\'ependant de $ k $, tel que,

$$ c(a) < v_k(y) \leq v_k(y) + |\nabla v_k(y)| + |\nabla^2 v_k(y)| \leq \dfrac{1}{c(a)}, \,\,\, \forall \,\, |y| \leq a \,\, {\rm et} \,\, k \geq \bar k_a. \,\,\, (16) $$

For $ x\in {\mathbb R}^n $ et $ \lambda > 0 $, soit,

$$ v_k^{\lambda,x}(y) : = \left (\dfrac{\lambda}{|y-x|} \right )^{n-2} v_k \left (x + \dfrac{\lambda^2(y-x)}{|y-x|^2} \right ) $$

la transformation de Kelvin de $ v_k $ pour une boule centr\'ee en $ x $ et de rayon $ \lambda $.

\bigskip

Nous allons comparer, pour  tout $ x $ fix\'e, $ v_k $ et $ v_k^{\lambda,x} $. Par soucis de simplifications, on prendera $ x = 0 $. Pour $ x \not = 0 $, les arguments restent les m\^emes. On notera $ v_k^{\lambda} $, la fonction $ v_k^{\lambda, 0} $, i.e.

$$ v_k^{\lambda}(y) : = \left (\dfrac{\lambda}{|y|} \right )^{n-2} v_k \left ( y^{\lambda} \right )   \,\, { \rm avec} \,\, y^{\lambda} =   \dfrac{\lambda^2y}{|y|^2} . $$

On pose, pour $ \lambda > 0 $,

$$ \Sigma_{\lambda} = B \left (0, \epsilon_k M_k^{2/(n-2)} \right ) - \bar {B(0, \lambda)}. $$

On se placera dans $  B ( 0, \epsilon_k M_k^{2/(n-2)} ) $. Si, $ x \not = 0 $, on remarquera que $ B ( x, \epsilon_k M_k^{2/(n-2)} )  \subset B ( 0, 2 \epsilon_k M_k^{2/(n-2)} ) $ pour $ k $ assez grand, ce qui revient \`a consid\'erer le cas $ x = 0 $. On peut \'ecrire:

$$ \min_{ |y| = \epsilon_k M_k^{2/(n-2)} }  (v_k(y) |y|^{n-2} ) \geq 2^{2-n}  \min_{ |y| = 2\epsilon_k M_k^{2/(n-2)} }  (v_k(y) |y|^{n-2} ) \to + \infty. \qquad (17) $$

Dans la suite, on utilise les notations suivantes: $ \lambda_1 > 0 $ est une constante assez grande fix\'ee, $ \lambda \in (0, \lambda_1) $, $ k $ est grand (d\'epend de $ \lambda_1 $), et $ C $ une constante positive ind\'ependante de $ k $ et $ \lambda $ ( mais d\'epend de  $ \lambda_1 $).

\bigskip

Comme,

$$  \Delta v_k^{\lambda} = - \left (\dfrac{\lambda}{|y|} \right )^{n+2} \Delta v_k \left ( y^{\lambda} \right ), $$

On a, en utilisant $ (11) $,

$$  \Delta v_k^{\lambda}(y) + { v_k^{\lambda}(y)}^{(n+2)/(n-2)} = E_1(y) \,\,\, y \in \Sigma_{\lambda}. \qquad (18) $$ 

Avec,

$$ E_1(y) =  - \left (\dfrac{\lambda}{|y|} \right )^{n+2} \left ( \bar b_i (y^{\lambda}) \partial_i v_k(y^{\lambda} ) + \bar d_{ij}(y^{\lambda}) \partial_{ij} v_k(y^{\lambda}) -\bar c(y^{\lambda}) v_k (y^{\lambda}) \right ). \,\,\, (19) $$

Il s'en suit que,

$$ | E_1(y) | \leq C_2 { \lambda }^{n+2} M_k^{-4/(n-2)}  |y|^{-n-2} , \,\,\, y \in \Sigma_{\lambda}   \qquad (20) $$

Soit,

$$ w_{\lambda} = v_k - v_k^{\lambda}. $$

Ici, par soucis de simplification, on omet $ k $ dans la notation de $ w_{\lambda} $.  En utilisant $ (11) $ et $ (18) $,

$$ ( \Delta w_{\lambda} + \bar b_i \partial_i  w_{\lambda} + \bar d_{ij} \partial_{ij} w_{\lambda}) -\bar c w_{\lambda} +\dfrac{n+2}{n-2} \xi^{4/(n-2)} w_{\lambda} = E_{\lambda}, \,\,\, \,\, {\rm dans } \,\, \Sigma_{\lambda}   \qquad (21) $$

avec $ \xi $ reste entre $ v_k $ et $ v_k^{\lambda} $, et,

$$ E_{\lambda} = - \bar b_i \partial_i v_k^{\lambda} - \bar d_{ij} \partial_{ij} v_k^{\lambda} -\bar c v_k^{\lambda} - E_1. \qquad (22)  $$

Dans la suite, on supposera $ n = 3 $. En utilisant $ (16) $, on a,

$$ | \partial_i v_k^{\lambda} | \leq  C \lambda |y|^{-2} | \partial_{ij} v_k^{\lambda} |  \leq   C \lambda |y|^{-3}, \,\,\, \,\, {\rm dans } \,\, \Sigma_{\lambda}   \qquad (23) $$

On utilise $ (14) $ et $ (23) $, on d\'eduit de $ (22) $ le lemme suivant:

{\it Lemme 2.1.}   Il existe une constante $ C_3 = C_3(\lambda_1) $ telle que,

$$   | E_{\lambda}(y) | \leq  C_3 M_k^{-4} \lambda |y|^{-1}, \,\,\, {\rm dans } \,\, \Sigma_{\lambda}.   \qquad (24) $$

Soit,

$$ h_{\lambda}(y) = - C_3 M_k^{-4} \lambda ( |y| - \lambda), \,\,\,\,\, {\rm dans } \,\,  \Sigma_{\lambda}. $$

{\it Lemme 2.2.} On a,

$$ w_{\lambda} + h_{\lambda}  \geq 0, \,\,\, \in \Sigma_{\lambda} \,\,\, \forall  \,\, 0 < \lambda \leq \lambda_1.   \qquad (25)   $$

{\it Preuve du lemme 2.2.} Il y a deux \'etapes dans la preuve de ce lemme:

\bigskip

{\it Etape 1.} Il existe $ {\lambda}_{0, k} > 0 $ telle que $ (25) $ est vraie pour $ 0 <  \lambda \leq {\lambda}_{0, k} $.

\bigskip

Pour le voir, on \'ecrit,

$$ w_{\lambda}(y) = v_k(y) -v_k^{\lambda}(y) = \dfrac{1}{\sqrt |y|} \left ( \sqrt {|y|} v_k(y) - \sqrt {|y^{\lambda }| } v_k(y^{\lambda}) \right ). $$

En coordonn\'ees polaires, on a,

$$ f(r, \theta) = \sqrt {r} v_k(r, \theta ) . $$

En utilisant $ (16) $, il existe $ r_0 > 0 $ et $ C > 0 $ ind\'ependant de $ k $ tels que,

$$ \partial_r f(r, \theta ) >  C r^{-1/2}  \quad {\rm pour } \quad  0 < r < r_0 . $$

Par cons\'equent, pour $ 0 < \lambda < |y| < r_0 $, on a,

$$  w_{\lambda }(y) + h_{\lambda}(y) = v_k(y) - v_k^{\lambda}(y) + h_{\lambda}(y) $$

$$ > \dfrac {1}{\sqrt {r_0}} C r_0^{-1/2} ( |y| - |y^{\lambda}|) + h_{\lambda} $$

$$ > ( \dfrac{C}{\sqrt r_0} - C_3 \lambda M_k^{-4} ) ( |y| - \lambda ) \quad {\rm car } \quad  | y | - | y^{\lambda } | >  |y| - \lambda  $$

$$ > 0. \qquad ( 26) $$

Comme,

$$ | h_{ \lambda } (y) | + v_k^{\lambda } (y) \leq C(k, r_0) \lambda,  \quad  r_0 \leq  |y| \leq \epsilon_k M_k^{-2}, $$

On peut choisir $ \lambda_{0, k} \in (0, r_0) $ assez petit ( d\'ependant de $ k $ et de $ r_0 $) tel que pour tout $ 0 < \lambda < \lambda_{0, k} $ on a,

$$ w_{\lambda } ( y ) + h_{ \lambda } (y) \geq   \min_{|y| \leq \epsilon_k M_k^{-2} } v_k(y) -C(k, r_0) \lambda_{0, k} > 0 , \quad \forall \,\, r_0 \leq |y| \leq \epsilon_k M_k^{2}. $$

La derni\` ere assertion et $ (26) $ donnent l'\'etape 1.

\bigskip

Soit,

$$ { \bar \lambda }^k = \sup \{ 0 <  \lambda \leq  \lambda_1, \,\,\, w_{\mu} + h_{\mu} \geq 0  \,\, {\rm dans } \,\, \Sigma_ {\mu}, \,\, {\rm pour \,\, tout }  \,\, 0 < \mu \leq \lambda \} .  \,\,\, (27) $$

\bigskip

{\it Etape 2.} $ {\bar \lambda }^k = \lambda_1 $, i.e. $ (25) $ est vraie.

\bigskip

Pour prouver cela, on a besoin d'estimer,

$$ ( \Delta + \bar b_i \partial_i  + \bar d_{ij} \partial_{ij} + \dfrac{n+2}{n-2} \xi^{4/(n-2)} - \bar c )(w_{\lambda } +  h_{\lambda} ) \leq  0, \,\, {\rm dans } \,\, \Sigma_ {\lambda}, \quad (28) $$  

i.e., d'apr\'es $ (21) $, il suffit de v\'erifier,

$$  \Delta  h_{\lambda }+ \bar b_i \partial_i  h_{\lambda } + \bar d_{ij} \partial_{ij} h_{\lambda } + E_{\lambda } + ( 5 \xi^{4} - \bar c ) h_{\lambda } \leq  0, \,\, {\rm dans } \,\, \Sigma_ {\lambda}, \quad (29) $$  

Comme $  h_{\lambda } < 0 $ dans $ \Sigma_ {\lambda} $,

$$ 5 \xi^{4} h_{\lambda }  < 0 \,\, {\rm dans } \,\, \Sigma_ {\lambda}. $$

Le terme dominant dans $ (29) $ est,

$$  \Delta  h_{\lambda }(y) = -2C_3 \lambda M_k^{-4} |y|^{-1} . $$

Les termes restant sont d'ordre tr\'es grand. En effet,

$$  | \partial_i  h_{\lambda } | \leq C \lambda M_k^{-4},  \,\,\,  |   \partial_{ij} h_{\lambda } | \leq C \lambda M_k^{-4} |y|^{-1}, $$

et, en utilisant $ (14) $,

$$ | \bar b_i(y) \partial_i  h_{\lambda } | +  | \bar d_{ij}(y) \partial_{ij} h_{\lambda }  | +  | \bar c h_{\lambda }  | \leq C \lambda M_k^{-8} |y| \leq C \lambda \epsilon_k^2 M_k^{-4} |y|^{-1} $$

$$ \leq C_3 \lambda M_k^{-4} |y|^{-1}, \,\, {\rm dans } \,\, \Sigma_ {\lambda}. $$ 

Alors, en usant de $ (24) $ et les estimations pr\'ec\'edentes, on obtient,

$$  \Delta  h_{\lambda }+ \bar b_i \partial_i  h_{\lambda } + \bar d_{ij} \partial_{ij} h_{\lambda } + E_{\lambda } + ( 5 \xi^{4} - \bar c ) h_{\lambda } \leq  $$

$$  \leq  \Delta  h_{\lambda } + C_3 \lambda M_k^{-4} |y|^{-1} +  |E_{\lambda } |  \leq $$

$$  \leq - C_3 \lambda M_k^{-4} |y|^{-1} +  |E_{\lambda } |  \leq 0 \,\, {\rm dans } \,\, \Sigma_ {\lambda}.  $$

On voit qu'\`a partir de $ (16) $ et les d\'efinitions de $ v_k^{\lambda }  $ et $ h_{\lambda} $ que,

$$  | v_k^{ { \bar \lambda }^k  }(y) | + | h_{\bar \lambda^k} (y)| \leq \dfrac{C(\lambda_1)}{|y|},  \,\,\, \forall \,\,  |y| = \epsilon_k M_k^2. $$

Alors, en utilisant la condition au bord $ (17) $,

$$ (w_{{\bar \lambda }^k} + h_{{\bar \lambda }^k})(y) > 0 \,\,\, \forall \,\, |y| = \epsilon_k M_k^2, $$

Comme $ w_{{\bar \lambda }^k} + h_{{\bar \lambda }^k} $ est positive ou nulle et v\'erifie $ (29) $ avec $ \lambda = {\bar \lambda }^k $, on peut appliquer le principe du maximum fort et le lemme de Hopf pour obtenir,

$$ w_{{\bar \lambda}^k} + h_{{\bar \lambda}^k} > 0 \,\,\, {\rm dans} \,\,\, \Sigma_{{\bar \lambda }^k}, $$

et,

$$ \partial_{\nu} (w_{{\bar \lambda}^k} + h_{{\bar \lambda }^k} ) > 0 \,\,\, {\rm sur } \,\,\, \partial B(0, {\bar \lambda }^k ), $$

o\`u $ \partial_{\nu} $ d\'enote la d\'eriv\'ee selon la normale exterieure.

\bigskip

D'apr\'es les trois estimations pr\'ec\'edentes, on a, $ {\bar \lambda }^k = \lambda_1 $ et l'\'etape 2 est prouv\'ee. Le lemme 2.2 est prouv\'e.

\bigskip

En se donnant $ \lambda > 0 $, comme $ (v_k) $ converge vers $ U $ (apr\'es passage aux une sous-suites) et $ h_{\lambda} $ converge vers $ 0 $ sur tout compact de $ {\mathbb R}^n $, on a, en faisant tendre $ k $ vers l'infini dans $ (25) $,

$$ U(y) \geq U^{\lambda }(y), \,\,\, {\rm pour \,\, tout } \,\, |y| \geq \lambda, \,\,\, 0 < \lambda < \lambda_1. $$

Comme $ \lambda_1 > 0 $ est arbitaire et le fait qu'on peut appliquer le m\^eme argument pour comparer $ v_k $ et $ v_k^{\lambda, x} $, on a,

$$ U(y) \geq U^{\lambda, x}(y), \,\,\, {\rm pour\,\, tout} \,\, |y-x| \geq \lambda > 0. $$

Ceci implique, par un lemme de calcul donn\'e dans un article pr\'ec\'edent, qe $ U $ est constante, ceci est une contradiction de $ (15) $.

\bigskip

\underbar{\it Cas $ n= 4 $ :}

\bigskip

La preuve suit le m\^eme type d'arguments et \'etapes que pour le cas $ n = 3 $. Le premier changement est en relation avec la courbure scalaire, $  R $ et la courbure de Ricci, $ Ricci $. Les \'estimations se font autour d'un pount $ x $ qu'on supposera $ 0 $ ici, on fera un changement de m\'etrique conforme de telle mani\`ere que $ R(0) = Ricci(0) = 0 $. On a,

$$ g = det(g_{ij}) = 1-\dfrac{1}{3} R_{ij}z^iz^j + O(r^3), \qquad (30) $$

$$ g_{pq}(z) = \delta_{pq} + \dfrac{1}{3} R_{pijq}z^iz^j + O(r^3), \qquad (31) $$

$$ \Delta_g u = \dfrac{1}{\sqrt g} \partial_i \left (\sqrt g g^{ij} \partial_j u \right ) = - \Delta u + b_i \partial_i u + d_{ij} \partial_{ij} u $$

o\`u,

$$ b_j = \dfrac{1}{2g} \partial_i g g^{ij} + \partial_i g^{ij}, \,\,\, d_{ij} = g^{ij} - \delta^{ij}. \qquad (32) $$

Comme $ R_{jp} = 0 $, on a, en utilisant $ g $ et $ g_{ij} $, que,

$$ g^{ij} = \delta_{ij} -\dfrac{1}{3}R_{ipqj} z^pz^q + O(r^3), $$

$$ \partial_i g^{ij} = -\dfrac{1}{3} R_{ipij}z^p - \dfrac{1}{3}R_{iiqj}z^q + O(r^2) = O(r^2) $$

et,

$$ \partial_i g = -\dfrac{2}{3} R_{ip} z^p + O(r^2) = O(r^2). $$

Des expressions pr\'ec\'edents et $ R(0) = 0 $, on obtient,

$$ b_i = O(r^2), \,\,\, d_{ij} = -\dfrac{1}{3} R_{ipqj} z^p z^q + O(r^3), \,\,\, R =O(r). \qquad (33) $$

On garde les m\^emes notations que pour le cas $ n=3 $, pour $ \bar b_i, \bar d_{ij} $ et $ \bar c $. On utilise les estimations pr\'ec\'edentes pour avoir,

$$ |\bar b_i(y)|\leq CM_k^{-2/(n-2)}|y|^2,\,\,\, |\bar c|\leq CM_k^{-6/(n-2)}|y| \qquad (34) $$

et,

$$ \bar d_{ij}(y) =-\dfrac{1}{3} M_k^{-4/(n-2)}R_{ipqj}y^p y^q + O(1)M_k^{-6/(n-2)} |y|^3, \qquad (35) $$

En utilisant la fonction de Green du Laplacien conforme sur $ B(x_k, 2\epsilon) $ avec $ \epsilon > 0 $ assez petit pour que le principe du maximum sout vrai, on peut affirmer qu'il existe une constante positive $ C_4 > 0 $ ind\'ependante de $ k $ telle que pour $ k $ assez large,

$$ v_k(y) \geq C_4 |y|^{2-n}, \,\,\, {\rm pour }\,\,\, 1 \leq |y| \leq 2 \epsilon_k M_k^{2/(n-2)}. \qquad (36) $$

Dans la suite, on garde les m\^emes notations que pour le cas $ n = 3 $, $ v_k, v_k^{\lambda}, w_{\lambda}, E_{\lambda} $ et $ \lambda_1 > 0 $.

\bigskip

{\it Lemme 3.1.} Il existe une constante positive $ C_3 = C_3(\lambda_1) $ telle que,

$$ |E_{\lambda}(y)| \leq C_3 \lambda^2 M_k^{-2} |y|^{-3} + C_3 \lambda^2 M_k^{-3} |y|^{-1}, \,\,\, {\rm dans } \,\,\, \Sigma_{\lambda}. \qquad (37) $$

{\it Preuve du lemme 3.1}. Dans $ \Sigma_{\lambda} $ on a,

$$ \partial_i v_k^{\lambda} = (2-n) \dfrac{\lambda^{n-2}}{|y|^n} y^i v_k(y^{\lambda}) + O(1) \lambda^{n-2} |y|^{1-n}, \qquad (38) $$

et,

$$ \partial_{ij} v_k^{\lambda}(y) = (n-2) (n\dfrac{\lambda^{n-2}}{|y|^{n+2}} y^iy^j -\dfrac{\lambda^{n-2}}{|y|^n}\delta_{ij}) v_k(y^{\lambda}) + O(1) \lambda^n |y|^{-1-n}, \qquad (39) $$

avec $ O(1) $ d\'epend $ \lambda_1 $ mais ind\'ependant de $ k $.

\bigskip

Une cons\'equence de $ (22),(20), (34), (35), (38) $  et $ (39) $ que,

$$ E_{\lambda}(y) = O(1) M_k^{-6/(n-2)} |y|^2 (\lambda^{n-2} |y|^{1-n}) + $$

$$ + \left [ \dfrac{1}{3} M_k^{-4/(n-2)} R_{ipqj} y^p y^q + O(1) M_k^{-6/(n-2)} |y|^3 \right ] \times $$

$$ \times  \left [ (n-2) \left [n\dfrac{\lambda^{n-2}}{|y|^{n+2}} y^iy^j - \dfrac{\lambda^{n-2}}{|y|^n} \delta_{ij} \right ] v_k(y^{\lambda}) + O(1) \lambda^n |y|^{-n-1} \right ] + $$

$$ + O(1) M_k^{-6/(n-2)} |y| \left (\dfrac{\lambda}{|y|} \right )^{n-2} v_k(y^{\lambda}) + O(1) M_k^{-4/(n-2)} \dfrac{\lambda^{n+2}}{|y|^{n+2}}. $$

On utilise les propri\'et\'es d'antisymmetrie de $ R_{ipqj} $ et le fait que $ R_{pq} = 0 $, on a, $ R_{ipqj}y^py^q \delta_{ij} = -R_{pq}y^py^q = 0 $ et $ R_{ipqj} y^p y^q y^i y^j = 0 $. De ce qui pr\'ec\`ede on d\'eduit l'estimation du lemme.

\bigskip

Pour $ \alpha < 4 $ et $ \alpha \not = 2 $, soit,

$$ f_{\alpha}(z) =-\dfrac{1}{(4-\alpha)(2-\alpha)}[|z|^{2-\alpha} -1] - \dfrac{1}{2(4-\alpha)}[|z|^{-2} -1], \,\,\, |z|\geq 1. $$

Alors,

$$ f_{\alpha}(z) = 0, \,\,\, |z| =1, $$

et, pour $ |z| \geq 1 $,

$$ \Delta f_{\alpha}(z) =-|z|^{-\alpha}, $$

$$ f_{\alpha}(z) \leq 0, \,\,\, |f_{\alpha}(z)| \leq C(\alpha) |z|^{\max \{0, 2-\alpha \}}, \qquad (40) $$

$$ |\nabla f_{\alpha}(z)|\leq C(\alpha)(|z|^{1-\alpha} + |z|^{-3}) \leq C(\alpha) |z|^{1-\alpha}, $$

et,

$$ |\nabla^2 f_{\alpha}(z)|\leq C(\alpha)(|z|^{-\alpha} + |z|^{-4}) \leq C(\alpha)|z|^{-\alpha}, \qquad (41) $$

On d\'efinit,

$$ h_{\lambda}(y) = 2C_3 \lambda M_k^{-2} f_3(\dfrac{y}{\lambda}) + 2C_3 \lambda^3 M_k^{-3} f_1(\dfrac{y}{\lambda}), \,\,\, y \in \Sigma_{\lambda}. $$

Alors,

$$ h_{\lambda}(y) \leq 0 \,\,\, {\rm dans} \,\,\, \Sigma_{\lambda}, \qquad (42) $$

et,

$$ -\Delta h_{\lambda}(y) = -2 C_3 \lambda^2 M_k^{-2} |y|^{-3} -2C_3  \lambda^2 M_k^{-3} |y|^{-1}, \,\,\, y \in \Sigma_{\lambda}, \qquad (43) $$

\bigskip

{\it Lemme 3.2.} On a,

$$ w_{\lambda} + h_{\lambda} > 0 \,\,\, {\rm dans} \,\,\, \Sigma_{\lambda}. \qquad (44) $$

{\it Etape 1.} Il existe $ \lambda_0 > 0 $ ind\'ependant de $ k $ tel que $ (44) $ est vraie pour tout $ 0 < \lambda < \lambda_0 $,

\bigskip

Pour le voir, on \'ecrit,

$$ w_{\lambda}(y) = v_k(y) -v_k^{\lambda}(y) = |y|^{-1} [ |y|v_k(y) -|y^{\lambda}|v_k(y^{\lambda}) ]. $$

On pose, en coordonn\'ees polaires,

$$ f(r, \theta) = rv_k(r, \theta). $$

A l'aide $ (16) $, il existe $ r_0 > 0 $ et $ C > 0 $ ind\'ependants de $ k $ tels que,

$$ \partial_r f (r, \theta) > C > 0, \,\,\, {\rm pour} \,\,\, 0 < r < r_0. $$

Par cons\'equent,

$$ w_{\lambda}(y) \geq C^{-1} |y|^{-1} |y-y^{\lambda}| \geq \dfrac{1}{C r_0} (|y|-\lambda),\,\,\, {\rm pour } \,\,\, 0 < \lambda < |y| < r_0. $$

D'autre part, pour $ y \in \Sigma_{\lambda} $,

$$ |h_{\lambda}(y)| \leq C \lambda M_k^{-2}|f_3(\dfrac{y}{\lambda})| +C\lambda^3M_k^{-3}|f_1(\dfrac{y}{\lambda})| $$

$$ \leq C \lambda M_k^{-2}| |\dfrac{y}{\lambda}|^{-2} -1 | + C\lambda^3 M_k^{-3} | |\dfrac{y}{\lambda} | -1| \leq CM_k^{-2}(|y|-\lambda). $$

Il s'en suit que,

$$ w_{\lambda} + h_{\lambda} \geq ( \dfrac{1}{C r_0}- \dfrac{C}{M_k^2} )(|y|-\lambda) > 0 \,\,\, 0 < \lambda < |y| < r_0. \qquad (45) $$

Pour $ \epsilon_k M_k \geq |y| \geq r_0 $, on a, en utilisant $ (40) $, $ (16) $ et $ (36) $, que, $ |h_{\lambda}(y) | \leq CM_k^{-2} \leq C \epsilon_k^2 |y|^{-2} < \dfrac{1}{2} v_k(y) $. Ainsi, pour $ \epsilon_k M_k \geq |y| \geq r_0 $,

$$ v_k(y) - v_k^{\lambda}(y) + h_{\lambda}(y) > \dfrac{1}{2} v_k(y) -v_k^{\lambda}(y) > \dfrac{1}{2} v_k(y) - \left (\dfrac{\lambda}{|y|} \right )^2 \max_{B(0,r_0)} v_k. \,\,\, (46) $$

De $ (16) $, $ (36) $, $ (45) $ et $ (46) $, on d\'eduit qu'on peut choisir $ \lambda_0 > 0 $ ind\'ependant de $ k $ tel que $ (44) $ soit vraie pour $ 0 < \lambda < \lambda_0 $.

\bigskip

On d\'efinit $ {\bar \lambda }^k $ comme dans $ (27) $.

\bigskip

{\it Etape 2.} $ {\bar \lambda}^k = \lambda_1 $, i.e. $ (44) $ est vraie.

\bigskip

D'apr\'es l'\'etape 1, on sait que, $ \lambda_0 \leq {\bar \lambda}^k \leq \lambda_1 $. On veut prouver que, pour $ \lambda_0 \leq \lambda \leq \lambda_k \leq \lambda_1 $, que,

$$ (\Delta + \bar b_i \partial_i + \bar d_{ij} \partial_{ij} + 3\xi^2 - \bar c) (w_{\lambda} + h_{\lambda}) \leq 0, \,\,\, {\rm dans} \,\,\, \Sigma_{\lambda}. \qquad (47) $$

En s'aidant de $ (21) $, cela revient \`a prouver que,

$$ (\Delta + \bar b_i \partial_i + \bar d_{ij} \partial_{ij} + 3 \xi^2 - \bar c) h_{\lambda} + E_{\lambda} \leq 0 \,\,\, {\rm dans} \,\,\, \Sigma_{\lambda}, \qquad (48) $$

pour $ \lambda_0 \leq \lambda \leq {\bar \lambda}^k \leq \lambda_1 $.

\bigskip

Dans ce qui suit, on suppose que, $ \lambda_0 \leq \lambda \leq {\bar \lambda}^k \leq \lambda_1 $. Rappelons que $ \lambda_0 > 0 $ est ind\'ependant de $ k $, un fait qui sera utilis\'e ci-apr\'es. En utilisant $ (34) $, $ (35) $ et $ (40) $, on a, pour $ y \in \Sigma_{\lambda} $, que,

$$ |\bar b_i||\partial_i h_{\lambda}(y)| \leq CM_k^{-3}|y|^{-3}|y|^2 [M_k^{-2} |\nabla f_3 \left ( \dfrac{y}{\lambda} \right) | + M_k^{-3} |\nabla f_1 \left ( \dfrac{y}{\lambda} \right )|] $$

$$ \leq C M_k^{-5} + CM_k^{-6}|y|^2 \leq C \epsilon_k^3 M_k^{-2} |y|^{-3} + C\epsilon_k^3 M_k^{-3} |y|^{-1}, $$

$$ |\bar c(y)h_{\lambda}(y)| \leq CM_k^{-5} |y| \leq C \epsilon_k M_k^{-3} |y|^{-1}, $$

et, avec $ (35) $,

$$ |\bar d_{ij}(y)||\partial_{ij}h_{\lambda}(y)| \leq CM_k^{-2} |y|^2 |\nabla^2 h_{\lambda}(y)| $$

$$ \leq C { \epsilon_k}^2  \left [  M_k^{-2} \left | \nabla^2 f_3 \left ( \dfrac{y}{\lambda}  \right ) \right | + M_k^{-3} \left |  \nabla^2 f_1  \left (  \dfrac{y}{\lambda} \right )  \right |   \right ]  $$

$$ \leq C \epsilon_k^2 M_k^{-2} |y|^{-3} + C \epsilon_k^2 M_k^{-3} |y|^{-1}, $$

En combinant les estimations pr\'ec\'edentes, $ (43) $ et le lemme 3.1, on a, dans $ \Sigma_{\lambda} $, que,

$$ (\Delta + \bar b_i \partial_i +\bar d_{ij} \partial_{ij} - \bar c) h_{\lambda}(y) $$

$$ \leq -2C_3\lambda^2 M_k^{-2}|y|^{-3} -2C_3\lambda^2M_k^{-3}|y|^{-1} + C\epsilon_k \left ( M_k^{-2}|y|^{-3} + M_k^{-3}|y|^{-1} \right ) $$

$$ \leq - C_3 \lambda^2 M_k^{-2} |y|^{-3} - C_3 \lambda^2 M_k^{-3} |y|^{-1} \leq -|E_{\lambda}|. $$

L'estimation $ (48) $ et $ (47) $ sont les cons\'equences de ce qui pr\'ec\`ede et du fait que $ h_{\lambda} \leq 0 $ dans $ \Sigma_{\lambda} $. Avec $ (16) $ et la forme explicite de $ h_{\lambda} $, on a,

$$ |v_k^{{\bar \lambda}^k}(y)| + |h_{{\bar \lambda}^k}(y)| \leq \dfrac{C}{|y|^2}, \qquad \forall \,\, |y|=\epsilon_kM_k. $$

et avec la condition au bord $ (17) $,

$$ (w_{{\bar \lambda}^k} + h_{{\bar \lambda}^k})(y), \qquad \forall \,\, |y|=\epsilon_kM_k. $$

L'\'etape 2 est la cons\'equence des estimations pr\'ec\'edentes, comme pour le lemme 2.2. Le reste de la preuve du lemme 3.2 est identique \`a celle du lemme 2.2. La preuve du cas $ n = 4 $ est alors identique \`a celle du cas $ n = 3 $.

\underbar {\bf Autour de l'\'equation de la courbure sclaire prescrite.}

\bigskip

{\bf Probl\`eme 1}. Etant donn\'ee sur un ouvert $ \Omega $ de $ {\mathbb R}^n $
l'\'equation suivante:

$$ \Delta u =Vu^{q-1} \,\, {\rm et } \,\, u>0 \,\,\, {\rm dans } \,\,\,
{\Omega },
\qquad  (E_1)$$

o\`u  $ 2 < q \leq N=\dfrac{2n}{n-2} $ et $ V $ est une fonction qui v\'erifie,
pour trois r\'eels positifs $ a, b, A $, les conditions suivantes:

\bigskip

$  \qquad \qquad \qquad \qquad 0< a \leq V(x) \leq b, \qquad \forall
\, x \in \Omega  \,\,\,\,  {\rm et }$

\bigskip

$ \qquad \qquad \qquad \qquad |V(x)-V(y)|\leq A || x-y || \qquad \forall \, x,y\in
\Omega .$

\bigskip

On se pose la question de savoir si, pour chaque compact $ K $ de $
\Omega $, il existe une constante $
c>0 $, ne d{\'e}pendant que de $  a, b, A, K, \Omega  $, telle que pour toute fonction $ u $ solution
de $ (E_1) $, on ait:

\bigskip

$$ \sup_{{}_K} u \times \inf_{{}_{\Omega }} u\leq c. $$

\underbar{\it Estimations Asymptotiques}

\bigskip

{\bf Th\'eor\`eme 1}. Consid\'erons deux  suites $ (u_{\epsilon_i}), (V_{\epsilon_i})  $ de fonctions relatives
au probl\`eme concernant $ (E_1) $ avec $ q_{\epsilon_i} =N-\epsilon_i
\to N $, alors on a:

\smallskip

Pour tout compact $ K $ de $ \Omega $, il existe une constante $ c>0 $ ne
d\'ependant que de $ a, b, A, K, \Omega $ telle que pour tout $
u_{\epsilon_i} $ :

$$ {\epsilon_i}^{(n-2)/2} (\sup_{{}_K} u_{\epsilon_i})^{1/4} \times
\inf_{{}_{\Omega}} u_{\epsilon_i} \leq c. $$

\underbar {\bf Preuve.}

\bigskip

{\bf \underbar { D{\'e}monstration du Th{\'e}or{\`e}me 1}}

\smallskip 

\underbar { $ q_i=N-\epsilon_i \to N $ }

Par souci de compr{\'e}hension, nous allons  d{\'e}tailler cette
partie. On aura {\`a}
utiliser la technique "moving plane " qui utilise essentiellement le
principe du maximum.

\smallskip

On suppose, pour simplifier, que $ \Omega=B_2(0) $, et on raisonne par
l'absurde, en essayant de d{\'e}montrer qu'il existe pour un certain $ \beta
\in ]0,\dfrac{1}{3} [ $, une constante $ c $ ne d{\'e}pendant que de $
a, b, A, \beta $ et un r{\'e}el $ R\in]0,1[ $ tels que pour tout $
u_{\epsilon} >0 $ solution de $ (E_1) $ avec $ V=V_{\epsilon} $ v{\'e}rifie:

$$  {\epsilon }^{[2/(n-2)-\epsilon/2]^{-1}}{(\sup_{{}_{B_R(0)}}
  u_{\epsilon})}^{\beta }  \times \inf_{{}_{B_1(0)}} u_{\epsilon} \leq
\dfrac{c}{R^{4/(N-\epsilon-2)}} \,\, \forall\,\, \epsilon>0 .$$

Le fait de prendre l'inf sur la boule unit{\'e}  dans  la boule de
rayon 2 est li{\'e} aux calculs qui vont suivre, car nous  serons oblig{\'e}s
d'effectuer des translations et il nous faut une marge de manoeuvre.

\smallskip
 
On a remplac{\'e} l'exposant $ \dfrac{1}{4} $ du sup par $ \beta $, on verra que le
r{\'e}sultat de cette deuxi{\`e}me partie du th{\'e}or{\`e}me est valable pour tout $
0< \beta <\dfrac{1}{3} $.

\smallskip

Supposons donc que pour tout $ c>0 $ et $ R \in]0,1[ $, il existe $
V_{\epsilon} $ et $ u_{\epsilon} $  v{\'e}rifiant:

$$  {\epsilon}^{[2/(n-2)-\epsilon/2]^{-1}} {(\sup_{{}_{B_R(0)}}
  u_{\epsilon})}^{\beta } \times \inf_{{}_{B_1(0)}} u_{\epsilon} \geq
\dfrac{c}{R^{4/(N-\epsilon-2)}} \,\,\,\, {\rm et } $$

$$ \Delta u_{\epsilon }=V_{\epsilon} {u_{\epsilon}}^{N-\epsilon-1}. $$

On choisira : $ R=R_i \rightarrow 0  $ et $ c=c_i\rightarrow
+\infty $. Notre hypoth{\`e}se est: il existe deux suite $   \{
u_{\epsilon_i } \}  $ et $ \{ V_{\epsilon_i}\}  $ not{\'e}es, pour
simplifier l'{\'e}criture $ \{u_i\} $ et $ \{V_i\} $, telles que pour
tout $ i \in {\mathbb N} $:

$$ \Delta u_i=V_{i}
{u_i}^{N-\epsilon_i-1},  $$

$$ {\epsilon_i}^{[2/(n-2)-\epsilon_i/2]^{-1}} {(\sup_{{}_{B_{R_i}(0)}}
  u_{i})}^{\beta } \times \inf_{{}_{B_1(0)}} u_{i} \geq
\dfrac{c_i}{{R_i}^{4/(N-\epsilon_i-2)}}. $$ 

D'une mani{\`e}re {\'e}vidente (on suppose $ \sup u_i >1 $ et $ \epsilon_i \to 0 $ ),

$$ {(\sup_{{}_{B_{R_i(0)}}} u_i)}^{4/3}\geq {(\sup_{{}_{B_{R_i(0)}}}
  u_i)}^{\beta }\times  \inf_{{}_{B_1(0)}} u_i \geq {\epsilon_i}^{{[2/(n-2)
    -\epsilon_i]}^{-1}} {(\sup_{{}_{B_{R_i(0)}}}
  u_i)}^{\beta }\times  \inf_{{}_{B_1(0)}} u_i, $$

donc:

$$ {(\sup_{{}_{B_{R_i(0)}}} u_i)}^{1+\beta } \geq
\dfrac{c_i}{{R_i}^{4/(N-\epsilon_i-2)}} \to +\infty .$$

En particulier,

$$ (\sup_{B_{R_i(0)}} u_i)\times
{R_i}^{2/(N-\epsilon_i-2)} \geq \sqrt {c_i} \to +\infty .$$

Consid{\'e}rons alors:

$ s_i(x)=u_i(x){(R_i-|x-x_i|)}^{2/ (N-\epsilon_i-2)} $ avec $ u_i(x_i)=\max_{\bar B_{R_i(0)}} u_i \, .$

\smallskip

Soit  $ a_i $, tel que : 

\smallskip

$ \qquad \qquad s_i(a_i)=\max_{B_{R_i(x_i)}} s_i =  u_i(a_i)
{(R_i-|a_i-x_i|)}^{2/ (N-\epsilon-2)} $.

\smallskip

Nous avons : 

\smallskip

$  \qquad  \qquad s_i(a_i) \geq s_i(x_i) =
u_i(x_i){R_i}^{2/ (N-\epsilon_i-
2)  }\geq \sqrt{c_i}  \,\,\,\, (*) $ avec $ c_i \to +\infty $.

\smallskip

Posons : 

 $$ l_i=(R_i-|a_i-x_i|),\,\, {\rm et } \,\, 
L_i=\dfrac{l_i}{ \root 4 \of {c_i} }{[u(a_i)]}^{(N-\epsilon_i-2)/2} $$

et remarquons que 
  
$$ 0 < l_i\leq R_i\longrightarrow 0 ,\,\,\Rightarrow u_i(a_i)\to
+\infty \,, $$

d'apr{\`e}s $ (*) $ :

\smallskip

$\qquad \qquad \qquad \qquad \qquad  L_i \to +\infty  $

\smallskip

Posons, lorsque  $ |y|\leq L_i $,

$$ v_i(y)=\dfrac{1}{u_i(a_i)} u_i\{ y
[{u_i(a_i)}^{(2+\epsilon_i-N)/2}]+a_i \}       $$

et v{\'e}rifions que si, $ |y|\leq L_i $, alors $ x= y
[{u_i(a_i)}^{(2+\epsilon_i-N)/2}]+a_i \in B_{R_i}(x_i) .   $

\smallskip

Gr{\^a}ce {\`a} l'in{\'e}galit{\'e} triangulaire  :

$$ R_i-|x-x_i|=R_i-|a_i-x_i+{[u_i(a_i)]}^{(2+\epsilon_i-N)/2}y|\geq 
R_i-|a_i-x_i|-|{[u_i(a_i)]}^{(2+\epsilon_i-N)/2}y| $$

et donc: 

\smallskip

$  \qquad R_i-|x-x_i|\geq l_i-l_i\dfrac{1}{\root 4 \of {c_i} }=  l_i \left
  (1-\dfrac{1}{ { \root 4 \of {c_i} } } \right ) > 0 , \qquad (**) $

\smallskip

$\qquad  |x-x_i|\leq R_i- l_i \left
  (1-\dfrac{1}{ { \root 4 \of {c_i} } } \right ) < R_i ,  $

\smallskip

$ v_i $ est ainsi bien d{\'e}finie et v{\'e}rifie pour tout $ i \in {\mathbb N
  }, $

$$\Delta v_i = V_i\{a_i+ y
[{u_i(a_i)}^{(-N+\epsilon_i+2 )/2}]\} {v_i}^{N-\epsilon_i-1}
  ,\,\, {\rm et } \,\,  v_i(0)=1\,.  $$

D'autre part:

$$ v_i(y)=\dfrac{s_i(x)}{s_i(a_i)} \times  \dfrac {
  {(R_i-|a_i-x_i|)}^{2/(N-\epsilon_i-2  )} }{
  {(R_i-|x-x_i|)}^{2/(N-\epsilon_i-2   )} } \leq \left (
  \dfrac{l_i}{R_i-|x-x_i|} \right )^{2/(N-\epsilon_i-2)} . $$

D'o{\`u} pour tout $ i\in{\mathbb N} $ et $ |y|\leq L_i $:

$$ \qquad \qquad 0< v_i(y)\leq  \left (1-\dfrac{1}{ \root 4
        \of {c_i}}  \right ) ^{-2/(N-\epsilon_i-2)}  . \quad (***) $$

Comme dans l'{\'e}tape 2 de la preuve du th{\'e}or{\`e}me 1, gr{\^a}ce aux th{\'e}or{\`e}mes de Ladyzhenskaya et Ascoli, de la suite de
fonctions $ v_i $  on peut  extraire une sous-suite qui converge
uniform{\'e}ment vers une fonction $ v\geq 0$ et qui v{\'e}rifie :

$$ \Delta v=V(0)v^{N-1},\,\, v(0)=1,\,\, \,\, 0<a\leq V(0)\leq
b<+\infty  .$$

En faisant un changement d'{\'e}chelle on peut se ramener au cas: $
V(0)=n(n-2) $.

\smallskip

Les solutions positives  de: $\Delta v=n(n-2)v^{N-1} $, sur  ${\mathbb
  R}^n $ sont  les fonctions (voir le r{\'e}sultat de
  Caffarelli-Gidas-Spruck [5]) :

$$ v(y)=\dfrac{\mu}{{(\mu^2+|y-x_0|^2)}^{(n-2)/2}} \,\, {\rm avec }
\, \mu \in {\mathbb R}^+ \,\, {\rm et } \,\, x_0 \in {\mathbb R}^n \,. $$

D'apr{\`e}s $ (***) $ et comme $ c_i \to +\infty $,

$$ v(y) \leq 1 \,\,{\rm pour \, tout } \,\, y \in {\mathbb
  R}^n  $$

d'o{\`u},

$$ v(0)=\max_{{\mathbb R}^n} v=1 \,\, {\rm et } \,\,  \nabla v(0)=0 \Rightarrow
x_0=0\,, $$

enfin

$$ v(0)=1 \Rightarrow \mu=1 \,.$$

Remarquons aussi que

$$
{l_i}^{2\beta/(N-\epsilon_i-1)}[u_i(a_i)]^{\beta }=[s_i(a_i)]^{\beta }\geq
{[s_i(x_i)]}^{\beta }=[u_i(x_i){R_i}^{2/(N-\epsilon_i-1)}]^{\beta } .$$

D'o{\`u} d'apr{\`e}s notre hypoth{\`e}se,

$$ {l_i}^{2\beta /(N-\epsilon_i-1)}[u_i(a_i)]^{\beta }\times
\inf_{B_1(0)}  u_i\geq \dfrac{{c_i}^{\beta }}{
{R_i}^{(4-2\beta ) /(N-\epsilon_i-1)}} .$$

Rappelons que, $ u_i(x_i)=\max_{\bar B_r(0)} u_i $. Comme $
l_i,R_i\rightarrow 0 $ et $  0< \beta <\dfrac{1}{3} $, on a:

\smallskip

$  \qquad \qquad \qquad \qquad [u_i(a_i)]^{\beta }\times
\inf_{B_1(0)} u_i \rightarrow +\infty .$

\smallskip

\underbar{ Conclusion de l'Etape 1:}

\smallskip

$ v_i(y)=
v_{\epsilon_i}
(y)=
\dfrac{u_{\epsilon_i
}[a_{\epsilon_i
}+y[u_{\epsilon_i
}(a_{\epsilon_i
})]^{\epsilon_i
/2-2/(n-2)}]}{u_{\epsilon_i}
 (a_{\epsilon_i
})} $ \,\,avec $  a_{\epsilon_i}=a_i $,   v{\'e}rifie:

$$ \Delta v_{\epsilon_i}
=V_{\epsilon_i
}{v_{\epsilon_i
}}^{N-\epsilon_i
-1} \,\, {\rm et } \,\, v_{\epsilon_i}\rightarrow {\left
  (\dfrac{1}{1+|y|^2}\right
  )}^{(n-2)/2} \,\,. $$

Cette convergence {\'e}tant  uniforme  sur  tout compact de  $ {\mathbb R}^n
$:

$$ {[u_i(a_i)]}^{\beta } \times \inf_{B_1(0)} u_i \to
+\infty $$
 avec $ \beta < \dfrac{1}{3} $ et $ a_i \to 0 $.

\smallskip

{\bf Etape 2: Passage en coordonn{\'e}es polaires et utilisation de la
  m{\'e}thode " moving plane "} 

\smallskip

\underbar { Lemme  :}

\smallskip

On pose pour $ t\in]-\infty,0],\, \theta \in {\mathbb S}_{n-1} $:

$$ w_i(t,\theta)=e^{(n-2)t/2} u_i(a_i+
e^t\theta) \,\,\, {\rm et } \,\,\, 
V_i(t,\theta)=V_{\epsilon_i}
(a_i+
e^t\theta). $$

Et on  consid{\`e}re l'op{\'e}rateur suivant :

$$ L=\partial_{tt}-\Delta_{\sigma}-\dfrac{(n-2)^2}{4},\,\,\, {\rm sur
  } \,\,
]-\infty,0]\times {\mathbb S}_{n-1} , $$

avec $ \Delta_{\sigma} $ l'op{\'e}rateur de Baltrami-Laplace sur la sph{\'e}re  ${\mathbb
    S}_{n-1} $
 
alors :

$$ -Lw_i=e^{[(n-2)\epsilon_i t]/2} V_i
 {w_i}^{N-\epsilon_i-1},\,\,
{\rm pour \, tout }\,  i .$$

\underbar { D{\'e}monstration du lemme  :}

$$ \partial_t w_i=\dfrac{(n-2)}{2}
e^{(n-2)t/2} u_i(a_i+
e^t\theta)+e^{nt/2}\partial_ru_i(a_i+e^t\theta), $$

donc:

$$  \partial_{tt}
w_i=\dfrac{(n-2)^2}{4}w_i+e^{(n+2)t/2}[\partial_{rr}u_i(a_i+
e^t\theta)+\dfrac{(n-1)
 }{e^t}\partial_ru_i(a_i+
e^t\theta)]. $$ 

Par d{\'e}finition du $ \Delta_{\sigma} $,

$$ \Delta_{\sigma} w_i= e^{(n-2)t/2}\Delta_{\sigma}
u_i(a_i+
e^t\theta)= e^{(n-2)t/2}\Delta_{\sigma} u_i(a_i+
e^t\theta), $$

d'o{\`u}:

$$ \partial_{tt} w_i-\Delta_{\sigma}
w_i=\dfrac{(n-2)^2}{4}w_i+e^{(n+2)t/2}[\partial_{rr}u_i(a_i+e^t\theta)+\dfrac{(n-1)
 }{e^t}\partial_r(u_i)(a_i+e^t\theta)-\dfrac{1}{e^{2t}}
\Delta_{\sigma} u_i(a_i+e^t\theta)].$$

En remplacant $ e^t $ par $ r > 0 $, sachant que l'expression du Laplacien en coordonn{\'e}es polaires
 est:

$$
-\Delta=\partial_{rr}+\dfrac{(n-1)}{r}\partial_r-\dfrac{1}{r^2}\Delta_{\sigma},
$$ 

en cons{\'e}quence:

$$ -Lw_i=-[\partial_{tt} w_i-\Delta_{\sigma}
w_i-\dfrac{(n-2)^2}{4}w_i]=V_i e ^{(n-2)\epsilon_i
  t/2}{w_i}^{N-\epsilon_i-1}. $$

\smallskip

\underbar{ Etape 2-2: Quelques propri{\'e}t{\'e}s concernant les fonctions $ w_i
  $ }

\smallskip

Posons : $ \eta_i=\dfrac{1}{{u_i(a_i)}^{(N-\epsilon_i-2)/2}} $, alors :

$$ \log \eta_i=-\dfrac{N-\epsilon_i-2}{2}\log u_i(a_i)=-(\dfrac{2}{n-2}-\dfrac{\epsilon_i}{2})\log
  u_i(a_i) .$$ 

\underbar { Lemme :}

On a:

\smallskip

1) \,  $ w_i(\log\eta_i,\theta)-w_i(\log\eta_i+4,\theta)>0, \, i\geq i_0 $

\smallskip

2) \,  $ \forall\, \delta \geq 0,\, \exists \, c(\delta)>0,\,
   i_0=i(\delta)\in {\mathbb N}, $ tels que :

$$   \dfrac{1}{c(\delta)} e^{(n-2)t/2}\times
{u_i(a_i)}^{(n-2)\epsilon_i/4}\leq w_i(t+\log \eta_i,\theta)\leq c(\delta) e^{(n-2)t/2}\times
{u_i(a_i)}^{(n-2)\epsilon_i/4}, $$

 pour tout $  \theta\in {\mathbb
  S}_{n-1},$ $i\geq i_0 $ et $ t \leq \delta $.

\smallskip

\underbar { D{\'e}monstration de 1) :}

\bigskip

En utilisant la d{\'e}finition de $ v_i $, donn{\'e}e dans la conclusion
de l'{\'e}tape 1, nous pouvons {\'e}crire:

 $$ w_i(t+\log \eta_i,\theta)= e^{(n-2)t/2} u_i(a_i+e^t\theta \eta_i)
 {\eta_i}^{(n-2)/2} =e^{(n-2)t/2}{[u_i(a_i)]}^{(n-2) \epsilon_i/4 }\times
   v_i(e^t\theta) $$

Toujours d'apr{\`e}s  l'{\'e}tape 1:

\smallskip

Pour tout $ \beta >0 $, $ z_i(t,\theta)= e^{(n-2)t/2} v_i(e^t\theta) $ converge uniform{\'e}ment
 sur
$]-\infty, \log\beta]\times {\mathbb S}_{n-1}   $
 vers la fonction, $ z(t)=\dfrac{e^{(n-2)t/2}}{(1+e^{2t)^{(n-2)/2}}}={\left
  (\dfrac{e^t}{1+e^{2t}}\right )}^{(n-2)/2} $.

\smallskip

Si on prend $ \log \beta=4 $ et donc, $ t\leq 4 $:

\smallskip

Pour tout $ \epsilon >0 $ il existe un entier  $ i_0 $ tel que   $ i\geq i_0 $
entraine  pour $ t\leq 4 $:  $ z_i(t,\theta)-z(t)< \epsilon .$

\smallskip

En cons{\'e}quence

$$ z_i(0,\theta)-z_i(4,\theta)=[z_i(0,\theta)-z(0)]-[z_i(4,\theta)-z(4)]+
z(0)-z(4) \geq -2\epsilon+z(0)-z(4) $$ 

Sachant que $ w_i $ est obtenue en multipliant $ z_i $ par $
[u_i(a_i)]^{(n-2)\epsilon_i /4} $, comme  $ z(t) $ est maximum en $ t=0 $, pour $ i\geq i_0 $ (on prend
$ 2\epsilon < z(0)-z(4) $), on a:

$$ w_i(\log\eta_i,\theta)-w_i(\log\eta_i+4,\theta)>0 .$$

\underbar { D{\'e}monstration de 2):}

\smallskip

Nous venons de voir qu'en utilisant la convergence uniforme des $ v_i $, on obtient  2).

\smallskip

\bigskip

\underbar {Etape 2-3: Utilisation de la m{\'e}thode " moving plane "}.

\smallskip

On pose lorsque $ \lambda \leq t $ :

$$ t^{\lambda}=2\lambda-t \,\, {\rm et } \,\,
{w_i}^{\lambda}(t,\theta)=w_i(2\lambda-t,\theta). $$

\underbar { Lemme 1:}

\smallskip

Soit $ A_{\lambda} $ la propri{\'e}t{\'e} suivante:

\smallskip

\qquad  $ A_{\lambda}=\{\lambda\leq 0,\,\,\exists \,\,
 (t_{\lambda},\theta_{\lambda })\in
   [\lambda,1/2]\times  {\mathbb S}_{n-1},\,\,   {
   w_i}^{\lambda}(t_{\lambda},\theta_{\lambda})-
   w_i(t_{\lambda},\theta_{\lambda} ) \geq 0\}$

alors:

\smallskip

\qquad $\exists \, \nu \leq 0$, tel que  pour $\lambda \leq \nu $, la
 propri{\'e}t{\'e} 
 $ A_{\lambda} $ n'est pas vraie.

\smallskip

\underbar { Lemme 2 :}

\smallskip

Pour $\lambda \leq 0 $ on a :

$$ {     w_i}^{\lambda} -{     w_i}<0 \Rightarrow -L({     w_i}^{\lambda}
-{     w_i})<0, \,\,\,                                $$

sur $ ]\lambda,t_i]\times {\mathbb S}_{n-1} $ o{\`u} $ t_i=\beta \log
\eta_i+\log \dfrac{(n-2)a}{2A},\,0<\beta<\dfrac{1}{3} $.

\smallskip

3) \underbar{ Un point utile: }

\smallskip

\qquad                ${\xi}_i=$ sup $\{\lambda \leq { \bar\lambda_i}=2+\log
  \eta_i, {     w_i}^{\lambda} -     w_i < 0 $, sur $
  ]\lambda,t_i]\times {\mathbb S}_{n-1}
  \} $
  existe. Avec $ t_i= \beta\log\eta_i+\log \dfrac{(n-2)a}{2A} $ et $ 0<\beta<\dfrac{1}{3} $. 

\smallskip

\underbar { Remarques :}

\smallskip

Dans le Lemme 1, il ne faut pas confondre $ t^{\lambda} $ et $
  t_{\lambda } $, le premier d{\'e}signe le sym{\'e}tris{\'e} de $ t $ alors
  que le second d{\'e}signe un point particulier pour lequel (avec
  $\theta_{\lambda} $), une propri{\'e}t{\'e} donn{\'e}e est
  v{\'e}rifi{\'e}e. 

\smallskip

Le Lemme 1 signifie qu'il existe un rang $ \nu $, petit, tel que pour $
\lambda \leq \nu $, on a: pour tout $ (t,\theta)\in ]\lambda,1/2]\times
{\mathbb S}_{n-1} $ $ w_i^{\lambda}(t,\theta)-w_i(t,\theta)<0 $.

\smallskip
 
Sur les ensembles  consid{\'e}r{\'e}s, le Lemme 2) permettera d'utiliser le
principe du maximum. On trouve  des fonctions $ h $ verifiant :

\smallskip

$  h\leq 0 $ et $ Lh\geq 0 $ avec $
L=\partial_{tt}-\Delta_{\sigma}-\dfrac{{(n-2)}^2}{4} $ o{\'u} $ \Delta_{\sigma} $  est le laplacien sur la sph{\'e}re $ {\mathbb S}_{n-1} $.

\smallskip

Localement  $ L $ s'{\'e}crit  : $ \Sigma_{ij}
a_{ij}\partial_{ij}+\Sigma_j b_j\partial_j -\dfrac{{(n-2)}^2}{4} $, et un
op{\'e}rateur de ce type v{\'e}rifie le principe du maximum de Hopf.

\smallskip

On choisira des domaines particuliers, pour pouvoir utiliser le Lemme
1 convenablement.

\smallskip

On voit aussi que le Lemme 2 est li{\'e} au lemme 1:  pour  $ \lambda
\leq \nu $, la diff{\'e}rence $
{w_i}^{\lambda}-w_i $ est n{\'e}gative. 

\smallskip

On verra l'utilit{\'e} du point 3) apr{\`e}s les d{\'e}monstrations des Lemmes 1 et
2.

\bigskip

\underbar { D{\'e}monstration du Lemme 1:}

\bigskip

D'abord,on fixe l'entier i on cherche le signe de ${\partial}_t{     w_i}$

\smallskip

\qquad ${\partial}_t{
  w_i}(t,\theta)=\dfrac{(n-2)}{
2}e^{(n-2)
t/2}u_i(a_i+e^t\theta)+e^{(n/2)t}{\partial}_ru_i(a_i+e^t\theta) $.

\smallskip

D'o{\`u},

\smallskip

\qquad ${\partial}_t{
  w_i}=  e^{(n-2)t/2}[\dfrac{n-2}{2}u_i(a_i+e^t\theta)+e^{t  }
  {\partial}_ru_i] $.

\smallskip

La fonction $ u_i $ est $ C^1 $, positive  et sous-harmonique,
on en d{\'e}duit
qu'il existe  $ A_i$ tel que $ { \parallel{\partial}_r u_i \parallel}_{\infty}\leq A_i $.

\smallskip

D'autre part, le principe du maximum indique que $ u_i $ atteint son
minimum sur le bord et ainsi,

\smallskip

\qquad $ \dfrac{n-2 }{2}            u_i(a_i+
e^t \theta)\geq \dfrac{n-2
     }{2} \min_{{}_{B_{\sqrt e}(a_i
)}}
 u_i
    \geq   \dfrac{n-2}{2} \min_{B_2(0)} u_i=   \beta_i >0 $.

\smallskip

Finalement,

$${\partial}_t{     w_i}\geq e^{(n-2)t/2} (\beta_i -e^ t   A_i). $$

Pour $
 t< \log \dfrac{\beta_i}{A_i} $, $\beta_i -e^ t   A_i >0 $. Ainsi  $ w_i$
 est strictement croissante sur $]-\infty, \log \dfrac{\beta_i}{A_i}]$
 uniform{\'e}ment en
$\theta \in {\mathbb S}_{n-1} $.

\smallskip

Supposons que Lemme 1 ne soit pas vrai :

\smallskip

Il existe une famille $\{\lambda_k\} $, telle que $\lambda_k
\to -\infty $ , des r{\'e}els  $ t_k \in
]\lambda_k,1/2],{\theta}_k \in {\mathbb S}_{n-1} $, tels que :

$$ { w_i}(2\lambda_k-t_k,{\theta}_k)-{
  w_i}(t_k,\theta_k)\geq 0.  \qquad (* )$$

On va voir que pour $ \lambda_k $  pris dans la famille pour
laquelle $ (*) $ est v{\'e}rifi{\'e}e, $ t_k\in
[\log(\beta_i/A_i),1/2] $. 

\smallskip

Supposons au  contraire que $ t_k<\log \dfrac{\beta_i}{A_i} $.

Lorsque  $\lambda_k$ est voisin de $-\infty $, nous avons : $ \lambda_k <  \log \dfrac{\beta_i}{A_i}$.

\smallskip

D'autre part, sachant qu'on a toujours  $t^{\lambda_k}<t$, en
prenant $ t=t_k $ dans $ ]\lambda_k  , \log \dfrac{\beta_i}{A_i} [$ et en utilisant
la croissance de $ w_i $ on obtient
l'in{\'e}galit{\'e} suivante:

        $$  {     w_i}(2\lambda_k-t_k,{\theta}          )-{
  w_i}(t_k,{\theta         })<0 \,\, {\rm pour\,\,tout}
\,\,    \theta \in {\mathbb S}_{n-1}. $$

En particulier pour $ \theta=\theta_k $ l'in{\'e}galit{\'e} obtenue,
  contredit $ (*) $.

\smallskip

Ainsi, pour tout $ \lambda_k \leq 0 $, pris dans la famille pour laquelle $ (*)
$ est v{\'e}rifi{\'e}e :

  $$ 1/2\geq t_k\geq   \log \dfrac{
    \beta_i}{A_i}. $$
(En particulier  $ \log \dfrac{
    \beta_i}{A_i} \leq \log\eta_i+4 $, car $ w_i(\log
  \eta_i,\theta)-w_i(\log \eta_i+4,\theta)>0 $).
\smallskip

 Par compacit{\'e} on obtient:

\smallskip

\qquad $ \lambda_k \to -\infty $, $
t_k\longrightarrow t_0 \in[\log \dfrac{\beta_i}{A_i},1/2]$ et $
{\theta}_k\longrightarrow \theta_0 \in {\mathbb S}_{n-1} $

\smallskip

Or,

\smallskip

\qquad  $ 0\leq {
  w_i}(2\lambda_k-t_k,{\theta}_k)-w_i(t_k,\theta_k) $, en faisant tendre  $ \lambda_k $ vers $  -\infty $,
  on obtient:

\smallskip

\qquad  $      u_i(a_i+e^{
t_0}\theta_0)\leq 0 $, or ceci est impossible car  $ u_i>0 $.

\smallskip

Ainsi, on a d{\'e}montr{\'e} que pour $ \lambda $ petit, voisin de $ -\infty $,
$ w_i ^{\lambda}(t,\theta)-w_i(t,\theta)<0,$ pour $ (t,\theta)\in
]\lambda,1/2]\times {\mathbb S}_{n-1} $.

\bigskip

\underbar {d{\'e}monstration du Lemme 2:}

\smallskip

On commence par prouver

$$ \partial_t V_i \geq  {\rm (termes\, positifs) } \times
{\left [ \dfrac{(n-2)a\epsilon_i}{2}-Ae ^t\right]}. \qquad (*) $$

En effet, comme $ V_i=V_i(t,\theta)=e^{[(n-2)\epsilon_i
  t]/2}V_{\epsilon_i}(a_i+e^t\theta) $, on a:

\smallskip

$ \partial_t V_i = \dfrac{(n-2)\epsilon_i}{2}      e^{[(n-2)\epsilon_i t]/2}V_{\epsilon_i}(a_i+e^t\theta)+e^{[(n-2)\epsilon_i t]/2}\times e ^t
  <\nabla V_{\epsilon_i} (a_i+e ^t\theta)|\theta> $

\smallskip

D'o{\`u},

$$\partial_t V_i\geq e^{[(n-2)\epsilon_i t]/2} \times [
\dfrac{(n-2)a\epsilon_i}{2}-Ae ^t], $$

o{\`u} $ A $ est un majorant de la norme infinie du gradient de $ V_i $.

Ainsi,

$$  t\leq \log \epsilon_i
+\log \dfrac{(n-2)a}{2A} \,\, \Rightarrow \,\, \partial_t V_i \geq 0. $$

Or,d'apr{\`e}s notre hypoth{\`e}se de d{\'e}part:

$$ {\epsilon_i}^{[2/(n-2)-\epsilon_i/2]^{-1}}      {[
  u_{\epsilon_i}(a_{\epsilon_i})]}^{\beta }\geq
\dfrac{c_{\epsilon_i}}{{(l_{\epsilon_i}\times
    R_{\epsilon_i})}^{(n-2)/2}}\geq 1 , $$

$$ \log  \epsilon_i \geq -(\dfrac{2}{n-2}-\dfrac{\epsilon_i}{2}) \beta       \log
[u_{\epsilon_i}(a_{\epsilon_i})]=\beta       \log \eta_{\epsilon_i} , $$

avec $
\eta_{\epsilon_i}={[u_{\epsilon_i}(a_{\epsilon_i})]}^{\epsilon_i/2-2/(n-2)}
$ ( $ a_{\epsilon_i}=a_i $ est le point d{\'e}fini dans l'{\'e}tape 1).

\smallskip

On voit alors que :

$$ \log \epsilon_i+\log \dfrac{(n-2)a}{2A}\geq
\beta\log\eta_{\epsilon_i
}+\log \dfrac{(n-2)a}{2A}=\bar {t_{\epsilon_i}}=t_i . $$

\smallskip

Ceci nous permet d'avoir la croissance en $ t $ de la fonction $
e^{(n-2)\epsilon_i t/2 } V_i $
sur l'intervalle $ ]-\infty, t_i ] $.

\smallskip

On d{\'e}montre maintenant le Lemme 2:

\smallskip

Supposons que pour un $ \lambda\leq 0 $ on ait:

$$ {w_i}^{\lambda}(t,\theta)-w_i(t,\theta)<0,\,\, \forall \, (t,\theta)\in
]\lambda,t_i]\times {\mathbb S}_{n-1}\,, $$

En notant $ \tilde V_i(t,\theta)=e^{(n-2)\epsilon_i t/2}V_i(t,\theta), 
{\tilde V_i}^{\lambda}(t,\theta)=\tilde V_i(t^{\lambda},\theta)=\tilde
V_i(2\lambda-t,\theta)
$, on peut {\'e}crire:

$$
-L({w_i}^{\lambda}-w_i)=({\tilde V_i}^{\lambda}-\tilde V_i){({w_i}^{\lambda})}^{N-\epsilon_i-
 1}+\tilde V_i[{({w_i}^{\lambda})}^{N-\epsilon_i-1}-{w_i}^{N-\epsilon_i-1}] \,.$$

\smallskip

On a vu que sur l'intervalle $ [\lambda,t_i] $, la fonction $
 t\rightarrow \tilde V_i(t,\theta)=e^{(n-2)\epsilon_i t/2} V_i(t,\theta)
 $ est uniform{\'e}ment
croissante et comme $ t\in [\lambda,t_i] $, $
   t^{\lambda}-t=2\lambda-t-t=2(\lambda-t)\leq 0 $, on en d{\'e}duit que :

$$  {\tilde V_i}^{\lambda}\leq \tilde V_i \,\,\,  {\rm
  sur }\,\,\, [
 \lambda,t_i] \times {\mathbb
  S}_{n-1} .$$

D'autre part, il existe un rang $ i_1 $ {\`a} partir duquel
  $ N-\epsilon_i-1>1>0 $ puisque $ \epsilon_i\to 0 $. Ainsi la fonction   $ t\mapsto t^{
  N-\epsilon_i-1} $ est croissante et on a finalement :

$$ {w_i}^{\lambda}<w_i \Rightarrow {(
{w_i}^{\lambda})}^{N-\epsilon_i-1} < {w_i}^{N-\epsilon_i-1} \,. $$

Le Lemme 2 est ainsi d{\'e}montr{\'e}.

\bigskip

\underbar { V{\'e}rification du point 3) :}

\smallskip

D'apr{\`e}s le lemme        de l'{\'e}tape 2-2:

$$      w_i(\log\eta_i,\theta)-     w_i(\log\eta_i+4,\theta)>0 . $$

\bigskip

On pose $ l_i=\log\eta_i+4 $ et  $  \bar \lambda_i=2+\log\eta_i $,  alors:

$$ 2\bar \lambda_i-l_i=2(\log\eta_i+2)-\log\eta_i-4=\log\eta_i . $$ 

Comme  $ \bar
\lambda_i<l_i<t_i $
, on obtient:

$$ { w_i}^{\bar \lambda_i}(l_i,\theta)- w_i(l_i,\theta)>0 $$

et finalement $ \xi_i $ existe bien. 

\bigskip

{\bf Etape 3: Utilisation du principe du maximum pour la conclusion }

\smallskip

Montrons que les fonctions $ {     w_i}^{\xi_i}-     w_i $ v{\'e}rifient
les propri{\'e}t{\'e}s suivantes: 

\smallskip

1) sur $ ]\xi_i,t_i]\times {\mathbb S}_{n-1}  $ , $ {     w_i}^{\xi_i}-
   w_i \leq 0 $,

\smallskip

2) sur $]\xi_i,t_i]\times  {\mathbb S}_{n-1}  $, $ -     L({     w_i}^{\xi_i}-
   w_i) \leq 0 $,

\smallskip

Pour le point 1), on utilise la  d{\'e}finition de $\xi_i $ : il existe
une suite $ \{\mu_{i,k}\} $ telle que

\smallskip

a) $ \mu_{i,k}<\xi_i $ pour tout entier $ k $,

\smallskip

b) $     w_i^{\mu_{i,k}} - w_i <0 $ sur $ ]\mu_{i,k},t_i]\times {\mathbb
  S}_{n-1}  $, ${\rm pour \, tout }\,\, k  $,

\smallskip

donc:

\smallskip

$  \qquad     w_i(2\mu_{i,k}-t,\theta)-     w_i(t,\theta)<0 $, $ {\rm
  pour \, } t\in  \,
]\mu_{i,k}, t_i[ $ et tout $ \theta\in {\mathbb S}_{n-1} $.

\smallskip

La fonction $      w_i $ est continue et tout $ t\in ]\xi_i,t_i] $ est
dans des $]\mu_{i,k},t_i]$ par a), en passant {\`a} la limite en $ k $
on obtient 1).

\smallskip

Pour le point 2), la d{\'e}monstration est identique {\`a} celle du 1), les fonctions $  w_i $
   sont $ C^2 $, il suffit d'ecrire   $ 
{ w_i}^{\mu_{i,k}}-w_i=w_i(2\mu_{i,k}-.,.)-     w_i(.,.)
   $.
\bigskip

\underbar { Lemme  :}

les fonctions $ {     w_i}^{\xi_i} $ et $
   w_i $ v{\'e}rifient:

$$  \max_{{}_{\theta \in {\mathbb S}_{n-1} }} {     w_i}^{\xi_i}(t_i,\theta)\geq
   \min_{{}_{\theta \in {\mathbb S}_{n-1} }} 
   w_i(t_i,\theta) . $$

\underbar{ D{\'e}monstration du Lemme :}

\smallskip

Supposons, par l'absurde, que:

$$ \max_{{}_{\theta \in {\mathbb S}_{n-1}
 }} {     w_i}^{\xi_i}(t_i,\theta)<
   \min_{{}_{\theta \in {\mathbb S}_{n-1} }} 
   w_i(t_i,\theta). $$

Alors:

\bigskip

$ \qquad \qquad \forall \, \theta\in {\mathbb S}_{n-1} $, $
{     w_i}^{\xi_i}(t_i,\theta)< 
   w_i(t_i,\theta) .\qquad (3) $

\smallskip

Notons : 

\smallskip

$ \qquad h(t,\theta)=  {     w_i}^{\xi_i}(t,\theta)-
w_i(t,\theta) $ sur   $ [\xi_i,t_i]\times
{\mathbb S}_{n-1} $.

\smallskip

En utilisant les propri{\'e}t{\'e}s 1), 2) et (3), la fonction  $ h $ v{\'e}rifie:

\smallskip

$ \qquad h(t,\theta)\leq 0 $ sur $ [\xi_i,t_i] \times {\mathbb S}_{n-1}  $ et $
h(t_i,\theta)<0, \, \forall\, \theta\in {\mathbb S}_{n-1} , $

\smallskip

$ \qquad Lh\geq 0 $ sur $ [\xi_i,t_i]\times{\mathbb S}_{n-1} $.

\smallskip

Par le principe du maximum de Hopf, on obtient:

\smallskip

$ h $ atteint son maximum sur le bord ou bien elle est constante. Si $ h $ n'est pas constante, elle v{\'e}rifie {\`a} l'int{\'e}rieur du
 domaine, $ h< \max h $. L{\`a} o{\`u} $ h $ atteint son maximum elle v{\'e}rifie: $ \partial_{\nu}h > 0
$. $ \nu $ la normale exterieure.

\smallskip

En $ (\xi_i,\theta) $, $ h $ est nulle et en $ (t_i,\theta ) $ elle est
strictement n{\'e}gative, elle ne peut pas {\^e}tre constante. Donc:

\smallskip

$ h<0 $ sur $ ]\xi_i, t_i]\times {\mathbb S}_{n-1} $ et $
\partial_{\nu}h(\xi_i,\theta) >0 $. Comme $ \partial_{\nu}=-\partial_t
$, on obtient: 

\smallskip

$\partial_{\nu}({     w_i}^{\xi_i}-
w_i)(\xi_i,\theta)=-\partial_t[     w_i(2\xi_i-t,\theta)-
w_i(t,\theta)]=2\partial_tw_i(\xi_i,\theta)>0 $.

\smallskip

En fixant $ i $, la d{\'e}finition de $ \xi_i $ comme borne sup{\'e}rieure
d'un certain ensemble pr{\'e}cedemment d{\'e}fini, donne:

\smallskip

 Pour tout $  k>0 $, il existe $ \mu_k,\sigma_k,\theta_k $ v{\'e}rifiant: $ \xi_i+\dfrac{1}{k}>\mu_k>\xi_i $, et \\$ \mu_k <\sigma_k\leq
t_i $, $ \theta_k \in {\mathbb S}_{n-1} $  tels que

\smallskip

$ \qquad { w_i}^{\mu_k}(\sigma_k,\theta_k)-
w_i(\sigma_k,\theta_k)= w_i(2\mu_k-\sigma_k,\theta_k)-w_i(\sigma_k,\theta_k)\geq 0 $.

\smallskip

\underbar{ 1er Cas :} si $ \sigma_k\rightarrow \sigma_0>\xi_i $ (ou au moins une valeur
d'adh{\'e}rence) : 

\smallskip

En passant {\`a} la limite
($ {\mathbb S}_{n-1}  $ est compacte, quitte {\`a} passer aux sous-suites,
$ \theta_k\rightarrow \theta_0 ) $ et en utilisant la continuit{\'e} de $
\bar w_i $ on
obtient:

\smallskip

$ \qquad w_i(2\xi_i-\sigma_0,\theta_0)-w_i(\sigma_0,\theta_0)\geq 0
$,

\smallskip

$\qquad {w_i}^{\xi_i}(\sigma_0,\theta_0)-w_i(\sigma_0,\theta_0)\geq
0 $ ce qui contredit le resultat trouv{\'e} plus haut sur $ h $ .

\smallskip

\underbar{ 2{\`e}me Cas :} si $ \sigma_k\rightarrow \xi_i $ :

\smallskip

Comme $  \dfrac{ w_i(2\mu_k-\sigma_k,\theta_k)-
  w_i(\sigma_k,\theta_k)}{2(\mu_k-\sigma_k)} \leq  0 $, en passant {\`a}
  la limite, on obtient:

\smallskip
 
$ \qquad \lim_{ {k\rightarrow +\infty}} \dfrac{
  w_i(2\mu_k-\sigma_k,\theta_k)-
  w_i(\sigma_k,\theta_k)}{2(\mu_k-\sigma_k)}
  =\partial_tw_i(\xi_i,\theta_0)\leq 0 $,

\smallskip

ce qui contredit l'in{\'e}galit{\'e} {\'e}tablie plus haut.

D'o{\`u} le Lemme  est d{\'e}montr{\'e}:

\smallskip

$ \alpha ) $  $ \min_{  {\mathbb S}_{n-1} } w_i(t_i,\theta)\leq \max_{ {\mathbb S}_{n-1} }w_i(2\xi_i-t_i,\theta) $.

\smallskip

De plus, comme $ t_i \to -\infty $, on obtient:

\smallskip

$ \beta )$  $ w_i(t_i,\theta)=e^{(n-2)t_i/2}u_i(a_i+e ^t\theta)\geq
e^{(n-2)t_i/2} \, \min_{B_i} u_i\geq
e^{(n-2)t_i/2}\, \min_{B_{1/2}(0)} u_i $, o{\`u} $ B_i $ est la boule de
centre $ a_i \to 0 $ et de rayon $ e^{t_i} < \dfrac{1}{2} $ .

\smallskip

Sachant que:

\smallskip

\qquad $ w_i(2{\xi_i}
-t_i,\theta)=e^{(n-2)(2\xi_i-t_i)/2}u_i(a_i+e^{2\xi_i-t_i} \theta) $,

\smallskip

et que:

\smallskip

 $ 2\xi_i-t_i=(2\xi_i-t_i-\bar \lambda_i)+\bar \lambda_i $ et $
 \xi_i\leq \bar \lambda_i \leq t_i $ $\Rightarrow $ $ s_i=2\xi_i-t_i-\bar
 \lambda_i \leq 0 $, 

\smallskip

nous pouvons {\'e}crire:

\smallskip

 $ w_i(2\xi_i-t_i,\theta)=w_i(2\xi_i-t_i-\bar \lambda_i+\bar
 \lambda_i,\theta)=w_i(s_i+2+\log \eta_i              ,\theta) $ avec $
 s_i\leq 0 $.

\bigskip

En utilisant une des propri{\'e}t{\'e}s des fonctons $ w_i $, vues dans
l'{\'e}tape 2,

\smallskip

$ w_i(2\xi_i-t_i,\theta)\leq c\, e ^{(n-2)(2\xi_i-t_i-\bar
  \lambda_i+2)/2} \,\, {u_i(a_i)}^{(n-2)\epsilon_i/4} $, o{\`u}  $ c $
  une constante positive ne d{\'e}pendant pas de $ i $.

\smallskip

Comme $ \xi_i\leq \bar \lambda_i $, on a:

\smallskip

$ \gamma $) $ w_i(2\xi_i-t_i,\theta)\leq c \,\,
{u_i(a_i)}^{(n-2)\epsilon_i/4} \,\, e ^{(n-2)(\bar \lambda_i-t_i)/2} $.

\smallskip

Ce qui peut s'{\'e}crire, en combinant $ \alpha), \beta), \gamma) $:

$$ e ^{(n-2)t_i/2} \times \min_{B_{1/2}(0)} u_i \leq c \,\,
{u_i(a_i)}^{(n-2)\epsilon_i/4}\,\, e ^{(n-2)(\bar \lambda_i-t_i)/2}. 
$$

Ou encore,

$$ e ^{(n-2)(-\bar \lambda_i+2t_i)/2}
 \times \min_{B_{1/2}(0)} u_i \leq c\,\,
 {u_i(a_i)}^{(n-2)\epsilon_i/4} , $$

ainsi

$$ {u_i(a_i)}^{     (1-2\beta)      (1-(n-2)\epsilon_i/2)}\,
\min_{B_{1/2}(0)} u_i \leq  c . $$

On voit qu'on s'est ramen{\'e} {\`a} une in{\'e}galit{\'e} du type  $
{[u_i(a_i)]}^{\delta} \times \min \, u_i \leq c $, avec $ \delta >0 $ (
car $ \beta <\dfrac{1}{2} \, \rm {et} \, \epsilon_i \to 0 $ ). 

\smallskip

Pour avoir la contradiction avec l'hypoth{\`e}se de d{\'e}part, il suffit que :

$$ (1-2\beta)(1-(n-2)\epsilon_i/2)\geq \beta \,\,\, {\rm pour \,\, tout }\,\, i .$$

On prend   $ \beta $ dans $ ]0,\dfrac{1}{3} [ $, on obtient une contradiction.

\smallskip

{\bf Etape 4:  preuve du Th{\'e}or{\`e}me 1 }

\smallskip

Soit $ x_0\in \Omega  $ alors  il existe un r{\'e}el $ r=r(\Omega)>0 $ tel que
, $ B_{r}(x_0)\in \Omega $.

\smallskip

Consid{\'e}rons la suite de fonctions :

$$ \bar u_i(x)= u_i(x_0+rx) \times r^{2/(N-\epsilon_i-2)},\,\,\, x\in
B_1(0) , $$

alors:

$$ \Delta \bar u_i= r^2 \Delta
  u_i(x_0+rx)r^{2/(N-\epsilon_i-2)} = V_i {u_i}^{N-\epsilon_i-1}
  r^{2(N-\epsilon_i-1)/(N-\epsilon_i-2)}=V_i {\bar u_i}^{N-\epsilon_i-1}, $$

d'apr{\`e}s le r{\'e}sultat qui pr{\'e}c{\`e}de ({\'e}tapes pr{\'e}c{\'e}dentes):

$$ \exists \, c,R>0,\,\, {\epsilon_i}^{(n-2)/2} {(\sup_{{}_{B_R(0)}}
  \bar u_i)}^{\beta } \times \inf_{{}_{B_1(0)}} \bar u_i \leq \dfrac{c}{R
  ^{4/(N-\epsilon_i-2)}} $$

et finalement:

$$ \forall\, x_0\in B_1(0),\,\, \exists c_{x_0},R_{x_0} > 0,\,
{\epsilon_i}^{(n-2)/2} {(\sup_{B_{R_{x_0}}(x_0)}  u_i)}^{\beta }
\times \inf_{\Omega }  u_i \leq c_{x_0} .$$

Soit $ K $ un compact de $ B_1(0) $, pour chaque $ x\in K $, on
consid{\`e}re le $ R_x $ comme pr{\'e}c{\'e}demment. Alors; $ K\subset \bigcup_{x\in K} B_{R_x}(x) $. Comme $ K $ est
compact, il existe  $ m\in {\mathbb N} $ tel que, $  K\subset \bigcup_{j=1}^m  B_{R_{x_j}}(x_j) $.

Donc:

$$ {\epsilon_i}^{(n-2)/2} {(\sup_{{}_{K}}  u_i)}^{\beta }
\times \inf_{{}_{\Omega }}  u_i \leq \sum_{j=1}^m 
{\epsilon_i}^{(n-2)/2} {(\sup_{B_{R_{x_j}}(x_j)}  u_i)}^{\beta }
\times \inf_{\Omega }  u_i \leq c(\beta,a,b,A,K,\Omega). $$

\bigskip

{\bf Th\'eor\`eme 2}. Consid\'erons deux suites de fonctions $ (u_{\epsilon_i}) $ et $ (V_{\epsilon_i}) $ relatives au probl\`eme $ (E_1) $, alors si
on suppose que:

\bigskip

$ q_i=N-\epsilon_i $ avec $ \epsilon_i \to 0 $ et $
||\nabla V_{\epsilon_i}||\leq k\epsilon_i $ ($ k>0 $),  alors il vient:

\bigskip

Pour tout compact $ K $ de $ \Omega $, il existe une constante $ c>0 $
ne d{\'e}pendant que de $ a,b,k,K,\Omega $, telle que :

$$ (\sup_{{}_K} u_{\epsilon_i})^{4/5 } \times \inf_{{}_{\Omega}} u_{\epsilon_i}
\leq c . $$

\underbar {\bf Preuve.}

\bigskip

La preuve utilise les m{\^e}mes t{\'e}chniques que celles mises en oeuvre dans
la d{\'e}monstration du Th{\'e}or{\`e}me 1. On suppose toujours que $ \Omega =B_2(0)\subset {\mathbb
  R}^n $ et on commence par {\'e}tablir des estimations locales telles que:

$$ \exists\, c=c(a,b,A)>0,\, \exists \, R>0, \, (\sup_{{}_{B_R(0)}}
u_{\epsilon_i })^{\beta }\times \inf_{{}_{B_1(0)}} u_{\epsilon_i} \leq
\dfrac{c}{R^{4/(N-\epsilon-2)}} .$$

Pour cela, on raisonne par l'absurde, les {\'e}tapes sont les m{\^e}mes que
celles de la d{\'e}monstration du Th{\'e}or{\`e}me 16, la diff{\'e}rence est 
que $ \epsilon^{(n-2)/2} $ absent du membre de droite et on
verra qu'on peut choisir l'exposant du sup aussi proche de $ 1 $ qu'on le veut.

\smallskip

On exhibe une suite de points $ (a_{\epsilon_i}) $, tendant vers $ 0 $,  telle
que:

$$ [u_{\epsilon_i}(a_{\epsilon_i})]^{\beta }\, \inf_{{}_{B_1(0)}}
u_{\epsilon_i}\to +\infty. \qquad (**) $$

Nous souhaiterions utiliser le principe du maximum. Pour cela, on regarde l'accroissement des fonctions $ V_i
$. Comme on a pos{\'e}, $ V_i(t,\theta)= e^{(n-2)\epsilon_i
  t/2}V_{\epsilon_i}(a_{\epsilon_i}+e^t\theta) $, on obtient:

$$ \partial_t V_i(t,\theta)\geq e^{(n-2)\epsilon_i
  t/2} \, [\dfrac{(n-2)a\epsilon_i}{2}-A_i e^t], $$

o{\`u} $ a $ est un minorant de $ V_i $ et $ A_i$ est un majorant de la
norme infinie du gradient des $ V_i $ (
La condition  $  A_i  \leq k\epsilon_i
$ ($ k>0 $), va {\^e}tre utilis{\'e}e ).  

donc:

$$ \partial_t V_i(t,\theta) \geq  e^{(n-2)\epsilon_i
  t/2}\, [\dfrac{(n-2)a\epsilon_i}{2}-k\epsilon_i e^t]\geq  e^{(n-2)\epsilon_i
  t/2}\, \epsilon_i[\dfrac{(n-2)a}{2}-ke^t] \,.$$

Ainsi on obtient la condition de croissance pour $ V_i $: 

\smallskip

$ {\rm Pour }\,\,  t\leq \log
\dfrac{(n-2)a}{2k}=t_0 \,\, \Rightarrow \partial_t V_i(t,\theta) \geq
0 \,\, {\rm pour \, tout } \, \theta \in {\mathbb S}_{n-1}, $ le $ t_0 $ ne d{\'e}pend pas de $ i $.

\smallskip

Soient $ w_i $ la fonction  $  w_i(t,\theta)=e^{(n-2)t/2}u_i(a_i+e^t\theta) $, $
\xi_i\leq \bar \lambda_i=2+\log \eta_i $ et 
$ \eta_i=\dfrac{1}{[u_{\epsilon_i}(a_i)]^{(N-\epsilon_i-2)/2}} $.
 
\smallskip

Comme dans la d{\'e}monstration du Th{\'e}or{\`e}me 1, en
supposant que $ \min_{\theta \in {\mathbb S}_{n-1}} w_i(2\xi_i-t_0,\theta) >
\max_{\theta \in {\mathbb S}_{n-1}} w_i(t_0,\theta) $ et en utilisant le
principe du maximum de Hopf, on aboutit {\`a} une contradiction.

\smallskip

Finalement on obtient:

$$ \min_{\theta \in {\mathbb S}_{n-1}} w_i(2\xi_i-t_0,\theta) \leq
\max_{\theta \in {\mathbb S}_{n-1}} w_i(t_0,\theta) .$$

En reprenant la cons{\'e}quence du Lemme  de l{\'e}tape 3 du Th\'eor\`eme 1, on obtient:

$$ w_i(t_0,\theta)\geq e^{(n-2)t_0/2}\, \min_{{}_{B_1(0)}} u_{\epsilon_i}, $$

$$ w_i(2\xi_i-t_0,\theta)\leq c\,
[u_i(a_{\epsilon_i})]^{(n-2)\epsilon_i/4}\, e^{(n-2)(\log
  \eta_i-t_0)/2}. $$

Donc:

$$  \min_{{}_{B_1(0)}} u_{\epsilon_i} \leq c \times
[u_i(a_{\epsilon_i})]^{(n-2)\epsilon_i/4}\,
\dfrac{1}{[u_{\epsilon_i}(a_i)]^{[1-(n-2)\epsilon_i/4]}} ,
$$

C'est {\`a} dire:

$$ [u_{\epsilon_i}(a_i)]^{1-(n-2)\epsilon_i/2} \,
\min_{B_1(0)} u_{\epsilon_i} \leq c ; $$

ceci contredit $ (**)$ car $ \beta <1-\dfrac{(n-2)\epsilon_i}{2} $
pour $ i\geq i_0 $ et $ [u_{\epsilon_i}(a_{\epsilon_i})]\to +\infty $.

\bigskip

\underbar{\it Les $ n= 3, 4 $ de l'\'equation de la courbure scalaire prescrite}

\bigskip

{\bf Th\'eor\`em 3}. Consid\'erons deux  suites de fonctions $ (u_i), (V_i)  $  relatives \`a l'\'equation $ (E_1) $ :

\bigskip

 Si $ n=3 $ et $ q=5 $, alors pour tout compact $ K $ de $ \Omega $, il existe une constante $ c>0 $ ne d\'ependant que de $ a,b,A,K,\Omega $ telle que:

$$ (\sup_{{}_K} u_i)^{1/3} \times
\inf_{{}_{\Omega}} u_i \leq c .$$

Si $ n=4 $, $ q=3 $ et si la
constante de Lipschitz $ A_i $, relative \`a $ V_i $, tend vers $ A\geq
0 $,  alors:

\bigskip

En supposant que $ \liminf_{i \to +\infty} \dfrac{\min_{\Omega}
  u_i}{A_i} \geq\dfrac{8e^2 \sqrt {2}}{3\,a \sqrt{a}}  $, on obtient:

\bigskip

Pour tout compact $ K $ de $ \Omega $, il existe une constante $ c>0 $, ne
d\'ependant que de $ a, b, (A_i)_{i\in {\mathbb N}}, K, \Omega $, telle que:

$$(\sup_{{}_K} u_i) \times
 \inf_{{}_{\Omega}} u_i \leq c .$$

\underbar {\bf Preuve.}

\bigskip

\underbar{ \bf 1er Cas: $ n=3, q=N=6 $  :}

\smallskip

La d{\'e}monstration est similaire {\`a} celle  du Th{\'e}or{\`e}me 1. Elle utilise  les
techniques " blow-up " et " moving-plane ".

\smallskip

\underbar { Etape 1: la technique blow-up }

\smallskip

Commen{\c c}ons par {\'e}tablir la propri{\'e}t{\'e} suivante:

$$\exists\,R\in]0,1[,\,\exists \,c>0, {(\sup_{B_R(0)} u_i)}^{1/3 } \times \inf_{B_1(0) } u_i \leq \dfrac{c}{R}.
$$

En supposant le contraire, on exhibe une sous-suite $\{u_j\} \subset
\{u_i\} $, une suite de points de la boule
unit{\'e} $ (a_j) $ et trois suites de r{\'e}els positifs $ (R_j), \,
(c_j), \, (l_j) $ telles que: 

$$ a_j\to 0, \,\,\,\, c_j\to +\infty, \,\,\,\,  R_j\to 0 \,\,\,\, et
\,\,\,\, l_j\to 0 $$

$$ \Delta u_{j}=V_{j}{u_{j}}^5 $$

$$ [u_{j}(a_j)]^{1/3}\times \inf_{ B_1(0) } u_{j} \geq
\dfrac{c_j}{R_j} $$

Comme on raisonne par l'absurde, on peut supposer que $ u_i=u_j $.

\smallskip

D'autre part, on a vu qu'on peut construire {\`a} partir de $ (u_i) $, une suite
$ (v_i) $ v{\'e}rifiant:

$$ v_i(y)=\dfrac{u_i {\left [
      a_i+\dfrac{y}{{u_i(a_i)}^{2}}\right ] }}{u_i(a_i)} \,\,\,\, si
\,\,   |y|\leq \dfrac{l_i}{ \root 4 \of {c_i} }{[u_i(a_i)]^2}, $$

$$ \Delta v_i=V_i{v_i}^5, $$

$ \qquad v_i\rightarrow v=\dfrac{1}{{(1+|y|^2)}^{1/2}} $, uniform{\'e}ment sur $
  B_{\beta }(0) $, $\forall \beta >0 $. 

\smallskip

\underbar { Etape 2: Passage en polaires et propri{\'e}t{\'e}s de certaines fonctions }

\smallskip

soit $L_0$ et $L$ les op{\'e}rateurs:

$$L_0={\partial}_{tt}+{\partial}_{t}-{\Delta}_{\sigma} \qquad
et\qquad   L={\partial}_{tt}-{\Delta}_{\sigma}, $$

avec ${\Delta}_{\sigma}$ l'op{\'e}rateur de Laplace-Beltrami sur ${\mathbb
  S}_2 $. 

\smallskip

En posant: $ h_i(t,\theta,\phi)=u_i(a_i+e^t \cos\theta
\sin\phi...,...,...) $, on obtient :

$$-L_0h_i= e^{2t} V_i {h_i}^{5}. $$

Notons, $ w_i=e^{t/2} h_i$, on a alors:

\smallskip

 $ \qquad \qquad -Lw_i=-\dfrac{1}{4}w_i
+V_i{w_i}^{5} $ avec   $ w_i > 0. $
 
\smallskip

Consid{\'e}rons l'op{\'e}rateur $ \bar L=L-\dfrac{1}{4} $, $ -\bar
Lw_i=V_i(a_i+e^t\theta){w_i}^5 $.

\smallskip
 
\underbar {Etablissons quelques propri{\'e}t{\'e}s des fonctions  $ w_i $ : }

\smallskip

En posant $ \eta_i=\dfrac{1}{
{u_i(a_i)}^2}$ on obtient :

a) La suite  $ w_i(t+\log \eta_i             ,\theta
   )= e^{t/2} \,
   \dfrac{u_i(a_i+\dfrac{e^t\theta}{{u_i(a_i)}^2})}{u_i(a_i)}= e^{t/2}
   \, v_i(e^t\theta)$ converge vers la fonction sym{\'e}trique $  w= {\left
  (\dfrac{e^{t}}{1+
  e^{2t}} \right )}^{1/2} $ uniform{\'e}ment sur $ ]-\infty
   ,\log\beta]\times {\mathbb S}_2 $, pour tout  $ \beta >0 $.

\smallskip

b)   Pour $ i\geq i_0,  w_i(\log {\eta}_i,\theta)- w_i(\log
{\eta}_i+4,\theta)>0 $ pour tout $ \theta $. 

\smallskip

c) Si $ \lambda >0 $, $ \bar w_i= w_i-\lambda e^t $ v{\'e}rifie :

$$ \bar w_i(\log \eta_i,\theta)-\bar w_i(4+\log\eta_i,\theta)>0,\,\,
{\rm pour \, tout } \,\, \theta. $$

\smallskip

d) On a:

\smallskip

 $ \qquad \forall \delta\leq 0, \, \exists c(\delta)>0, \, i_0=i(\delta)\in
  {\mathbb N}  $ tels que

$ t\leq \delta $ $ \Rightarrow $  $ \dfrac{1}{c(\delta)}e ^{t/2} \leq
w_i(t+\log\eta_i,\theta) \leq c(\delta) e^{t/2} $ pour $ i\geq i_0 $
et tout  $
\theta \in {\mathbb S}_2 $.

\smallskip

Les  in{\'e}galit{\'e}s  b) et c)  permetteront de preciser le sup des r{\'e}els pour
lesquels la propri{\'e}t{\'e} relative {\`a} un ensemble not{\'e} $ A_{\lambda}
$ (qu'on d{\'e}finira plustard),
est non vide .

\smallskip

L'in{\'e}galit{\'e}  d)  est tr{\'e}s importante et sera utilis{\'e}e  vers la fin,
pour aboutir {\`a} une contradiction.
 
\smallskip

\underbar { D{\'e}monstration de b) }

\smallskip

$ w_i(\log \eta_i,\theta )-w_i(\log \eta_i+4,\theta)=[w_i(0+\log
\eta_i,\theta)-w(0)] - [w_i(4+\log \eta_i,\theta)-w(4)]+[w(0)-w(4)] $.

\smallskip

La convergence uniforme des $ w_i $ nous permet d'avoir, pour tout $
\epsilon >0 $, un rang  $ i_0 $ {\`a} partir duquel :

$$  w_i(\log \eta_i,\theta )-w_i(\log \eta_i+4,\theta)\geq
-2\epsilon+[w(0)-w(4)]. $$

De plus :

$$ [w(0)]^2-[w(4)]^2=\dfrac{1}{2}-\dfrac{e^4}{1+e^8}=\dfrac{1+e^8-2e
 ^4}{1+e^8}=\dfrac{(e^4-1)^2}{1+e^8}>0. $$ 

En prenant $ \epsilon < \dfrac{w(0)-w(4)}{2} $, on obtient b).

\smallskip

\underbar { D{\'e}monstration de c): }

\smallskip

 Soit $ C=\bar w_i(\log \eta_i,\theta)-\bar w_i(4+\log
\eta_i,\theta) $, alors d'apr{\`e}s b),

\smallskip

 $ C=w_i(\log\eta_i,\theta)-w_i(\log\eta_i+4,\theta)-\lambda(e^{\log
  \eta_i} -e^{\log \eta_i+4} )> \lambda \eta_i (e ^4-1)>0 $ .

\smallskip

\underbar { D{\'e}monstration de d)}

\smallskip

Par d{\'e}finition de $ w_i $ et d'apr{\`e}s les propri{\'e}t{\'e}s de $ u_i
  $, 

\smallskip

$  \dfrac{w_i(t+\log \eta_i,\theta)}{e^{t/2}}=\dfrac{u_i(a_i+\dfrac{e^t\theta}{[u_i(a_i)]^2})}{u_i(a_i)}\rightarrow
  \dfrac{1}{\sqrt{1+e^{2t}}} $, et la convergence est uniforme  en $ \theta $ d'apr{\`e}s
  a),

\smallskip

c'est-{\`a}-dire:

\smallskip

$ \forall \epsilon>0 $, $ \exists i_0(\epsilon,\delta)>0 $, $ i\geq
i_0 $, $ -\epsilon \leq \dfrac{w_i(t+\log \eta_i,\theta)}{e^{t/2}}
-\dfrac{1}{\sqrt{1+e^{2t}}} \leq \epsilon $, pour tout  $\theta $ et
$ t\leq \delta $. 

\smallskip

Comme $ \dfrac{1}{\sqrt{1+e^{2\delta }}} \leq \dfrac{1}{\sqrt{1+e^{2t}}} \leq
1 $, nous avons pour tout $ \theta $ et $ t\leq \delta $ :

\smallskip

$ \qquad -\epsilon+\dfrac{1}{\sqrt{1+e^{2\delta }}}\leq \dfrac{w_i(t+\log
  \eta_i,\theta)}{e^{t/2}}\leq \epsilon+1 $.

\smallskip

le choix suivant de $ \epsilon $ permet d'avoir d):

$ \qquad \epsilon=\dfrac{1}{2\sqrt{1+e^{2\delta }}} $ et $ c(\delta)=1+\epsilon $
\smallskip

\underbar { Etape 3: Utilisation de la t{\'e}chnique " moving-plane "}

\smallskip

Posons :

\smallskip

$ \qquad \bar w_i=w_i -{\lambda}_0e^t $ o{\`u} ${\lambda}_0\geq 0 $ \, est {\`a}  choisir  convenablement,

\smallskip

$ \qquad t^{\lambda}=2\lambda-t $ et $ {\bar
  w_i}^{\lambda }(t,\theta )={\bar w_i}(2\lambda-t,\theta ) 
 $,
\smallskip

$ \qquad z_{i,\lambda }= {w_i}^{\lambda } - w_i $ avec  $ \lambda \leq 0 $.

\smallskip

\underbar { Quelques lemmes importants {\`a} propos des fonctions $ {\bar
   w_i}^{\lambda}-\bar w_i
   $ }

\smallskip

\underbar { Lemme 1 :} 

\smallskip

pour $ 0<\beta<1 $ il existe un ${\lambda}_0\geq {\mu}_0 $ tel que:

$\bar w_i(t,\theta )>0 \,\,\, {\rm si } \,\,\,  t\leq t_i={\beta }\log
\eta_i $ pour tout $ \theta $ et tout $  i $.

\smallskip

\underbar { Trois remarques  :}

\smallskip

i) Il est clair que $ t_i= \beta \log \eta_i >4+\log \eta_i $ pour $ i\geq
i_0 $.

\smallskip

ii) Le choix de l'intervalle $ ]-\infty, t_i] $, nous permet de conserver
la positivit{\'e} de la fonction, car le choix d'un $ \lambda_0 $, ne
permet pas forc{\'e}mment de conserver la positivit{\'e} de la fonction, si on
prend  0 au lieu de $ t_i $.

\smallskip

iii) Le choix d'un $ \beta \in]0,1[ $, nous permet de conserver  une certaine
marge de manoeuvre pour la suite (obtenir notre r{\'e}sultat en utilisant
notre hypoth{\`e}se de d{\'e}part). On prendera, $ \beta=\dfrac{1}{3} $.

\underbar{Lemme 2:} 

\smallskip

Soit $ A_{\lambda} $, la propri{\'e}t{\'e}  suivante:

\smallskip

\qquad  $ A_{\lambda}=\{\lambda\leq 0,\,\,\exists \,\,
 (t_{\lambda},\theta_{\lambda })\in
   ]\lambda,t_i]\times  {\mathbb S}_2,\,\,   {\bar
   w_i}^{\lambda}(t_{\lambda},\theta_{\lambda})-{\bar
   w_i}(t_{\lambda},\theta_{\lambda} ) \geq 0\} . $

\smallskip

Alors, il existe $ \nu \leq 0$, tel que  pour $\lambda \leq \nu $, $ A_{\lambda}  $ n'est pas vraie.

\smallskip

\underbar { Remarque  :}

\smallskip

 Ce lemme est tr{\'e}s important car il pr{\'e}cise le domaine d'existence des
r{\'e}els $ \lambda $ tels que $  {\bar w_i}^{\lambda}-\bar w_i<0 $, o{\`u}
le
Lemme 3, ci-dessous, peut {\^e}tre utilis{\'e}.

\smallskip

\underbar{Lemme 3:} 

\smallskip
 
Soit, $ \lambda $ un r{\'e}el quelconque inf{\'e}rieur {\`a} $ \bar
 {\lambda_i}=2+\log\eta_i $. Alors:

\smallskip
 
$ \exists \,\,  {\mu}_0>0, $ tel que si $   {\lambda}_0\geq {\mu}_0 $
et donc

$$ {\bar w_i}^{\lambda} -{\bar w_i}<0 \Rightarrow -\bar L({\bar w_i}^{\lambda}
-{\bar w_i})<0 \, .$$ 

\smallskip

Ne pas confondre  $ \lambda $ qui donne le sym{\'e}trique de la fonction
et le ${\lambda}_0$ qui nous permet de construire la fonction $\bar
{w_i}=w_i-\lambda_0e^t $.

\smallskip

\underbar {Une int{\'e}rpr{\'e}tation de ce Lemme 3:}

 Ce lemme implique, qu'il existe une valeur $ \mu_0 $ d{\'e}pendant que de
$ A $, telle que si on se donne n'importe quelle suite $ {\delta_i} $ 
 avec pour tout $ i, \, \delta_i \leq \bar { \lambda_i } $, alors:

$$ {\bar w_i}^{\delta_i}-\bar w_i <0 \Rightarrow    -\bar L({\bar
  w_i}^{\delta_i}-\bar w_i)<0 . $$

4)\underbar{ Un point utile:}

\smallskip

 ${\xi}_i=$ sup $\{\lambda \leq {\bar \lambda_i}=2+\log
  \eta_i, {\bar w_i}^{\lambda} -\bar w_i < 0 $,sur $
  ]\lambda,t_i]\times {\mathbb S}_2 \} $,

\smallskip

$  {\xi}_i $ existe toujours d'apr{\`e}s le Lemme 2.

\smallskip

\underbar { D{\'e}monstration du Lemme 1 :}
 
\smallskip

Ecrivons:

\smallskip

\qquad $ \bar
w_i(t,\theta)=e^{t/2}u_i(a_i+e^{t}\theta)-{\lambda_0}e^{t}=e^t\{e^{-t/2}u_i(a_i+e^t\theta)-\lambda_0\}$. Alors:

\smallskip

\qquad $\bar w_i(t,\theta)>0 $  $ \Leftrightarrow $  
$ e^{-t/2}u_i(a_i+e^t\theta)>\lambda_0 $

\smallskip

Rappelons que $ t_i=\beta \log \eta_i $ avec $
\eta_i=\dfrac{1}{[u_i(a_i)]^2} $. Nous avons :  

$$\lambda_0\leq e^{(t_i-t)/2}
e^{-t_i/2}u_i(a_i+e^t\theta) . $$

Or pour  $t\leq t_i $ on a, $ e^{(-t_i/2)}\leq e^{-(t/2)} \Rightarrow
e^{(-t_i/2)}\min u_i \leq  e^{(-t/2)} u_i(a_i+e^t\theta) $,

\smallskip

donc:

\smallskip

\qquad  $ \lambda_0\leq {u_i(a_i)}^{\beta } \min u_i \leq
e^{-(t/2)}u_i(a_i+e^t\theta)$ pour $ t\leq t_i $.

\smallskip

D'apr{\`e}s notre hypoth{\`e}se de d{\'e}part ( celle qui doit aboutir {\`a} une
absurdit{\'e}) en prenant $
\beta=\dfrac{1}{3} $, on a:

$$ {u_i(a_i)}^{\beta } \min u_i \to +\infty . $$

Le r{\'e}el $\lambda_0 $ peut etre choisit convenablement, on le choisira
de telle sorte qu'on ait :

\smallskip

\qquad $ \lambda_0\leq
(1/2){u_i(a_i)}^{\beta } \min u_i$ $ \leq
(1/2)e^{(t/2)}u_i(a_i+e^t\theta)$  pour 
$ t\leq t_i $ 

\smallskip

et en cons{\'e}quence,

\smallskip

\qquad $ {u_i(a_i)}^{\beta } \min u_i-\lambda_0\geq
(1/2){u_i(a_i)}^{\beta } \min u_i $.

\smallskip

Par exemple, on peut prendre $\lambda_0=\lambda_{0,i}=(1/2){u_i(a_i)}^{\beta } \min u_i $
(on a une d{\'e}pendance en fonction de i).

\smallskip

Pour alleger l'{\'e}criture, on mettra  $ \lambda_0 $ devant $ e^t $
au lieu de $ \lambda_{0,i} $ dans
l'expression de $ \bar w_i $. 

\smallskip

\underbar {D{\'e}monstration du Lemme 2: }

\smallskip

D'abord, on fixe l'entier i et on cherche le signe de ${\partial}_t{\bar w_i}$

\smallskip

${\partial}_t{\bar
  w_i}(t,\theta)=(1/2)e^{t/2}u_i(a_i+e^t\theta)+e^{(3/2)t}
  [\theta^1{\partial}_1 u_i(a_i+e^t\theta)+\theta^2{\partial}_2
  u_i(a_i+e^t\theta)]-\lambda_0e^t $,

\smallskip

\qquad ${\partial}_t{\bar
  w_i}=e^t\{(1/2)e^{-t/2}u_i(a_i+e^t\theta)-\lambda_0+e^{t/2}(\theta^1
  {\partial}_1u_i+\theta^2\partial_2u_i) $,

o{\`u} $ (\theta ^1,\theta^2)=\theta $ un point de la sph{\'e}re $ {\mathbb
  S}_2 $.

\smallskip

Comme $ u_i $ est suppos{\'e}e $ C^1 $, il existe  $ A_i$ tel que, $
 \parallel{\nabla u_i \parallel}_{\propto}\leq A_i $.

\smallskip

D'autre part, d'apr{\`e}s le choix de $\lambda_0$ ( fin de la preuve du Lemme 1):

\smallskip

\qquad \qquad $ (1/2)e^{-t/2}u_i(a_i+e^t\theta)-\lambda_0\geq
\beta_i=\dfrac{1}{2} [u_i(a_i)]^{\beta}\min u_i > 0 $  pour $ t\leq
t_i $.

En cons{\'e}quence, pour $ t\leq t_i $ on obtient, $  {\partial}_t{\bar w_i}\geq e^t(\beta_i -e^{t/2}A_i) $.

\smallskip

Ainsi pour $ t<2\log \dfrac{\beta_i}{A_i} $, $( \beta_i -e^{t/2}A_i
\geq 0 )$, la fonction $
\bar w_i$ est  strictement
croissante uniform{\'e}ment en
$\theta \in {\mathbb S}_2 $.

\bigskip

Comme $ \bar w_i(\log \eta_i,\theta)-\bar w_i(\log \eta_i+4,\theta)>0
$, on obtient $ 2\log \dfrac{\beta_i}{A_i} \leq \log \eta_i+4 <t_i$.

\smallskip

Supposons que Lemme 2 ne soit pas vrai :

\smallskip

Il existe une famille de $ \{\lambda_k\} $, telle que $ \lambda_k
\to -\infty $, $  b_k \in
]\lambda_k,t_i] $ et $ {\theta}_k\in {\mathbb S}_2 $, telles que 

$$ {\bar w_i}(2\lambda_k-b_k,{\theta}_k)-{\bar
  w_i}(b_k,\theta_k)\geq 0 . \qquad (*)$$

Pour $ \lambda_k $ voisin de $-\infty $,  $ \lambda_k $
verifie: $ \lambda_k < 2\log \dfrac{\beta_i}{A_i} $ et donc pour $
 t\in ]\lambda_k, 2\log \dfrac{\beta_i}{A_i}] $, la fonction  $\bar w_i$ est
strictement croissante.

\smallskip

Comme $ t^{\lambda }=2\lambda-t \leq t $ pour $ \lambda \leq t $, on obtient alors:

 $$  {\bar w_i}(2\lambda-t,\theta)-{\bar
  w_i}(t,\theta)<0 \,\, {\rm pour \,\, tout } \,\,  (t,\theta) \in
  ]\lambda, 2\log \dfrac{\beta_i}{A_i}] \times {\mathbb S}_2. $$

Le r{\'e}el $ b_k $ v{\'e}rifie, avec $ \theta_k $, l'in{\'e}galit{\'e}
$ (*) $, il v{\'e}rifie {\'e}galement l'in{\'e}galit{\'e} suivante:
 
 $$ t_i\geq b_k\geq  2\log {(
    \beta_i/A_i)}  \,\, {\rm pour \, tout } \, \,  k .$$

Par compacit{\'e}, on obtient $ \lambda_k \to -\infty $ ,
 $ b_k\to t_0 \in[-2\log \dfrac{\beta_i}{A_i},t_i]$ et\\ $ {\theta}_k\to\theta_0,\,\,\theta_0 \in {\mathbb S}_2 $.

\smallskip

Comme les fonctions $ \bar w_i $ sont continues :

\smallskip

\qquad  $  {\bar
  w_i}(2\lambda_k-b_k,{\theta}_k)-\bar w_i(b_k,\theta_k) $ $\longrightarrow \lambda_0e^{
t_0}
  -e^{t_0/2}u_i(a_i+e^{t_0}\theta_0) = -\bar w_i(t_0,\theta_0) $, quand $ \lambda \to -\infty $.

\smallskip

En utilisant $ (*) $, on obtient : 

\smallskip

\qquad  $ \bar w_i(t_0,\theta_0)\leq 0 $ et $t_0 \leq t_i $ ce qui
contredit le Lemme 1 (d'o{\`u} le choix de $ t_i $ et pas de 0, pour la
borne de droite ).

\smallskip

\underbar { D{\'e}monstration du Lemme 3 :}

\smallskip

Consid{\'e}rons l'op{\'e}rateur  $ \bar L= L-\dfrac{1}{4}=\partial_{tt}-\Delta_{\sigma
  }-\dfrac{1}{4} $, $ \Delta_{\sigma} $ le laplacien sur $ {\mathbb
  S}_2 $.

\smallskip

On a $ -\bar L \bar w_i= -\bar L(w_i-\lambda_0 e^t)=-\bar
Lw_i+\lambda_0\bar L\, e^t=V_i{w_i}^5+\dfrac{3}{4}\lambda_0 e^t $ avec
$ V_i(t,\theta)=V_i(a_i+e^t\theta) $.

\smallskip

De m{\^e}me, $ \bar L {\bar w_i}^{\lambda }
=V_i^{\lambda}{(w_i^{\lambda})}^5+\dfrac{3}{4}\lambda_0 e^{2\lambda-t}
$, o{\`u} on a pos{\'e} $ V_i^{\lambda}(t,\theta)=V_i(a_i+e^{2\lambda-t}\theta) $.

\smallskip

Ainsi,

\smallskip

$ -\bar L(w_i^{\lambda}-w_i)=\dfrac{3\lambda_0}{4}(e^{ 2 \lambda - t
    }-e^t)+(V_i^{\lambda}-V_i){(w_i^{\lambda})}^5+ V_i
    [{(w_i^{\lambda})}^5-w_i^5] $.

\smallskip

Or, $ {V_i}(a_i+e^{t^{\lambda }}\theta)
-V_i(a_i+ e^t\theta) \leq ||\nabla V_i||_{\infty} 
(e^t-e^{t^{\lambda}})\leq A(e^t-e^{t^{\lambda}}) $ si $ \lambda<t $
 (ce qui est toujours le cas ici), d'o{\`u}

\smallskip

$-\bar L({\bar w_i}^{\lambda} -{\bar w_i})\leq [\dfrac{3\lambda_0}{4}-A
{({w_i}^{\lambda })}^5](e^{t^{\lambda}}
-e^t)+V_i\{{({\bar w_i}^{\lambda}+{\lambda_0}e^{t^{\lambda}})}^5 -(\bar
w_i +{\lambda_0e^t)}^5\} . $

\smallskip

Alors pour avoir,

\qquad  ${\bar w_i}^{\lambda} -{\bar w_i}<0 \Rightarrow -{\bar L}({\bar
  w_i}^{\lambda} -{\bar w_i})<0 \qquad (*) $,

\smallskip

il suffit que

\smallskip

\qquad  $ \dfrac{3\lambda_0}{4}-A {({w_i}^{\lambda })}^5 \geq 0 $.

\smallskip

Comme  $ \lambda \leq \bar \lambda_i=2+\log \eta_i $, $ 2\lambda-t-\bar \lambda_i=(\lambda-\bar \lambda_i)+(\lambda-t)\leq 0
$, alors:

$$ w_i(2\lambda-t-\bar \lambda_i+\bar \lambda_i,\theta)\leq (1+\epsilon
 )e^{(2\lambda-t-\lambda_i ) }\leq 1+\epsilon $$

 o{\`u} $ \epsilon $
 est un r{\'e}el positif fix{\'e}.

\smallskip

Pour {\'e}tablir cette in{\'e}galit{\'e},  on utilise la propri{\'e}t{\'e} d) de l'{\'e}tape
2:

\smallskip

\qquad  $ \forall \beta>0 ,\,\, w_i(t+\bar
 \lambda_i,\theta)=w_i(t+\delta+\log \eta_i,\theta) $  converge
 uniform{\'e}ment vers  $ w $ sur $ ]-\infty,\log \beta]\times {\mathbb S}_2
 $. On prendera $ t=2\lambda -t-\bar \lambda_i\leq 0 $ et $ \beta=1 $.

\smallskip

Finalement pour avoir $ (*) $, il suffit de choisir $ \lambda_0 \geq (3A/4)(1+\epsilon)^5 =\mu_0$ et on note que $\mu_0 $
 ne d{\'e}pend pas de $ \lambda\leq \bar{ \lambda_i} $.

\smallskip

\underbar { D{\'e}monstration du point utile 4:}

\smallskip

D'apr{\`e}s la propri{\'e}t{\'e} d) de l'{\'e}tape 2-1:

$$ \bar w_i(\log\eta_i,\theta)-\bar w_i(\log\eta_i+4,\theta)>0 .$$

Posons, $ l_i=\log\eta_i+4 $, on a alors ;

\smallskip

$ \qquad 2\bar \lambda_i-l_i=2(\log\eta_i+2)-\log\eta_i-4=\log\eta_i $ et $ \bar
\lambda_i<l_i<t_i .$

\smallskip

Donc: 

\smallskip

$$ {\bar w_i}^{\bar \lambda_i}(l_i,\theta)-\bar w_i(l_i,\theta)>0, $$

\smallskip

 $ \xi_i $ existe .

\smallskip

{\bf Etape 4: Utilisation des lemmes pr{\'e}c{\'e}dents et  conclusion : } 

\smallskip

On choisit les  $ {\lambda_{0,i}} $ comme dans le lemme 1, puis on
d{\'e}termine les $ \xi_i  $ correspondant aux $\lambda_{0,i} $ du Lemme
2, et apr{\`e}s on peut utiliser le Lemme 3.

\smallskip

Les fonctions $ {\bar w_i}^{\xi_i}-\bar w_i $ v{\'e}rifient les propri{\'e}t{\'e}s
suivantes: 

\smallskip

1) sur $ ]\xi_i,t_i]\times {\mathbb S}_2 $ , $ {\bar w_i}^{\xi_i}-\bar
   w_i \leq 0 $,

\smallskip

2) sur $]\xi_i,t_i]\times  {\mathbb S}_2 $, $ -\bar L({\bar w_i}^{\xi_i}-\bar
   w_i) \leq 0 $.

\smallskip

D'o{\`u} par le principe du maximum, on a le:

\smallskip

\underbar{ Lemme :}

Les fonctions $ {\bar w_i}^{\xi_i} $ et $\bar
   w_i $ v{\'e}rifient:

$$  \max_{{}_{\theta \in {\mathbb S}_2}} {\bar w_i}^{\xi_i}(t_i,\theta)\geq
   \min_{{}_{\theta \in {\mathbb S}_2}} \bar
   w_i(t_i,\theta). $$

La preuve du Lemme est identique {\`a} celle du Lemme 3 du Th{\'e}or{\`e}me 1.

\smallskip

D'apr{\`e}s le choix de $ \lambda_0 $ dans la fin de la preuve du Lemme 2,
on a :

$$ {\bar w_i}(t_i,\theta)\geq (1/2)e^{t_i/2} \min
u_i . $$

D'autre part, d'apr{\`e}s le point d) de l'{\'e}tape 2):

$ w_i(t+\bar \lambda_i,\theta)=w_i(t+\delta+\log
 \eta_i,\theta )\rightarrow w(t+\delta)\leq e^{(t+\delta)/2} $,
 uniform{\'e}ment sur $ ]-\infty,\log \beta]\times {\mathbb S}_2(1). $

\smallskip

D'o{\`u} :

$$ w_i(2\xi_i-t_i,\theta)\leq (1+\epsilon )e^{\delta /2} e
^{(2\xi_i-t_i-\bar \lambda_i)/2}\leq c e^{(\bar \lambda_i-t_i)/2} .$$

Ce qui peut s'{\'e}crire:

 $$ e^{(1/2)(2t_i-\bar \lambda_i )} \min u_i \leq c, \,\,{\rm pour \,
   tout } \,  i .$$

Comme $ \bar \lambda_i=4+\log \eta_i $, $ t_i =\beta \log \eta_i =\dfrac{1}{3} \log
\eta_i $ et $ \eta_i=[u_i(a_i)]^{-2} $, on en d{\'e}duit que :

 $$ {u_i(a_i)}^{1/3 } \times \inf u_i \leq c. $$

Ceci contredit notre hypoth{\`e}se de d{\'e}part ({\'e}tape 1).

\bigskip

\smallskip
{\bf \underbar { Cas :  $ n=4,\, q=N=4 $  }}

\smallskip

Dans ce cas, on suppose les fonctions $ V_i $ lipschitziennes de
constantes $ A_i\to A\geq 0 $. 

\smallskip

La preuve est assez similaire {\`a} celle de la dimension 3. On se place
sur $ \Omega =B_2(0) $. 

\smallskip

Supposons que 

$$\liminf_{i \to +\infty } \dfrac{\min_{B_2(0)} u_i}{A_i}\geq
\dfrac{8e^2\sqrt 2}{3a\sqrt a} $$

et montrons alors :

$$ \forall \, R>0, \,\,\,\, \sup_{B_R(0)} u_i \times \inf_{B_2(0)} u_i \leq
c=c(a,b,(A_i)_{i\in \mathbb N},R) .$$

\underbar { Etape 1: technique blow-up }

\smallskip

On d{\'e}montre d'abord la propri{\'e}t{\'e} suivante:

$$ {\rm il \, existe \, }\, c>0 \, {\rm et }\, R\in ]0,1[ \, {\rm tels \,\,
  que } \,\, (\sup_{{}_{B_R(0)}} u_i) \times
\inf_{{}_{B_2(0)}} u_i \leq \dfrac{c}{R^2} .$$ 

Supposons le contraire, alors:

$$ {\rm pour \, tout } \, c, R>0, \, {\rm il \, existe }  \, i_j\in
{\mathbb N}, \, {\rm tels \, que } 
( \sup_{{}_{B_R(0)}} u_{i_j}) \times \inf_{{}_{B_2(0)}} u_{i_j} \geq
\dfrac{c}{R^2} .$$

Le but est d'arriver {\`a} une contradiction; on peut donc  supposer que la
suite extraite est la suite elle-m{\^e}me:

\smallskip

Etant donn{\'e}es deux suites $ c_i\to +\infty $ et $ R_i \to 0 $, il
existe une suite $ u_i $ telle que :

$$  \sup_{{}_{B_{R_i}(0)}} u_{i} \times \inf_{{}_{B_2(0)}} u_{i} \geq
\dfrac{c_i}{{R_i}^2} .$$

\smallskip

Comme, la suite des minima est born{\'e}e, on en d{\'e}duit:

$$  \sup_{{}_{B_{R_i}(0)}} u_{i}\times R_i^2\geq c_i .$$

Introduisons les fonctions suivantes :

$$ s_i(x)= {u_i(x)}(R_i-|x-x_i|), $$

o{\`u}, $ x_i $ est le point tel que, $ u_i(x_i)=\max_{B_{R_i}(0)} u_i $.

\smallskip

Comme  $ R_i\to 0 $, on a $ R_i>R_i^2 $ et $ u_i(x_i) \to +\infty
$, on obtient:

$$ \max_{B_R(x_i)} s_i=s_i(a_i)\geq s_i(x_i)= { u_i(x_i) }R_i=\sqrt {
[u_i(x_i)]^2R_i^2 } \to +\infty . $$

\smallskip

En posant, $ l_i=R_i-|a_i-x_i| $, $ (l_i\to 0 ) $, on montre, comme
dans la d{\'e}monstration  du Th{\'e}or{\`e}me 1 que 
$$  L_i=\dfrac{l_i}{\sqrt {c_i} }u_i(a_i) \to +\infty .$$

\smallskip
 
Soit alors $ v_i $ la fonction d{\'e}finie par:

$$ v_i(y)=\dfrac{u_i\left (a_i+\dfrac{y}{u_i(a_i)} \right
  )}{u_i(a_i)} \,\,\, {\rm pour } \,\, |y|\leq \dfrac{l_i}{\sqrt {c_i} } u_i(a_i),  $$

on montre aussi, comme dans la d{\'e}monstration du Th{\'e}or{\`e}me 16, que pour $ c_i
\geq 4 $ (mais $ c_i \to +\infty $ ):

$$ v_i(y)\leq \dfrac{1}{ \left ( 1-\dfrac{1}{\sqrt {c_i}} \right )}
\leq 2 \, . \qquad (*) $$

Cette fonction v{\'e}rifie  $  \Delta v_i=W_i {v_i}^3 $ avec $ W_i(y)=V_i[a_i+\dfrac{y}{u_i(a_i)}] $.

\smallskip

Et la suite $  (v_i) $ converge uniform{\'e}ment vers la fonction  $
v(y)=\dfrac{1}{1+\dfrac{V(0)}{8}|y|^2} $ sur toute boule $
B_{\beta}(0),\, \beta>0 $, avec $ V(0)=\lim_{ i \to +\infty} V_i(a_i) $.

\bigskip

Les {\'e}tapes suivantes, sont identiques {\`a} celles de la d{\'e}monstration du cas de
la dimension 3, mais des modifications importantes sont {\`a} noter.

\smallskip

\underbar { Etape 2: Passage en polaires et propri{\'e}t{\'e}s de certaines
  fonctions }

\smallskip

Comme dans le cas de la dimension 3, on consid{\'e}re les op{\'e}rateurs
suivants:

$$ L_0=\partial_{tt}+2\partial_t-\Delta_{\sigma } \qquad et \qquad
L=\partial_{tt} -\Delta_{\sigma }, $$

$ \Delta_{\sigma } $ est l'op{\'e}rateur de Laplace-Baltrami sur $ {\mathbb
  S}_3 $

Soit, $ w_i $ la fonction suivante:

$$ w_i(t,\theta)=e^t u_i(a_i+e^t\theta), $$

elle v{\'e}rifie:

$$ - L w_i+w_i=V_i(a_i+e^t\theta){w_i}^3 .$$

Comme dans le cas de la dimension 3, on montre que  les fonctions $
w_i $ ont les propri{\'e}t{\'e}s suivantes:

\smallskip
 
a)  $ w_i(t+\log \eta_i             ,\theta
   )= e^{t} \times
   \dfrac{u_i(a_i+\dfrac{e^t\theta}{{u_i(a_i)}})}{u_i(a_i)}= e^{t}
   \times v_i(e^t\theta) $ converge vers la fonction sym{\'e}trique $   w= {\left
  (\dfrac{e^{t}}{1+\dfrac{V(0)}{8}
  e^{2t}} \right )} $ uniform{\'e}ment sur $ ]-\infty
   ,\log\alpha]\times {\mathbb S}_{3} $ pour tout  $ \alpha >0 $.

\smallskip

b)   Pour $ i\geq i_0,  w_i(\log {\eta}_i+\dfrac{1}{2}\log \dfrac{8}{V(0)},\theta)- w_i(\log
{\eta}_i+4+\dfrac{1}{2}\log \dfrac{8}{V(0)},\theta)>0 $, pour tout $ \theta $, avec $
\eta_i=\dfrac{1}{[u_i(a_i)]^{2/(n-2)}}=\dfrac{1}{u_i(a_i)} $ (on est en dimension 4).

\smallskip

c) si $ \mu >0 $ et si on pose $ \tilde w_i= w_i-\mu e^t $
alors:

\smallskip

$ \tilde w_i(\log \eta_i+\dfrac{1}{2}\log \dfrac{8}{V(0)},\theta)-\tilde
w_i(4+\log\eta_i+\dfrac{1}{2}\log \dfrac{8}{V(0)},\theta)>0 $ pour
tout $  \theta $.

\bigskip

\underbar { Etape 3: Utilisation de la technique " moving-plane "}

\smallskip

On pose:

\smallskip

$ \tilde w_i(t,\theta)=w_i(t,\theta)-\dfrac{ \min_{{}_{B_{2}(0)}} u_i}{2}
e^t $ ( le $ \mu $ du  point c pr{\'e}c{\'e}dent est \\ $ \mu=-\dfrac{ \min_{{}_{B_{2}(0)}} u_i}{2} $), $ \tilde V_i(t,\theta)=V_i(a_i+e^t\theta) .$

D'autre part:

$$ t^{\lambda}=2\lambda-t, \,\,\, \,\,\, \tilde
w_i^{\lambda}(t,\theta)=\tilde w_i(2\lambda-t,\theta) \,\,\, et \,\,\,\tilde
V_i^{\lambda}(t,\theta)=\tilde V_i(2\lambda-t,\theta) .$$

Ici, comme dans  du Lemme 1 pour la dimension 3, on cherche {\`a} savoir
si les fonctions qu'on utilise sont positives, le choix de $ \mu=\dfrac{
  \min_{{}_{B_{2}(0)}} u_i}{2} $ dans le c) de l'{\'e}tape pr{\'e}c{\'e}dente sera tr{\`e}s important. Nous avons ici: 

\bigskip

$ t\leq 0\Rightarrow e^t\leq 1            \Rightarrow
u_i(a_i+e^t\theta) \geq \min_{{}_{B_{1  }(a_i)}} u_i\geq
\min_{{}_{B_2(0)}} u_i $, car $ a_i \to 0 $.

\smallskip

D'o{\`u} pour $ t\leq 0               $ et pour tout $ \theta $ dans $
{\mathbb S}_3 $: 

\smallskip

$ \tilde w_i(t,\theta)=e^t u_i(a_i+e^t\theta)-\dfrac{min_{{}_{B_2(0)}} u_i}{2}
e^t \geq  \dfrac{min_{{}_{B_2(0)}} u_i}{2}
e^t > 0 $.

\smallskip

Dans le cas de la dimension 3, la
borne de droite des intervalles sur lesquels on applique le principe
du maximum varie, ici plus simplement  $ t_i\equiv t_0 = 0              $ est fixe.

\bigskip

Concernant le Lemme 2 ainsi que le point utile 4, ils sont les m{\^e}mes, puis on montre que:

\smallskip

\qquad $ \xi_i=sup \, \{ \lambda \leq \bar \lambda_i+2+\dfrac{1}{2}\log \dfrac{8}{V(0)},\,\, \tilde
w_i^{\lambda}-\tilde w_i<0, \, $ sur $ ]\lambda,t_0]\times {\mathbb S}_3 \}
$ existe.

\smallskip

Enfin, par continuit{\'e} des fonctions $ \tilde w_i $, on obtient:

$$ \forall \, (t,\theta) \in ]\xi_i,t_0]\times {\mathbb S}_3, \,\, \tilde
w_i ^{\xi_i}-\tilde w_i \leq 0 . $$

\underbar {Lemme :}

\smallskip

$$ \tilde w_i^{\xi_i}-\tilde w_i <0 \Rightarrow -\bar L(\tilde
w_i^{\xi_i}-\tilde w_i)<0 .$$

\underbar {D{\'e}monstration :}

$$ -\bar L(\tilde w_i^{\xi_i}-\tilde w_i)=\tilde V_i^{\xi_i}{(
  w_i^{\xi_i})}^3 -\tilde V_i{ w_i}^3 .$$

D'o{\`u}:

$$ -\bar L(\tilde w_i^{\xi_i}-\tilde w_i)=(\tilde V_i^{\xi_i}-
 \tilde V_i) {(w_i^{\xi_i})}^3+[
  {(w_i^{\xi_i})}^3 -{ w_i}^3]\tilde V_i  .$$

Pour tous $ t\in [\xi_i,t_0] $ et $ \theta \in {\mathbb
  S}_3 $ :

$$ \tilde V_i^{\xi_i}(t,\theta)-\tilde
V_i(t,\theta)=V_i(a_i+e^{2\xi_i-t}\theta )-V_i(a_i+e^t\theta) \leq A_i
(e^t-e^{2\xi_i-t}) . $$ 

D'autre part, si $ \tilde
  w_i^{\xi_i} -{ \tilde w_i} <0 $, alors par d{\'e}finition de $ \tilde w_i $,
  on obtient:

$$  w_i^{\xi_i} - w_i \leq 
\dfrac{ \min_{{}_{B_2(0)}} u_i }{2}(e^{2\xi_i-t}-e^t) <0  .$$

Et en utilisant le fait que $ 0 < w_i^{\xi_i}< w_i $, on obtient:
$$ (w_i^{\xi_i})^3 -{ w_i}^3=(w_i^{\xi_i}
-w_i)[(w_i^{\xi_i})^2+w_i^{\xi_i} w_i+(w_i)^2] \leq 3 (w_i^{\xi_i}
-w_i) \times (w_i^{\xi_i})^2   .$$

Ces deux in{\'e}galit{\'e}s entrainent, pour tous $ t\in [\xi_i,t_0] $ et $ \theta \in {\mathbb
  S}_3 $ :

$$  \qquad (w_i^{\xi_i})^3 -{ w_i}^3 \leq 3 \dfrac{ \min_{{}_{B_2(0)}} u_i
  }{2}\, (w_i^{\xi_i})^2 (e^{2\xi_i-t}-e^t) . $$

En cons{\'e}quence, on obtient:

$$ -\bar L(\tilde w_i^{\xi_i}-\tilde w_i)\leq (w_i^{\xi_i})^2 \,
(\dfrac{3  min_{{}_{B_2(0)}} u_i
  }{2} \tilde V_i-A_i{w_i}^{\xi_i})\, (e^{2\xi_i-t }-e^t). \qquad (**)$$

Par d{\'e}finition de $ w_i $ et d'apr{\`e}s $ (*)$  de l'{\'e}tape 1, rappelons
que pour tout $  t\leq \log
(l_i)-\log 2 +\log \eta_i $, on a

$$ w_i(t,\theta)=e^t\times \dfrac{u_i\left ( a_i+\dfrac{e^t\theta}{u_i(a_i)}
  \right )}{u_i(a_i)} \leq 2 e^t . $$

Comme,

$$  {w_i}^{\xi_i}(t,\theta)= {w_i}(2\xi_i-t,\theta)=
{w_i}[(\xi_i-t)+(\xi_i-\log
\eta_i)+\log\eta_i,\theta], $$ 

nous trouvons que 

$$ {w_i}^{\xi_i}(t,\theta)= e^{(\xi_i-t)+(\xi_i-\log \eta_i)}\times
\dfrac{u_i\left ( a_i+\dfrac{ e^{(\xi_i-t)+(\xi_i-\log \eta_i)}}{u_i(a_i)}
   \theta \right )}{u_i(a_i)} \leq 2 e^2 \sqrt {\dfrac{8}{V(0)}} \leq 2e^2\sqrt {\dfrac{8}{a}},$$

car, $ \xi_i-\log \eta_i\leq 2+\dfrac{1}{2}\log \dfrac{8}{V(0)} $ et $ \xi_i\leq t\leq t_0 $. La
constante $ 2e^2\sqrt {\dfrac{8}{a}} $, peut {\^e}tre largement am{\'e}lior{\'e}e.

\smallskip

Revenons {\`a} $ (**) $ et regardons le signe de :

$$ \dfrac{3 \, min_{{}_{B_2(0)}} u_i
  }{2} \tilde V_i-A_i {w_i}^{\xi_i} \geq  \dfrac{3 \, a \,  min_{{}_{B_2(0)}} u_i
  }{2}-2e^2\sqrt {\dfrac{8}{a}} A_i =\dfrac{3 a \, A_i}{2}\times \left [\dfrac{min_{{}_{B_2(0)}} u_i
    }{A_i}-\dfrac{8e^2 \sqrt {2}}{3\,a \sqrt{a}}\right ]. $$

D'apr{\`e}s notre hypoth{\`e}se de d{\'e}part, $ \liminf \dfrac{ \min_{{}_{B_2(0)}} u_i
    }{A_i} \geq  \dfrac{8e^2 \sqrt {2}}{3\,a \sqrt{a}}$, on en conclut  que $ (**) $ est n{\'e}gative, et le
    Lemme est d{\'e}montr{\'e}. 

\smallskip

La fin de la d{\'e}monstration est semblable {\`a} celle du th\'eor\`eme 2. On a, apr{\`e}s avoir appliquer le principe du maximum :

$$ \min_{\theta \in {\mathbb S}^3} \tilde w_i(t_0,\theta) \leq \max_{\theta
  \in {\mathbb S}^3}  \tilde w_i(2\xi_i-t_0,\theta) . $$

Comme $ t_0=0                      $ et $ a_i\to 0 $ et $ B_1(a_i)
\subset B_2(0) $, on obtient:

$$ \tilde
w_i(t_0,\theta)=e^{t_0}[u_i(a_i+e^{t_0}\theta)-\dfrac{\min_{B_2(0)}
  \, u_i}{2} ] \geq \dfrac{e^{t_0}}{2} \min_{B_{2}(0)} u_i .$$

D'autre part, la convergence uniforme des $ w_i $ entraine:

$$ \tilde w_i(2\xi_i-t_0,\theta)=w_i(2\xi_i-t_0,\theta)-\dfrac{
  min_{{}_{B_2(0)}} u_i }{2} e^{2\xi_i-t_0} \leq
    w_i(2\xi_i-t_0,\theta)\leq c\times e^{\log \eta_i}   $$

Finalement,

$$  [u_i(a_i)] \times \inf_{{}_{B_{2}(0)}} u_i \leq \, c .$$
Et ceci, contredit notre hypoth{\`e}se de l'{\'e}tape 1 ( la constante $ c $
d{\'e}pend de $ t_0 $, elle est ind{\'e}pendante de $ i $).

\bigskip

{\bf Th\'eor\`eme 4}{\it (Bahoura)}. Si $ (u_i) $ et $ (V_i) $, sont deux suites de
  fonctions relatives {\`a} l'{\'e}quation $ (E_1) $, sur un ouvert $ \Omega $
  de $ {\mathbb R}^4 $, alors on a:

\smallskip

Si la constante de lipschitz $ A_i $ relative {\`a} $ V_i $ tend vers $
0 $ et si $ \min_{\Omega } u_i\geq m>0 $ pour tout $ i $, alors:

\smallskip

Pour tout compact $ K $ de $ \Omega $, il existe une constante $ c>0
$, ne d{\'e}pendant que de $ a, b, (A_i)_{i\in {\mathbb N}}, m, K $, telle
que:

$$ \sup_K u_i \leq c .$$ 

\underbar {\bf Preuve.}

\bigskip

Comme les constantes de Lipschitz $ A_i $ relative {\`a} $ V_i $ tendent
vers 0, on obtient:

$$ \dfrac{ \min_K u_i }{A_i}\geq \dfrac{m}{A_i}\to +\infty
>> \dfrac{8e^2 \sqrt {2}}{3\,a \sqrt{a}} , $$

avec $ K $  compact de $ \Omega $ et $ m
$ un minorant uniforme de la suite $ u_i $.

\smallskip
 
En appliquant le Th{\'e}or{\`e}me 3, on obtient:

$$ \sup_K u_i \leq \dfrac{c(a,b,(A_i)_{i\in \mathbb N},K,\Omega)}{m}. $$
 
\underbar{\it Le cas radial: platitude d'ordre  $ {\tiny (n-2)/2} $}

\bigskip

Sur la boule unit\'e de $ {\mathbb R}^n $ si on consid\'ere des conditions suppl\'ementaires sur $ (u_i) $ et $ (V_i) $, \`a savoir:

\bigskip

$ u_i $ et $ V_i $ sont radiales et

\smallskip

$$ |V_i(r)-V_i(r')|\leq A|r^{[(n-2)/2]+\epsilon}-r'^{[(n-2)/2]+\epsilon}|, \,\,\,
   \forall \, 0\leq r,r'\leq 1, \, \epsilon>0. $$

On obtient le:

\bigskip

{\bf Th\'eor\`eme 5}. Sous les conditions pr\'ec\'edentes, on a:

$$ [u_i(0)]^{\epsilon/  [(n-2)+\epsilon ]} \times u_i(1) \leq c $$

o\`u $ c>0 $ est une constante qui ne d\'epend que de $ a, b, A, \epsilon $.

\bigskip

On suppose que, $ \epsilon_0 = 0 $, $  u_i(1) \geq m > 0 $ et $ A = A_i \to 0 $ pour chaque indice $ i $. On a:

\bigskip

{\bf Th\'eor\`eme 6}{\it (Bahoura)}. Il existe une constante positive $ c = c[a, b, (A_i), m, n] $ telle que:

$$ u_i(0) \leq c. $$

\bigskip

\underbar{ \it Le cas d'une pertirbation nonlin\'eaire}

\bigskip

{\bf Probl\`eme 2}. Sur un ouvert $ \Omega $ de $ {\mathbb R}^n $, on consid\`ere l'\'equation:

$$ \Delta u =Vu^{N-1}+W{u}^{\alpha} \,\, {\rm et } \,\,  u>0, \qquad (E_2) $$

avec $ \dfrac{n}{n-2} \leq \alpha < N-1=\dfrac{n+2}{n-2} $. 

$ V $ et $ W $ v\'erifient, pour des r\'eels positifs donn\'es $ a,b,c,d,A,B $,

$$ 0< a \leq V(x) \leq b \,\,{\rm et }\,\, 0 <c \leq W(x) \leq d $$

$$ ||\nabla V||_{L^{\infty}} \leq A \,\, {\rm et} \,\, ||\nabla
W||_{L^{\infty}} \leq B. $$

On se pose la question de savoir si pour chaque compact $ K $ de $
\Omega $, il existe une constante $ c
$, ne d\'ependant que de $\alpha, a, b, c, d, A, B, K, \Omega $ telle qu'on ait, pour
toute solution $ u $ de $ (E_2) $:

\bigskip

$$ \sup_K u \times \inf_{\Omega} u \leq c. $$

{\bf Th\'eor\`eme 7}. On consid{\`e}re trois suites de fonctions $ (u_i) $, $ (V_i) $ et $
(W_i) $ solutions de $ (E_2) $, alors on a:

\bigskip

Pour tout compact $ K $ de $ \Omega $, il existe une constante $ c'>0
$, ne d\'ependant que de $ \alpha, a, b, c, d, A, B, K, \Omega $, telle qu'on
ait:

$$ \sup_K u_i \times \inf_{\Omega} u_i \leq c'. $$

\underbar {\bf Preuve.}

\bigskip

Soient $ \{ u_i\} $, $ \{V_i\} $ et $ \{W_i\} $ trois suites de
fonctions telles que:

$$ \Delta u_i=V_i{u_i}^{N-1}+W_i{u_i}^{\alpha} \,\,\, {\rm dans }\,\,\,
\Omega,  $$

avec $ u_i>0 $ et $ \alpha \in ]\dfrac{n}{n-2},\dfrac{n+2}{n-2}[ $, $
V_i $ et $ W_i $ v{\'e}rifiant les hypoth{\`e}ses du Probl{\`e}me 2.

\smallskip

Le sch{\'e}ma de la d{\'e}monstration est le m{\^e}me que celui du
Th{\'e}or{\`e}me 1. On commence par d{\'e}montrer une estimation locale en
utilisant les techniques blow-up et "moving-plane ".

\smallskip

On suppose $ \Omega=B_2(0) $ et on  cherche {\`a} d{\'e}montrer qu'il existe
deux constantes positives $ c $ et $ R< 2 $ telles que pour tout entier $
i $, on ait :

$$ \sup_{B_R(0)} u_i \times \inf_{B_2(0)} u_i \leq \dfrac{c}{R^{n-2}}.
$$

On raisonne par l'absurde en s'inspirant de la d{\'e}monstration du
th{\'e}or{\`e}me 16, on exhibe une suite de points $ (a_i) $ tendant vers
0, deux suites de r{\'e}els positifs $ (R_i) $, $ (l_i )$ tendant aussi
vers 0 et enfin une suite de fonctions  $ (v_i) $ born{\'e}es qui
convergent uniform{\'e}ment vers une certaine fonction positive $ v $.

\smallskip

Plus pr{\'e}cis{\'e}ment, on a:

$$ u_i(a_i) \times \inf_{B_2(0)} u_i \to +\infty, \qquad (*) $$

$$ v_i(y)=\dfrac{u_i[a_i+y[u_i(a_i)]^{-2/(n-2)}]}{u_i(a_i)},
\,\,\, {\rm pour } \,\, |y| \leq \dfrac{l_i}{2}[u_i(a_i)]^{2/(n-2)}=L_i, $$  

avec $ u_i(a_i)
\to +\infty $ et $L_i \to +\infty $.

\smallskip

Chaque  fonction $ v_i $ v{\'e}rifie, pour tout entier $ i $ et tout $ y
$, tel que  $ |y|\leq L_i $,

$$ 0< v_i(y) \leq \beta_i \leq 2^{(n-2)/2}\,\,\, {\rm avec } \,\,\,  \beta_i\to 1. $$

De plus,

$$ \Delta v_i=\bar V_i {v_i}^{N-1}+\dfrac{1}{[u_i(a_i)]^{N-1-\alpha}}
\bar W_i {v_i}^{\alpha}. $$

O{\`u} $ \bar V_i(y)=V_i[a_i+y[u_i(a_i)^{-2/(n-2)}] $ et
$\bar W_i(y)= W_i[a_i+y[u_i(a_i)]^{-2/(n-2)}] $.

\smallskip

Comme $ \alpha \in ]\dfrac{n}{n-2}, N-1[ $, on voit alors en utilisant les
th{\'e}or{\`e}mes de Ladyzhenskaya et d'Ascoli, que de la suite $ (v_i) $, on
peut extraire une sous-suite convergeant vers une fonction $ v\geq 0 $
v{\'e}rifiant:

$$ \Delta v=kv^{N-1} \,\, {\rm sur } \,\, {\mathbb R}^n 
\,\, v(0)=1 \,\, {\rm et} \,\, 0\leq v(y)\leq 1 \,\, \forall  \, y\in
{\mathbb R}^n, $$ 

avec $ 0 < a \leq k \leq b $ 

Par un changement d'{\'e}chelle, on peut toujours supposer que $ k=n(n-2)
$, et on sait que la fonction $ v $ d{\'e}finie pr{\'e}c{\'e}demment, ne
peut {\^e}tre que la suivante:

$$ v(y)=\left ( \dfrac{1}{1+|y|^2} \right )^{(n-2)/2} .$$

Maintenant, on peut consulter la d{\'e}monstration du Th{\'e}or{\`e}me 16 et utiliser la technique "moving-plane ".

\smallskip

On remarque que seul le Lemme 2 est {\`a} v{\'e}rifier. On commence par
pr{\'e}ciser quelques notations.

\smallskip

Posons pour $ t\in ]-\infty, \log2 ] $ et $ \theta \in {\mathbb S}_{n-1}
$ :

\smallskip

$ w_i(t,\theta)=e^{(n-2)t/2}u_i(a_i+e^t\theta), \,\,
\bar V_i(t,\theta)=V_i(a_i+e^t\theta) \,\,\, {\rm et }\,\,\,
\bar W_i(t,\theta)=W_i(a_i+e^t\theta). $

\smallskip

Par ailleurs, soit $ L $ l'op{\'e}rateur $
L=\partial_{tt}-\Delta_{\sigma}-\dfrac{(n-2)^2}{4} $, avec $
\Delta_{\sigma} $ op{\'e}rateur de Laplace-Baltrami sur $ {\mathbb
  S}_{n-1} $. 

\smallskip

La fonction $ w_i $  est solution de l'{\'e}quation suivante :

$$ -Lw_i=\bar V_i{w_i}^{N-1}+e^{[(n+2)-(n-2)\alpha] t/2} \times \bar W_i
{w_i}^{\alpha}. $$

On pose pour $ \lambda \leq 0 $ :

\smallskip

$ t^{\lambda}=2\lambda-t  $ $ w_i^{\lambda}(t,\theta)=w_i(t^{\lambda},\theta) $, $ \bar
V_i^{\lambda}(t,\theta)=\bar V_i(t^{\lambda},\theta) $ et $ \bar
W_i^{\lambda}(t,\theta)=\bar W_i(t^{\lambda},\theta) .$

\smallskip

Alors, pour pouvoir v{\'e}rifier si le Lemme 2 du Th{\'e}or{\`e}me 16 reste
valable, il suffit de noter que  la quantit{\'e} $ -L( w_i^{\lambda}-w_i) $
est n{\'e}gative lorsque $ w_i^{\lambda}-w_i $ l'est. En fait, pour chaque
indice $ i $, $ \lambda=
\xi_i \leq \log \eta_i+2 $, ($ \eta_i=[u_i(a_i)]^{(-2)/(n-2)}) $. 

\smallskip

Tout d'abord:

$$ w_i(2\xi_i-t,\theta)=w_i[(\xi_i-t+\xi_i-\log\eta_i-2)+(\log
\eta_i+2)] , $$

par d{\'e}finition de $ w_i $ et pour $ \xi_i \leq t $:

 $$ w_i(2\xi_i-t,\theta)=e^{[(n-2)(\xi_i-t+\xi_i-\log\eta_i-2)]/2}e^{n-2}v_i[\theta e^2e^{(\xi_i-t)+(\xi_i-\log\eta_i-2)}]
\leq 2^{(n-2)/2}e^{n-2}=\bar c. $$

On sait que

$$ -L( w_i^{\xi_i}-w_i)=[\bar V_i^{\xi_i
  }(w_i^{\xi_i})^{N-1}-\bar V_i
{w_i}^{N-1}]+[e^{\delta t^{\xi_i}}{\bar W_i}^{\xi_i
  }(w_i^{\xi_i})^{\alpha }-e^{\delta t}\bar W_i
{w_i}^{\alpha} ] ,  $$

avec $ \delta = \dfrac{(n+2)-(n-2)\alpha}{2} $.

Les deux termes du second membre, not{\'e}s  $ Z_1 $ et $ Z_2 $, peuvent s'{\'e}crire:

$$ Z_1=(\bar V_i^{\xi_i }-\bar V_i)(w_i^{\xi_i })^{N-1}+\bar
V_i[(w_i^{\xi_i })^{N-1}-{w_i}^{N-1}],$$

et

$$ Z_2=(\bar
W_i^{\xi_i }-\bar W_i)(w_i^{\xi_i })^{\alpha }e^{\delta
  t^{\xi_i }}+e^{\delta
  t^{\xi_i }}\bar W_i[(w_i^{\xi_i })^{\alpha }-{w_i}^{\alpha
  }]+\bar W_i {w_i}^{\alpha }(e^{\delta t^{\xi_i }}-e^{\delta t} ). $$

D'autre part, comme dans la d{\'e}monstration du Th{\'e}or{\`e}me 3:

$$ {w_i}^{\xi_i} \leq w_i \,\,\,{\rm et } \,\,\,
w_i^{\xi_i}(t,\theta)\leq \bar c \,\,\, {\rm pour \, tout } \,\,\,
(t,\theta)\in [\xi_i,\log 2] \times {\mathbb S}_{n-1}  , $$

o{\`u}  $ \bar c $ est une constante positive ind{\'e}pendante de $ i
$ de $ w_i^{\xi_i} $ pour $ \xi_i \leq \log \eta_i+2 $;

$$ |\bar V_i^{\xi_i }-\bar V_i|\leq A (e^t-e^{ t^{\xi_i }}) \,\,\,
{\rm et } \,\,\, |\bar W_i^{\xi_i }-\bar W_i|\leq B (e^t-e^{ t^{\xi_i
    }}), $$

D'o{\`u}

\smallskip

$ Z_1 \leq A\, ({w_i^{\xi_i}})^{N-1} \,  (e^t-e^{ t^{\xi_i }}) \,\,\,
{\rm et} \,\,\,  Z_2 \leq  B \, ({(w_i ^{\xi_i})}^{\alpha }\,  (e^t-e^{
  t^{\xi_i }})+ c\, {(w_i^{\xi_i})}^{\alpha} \times   (e^{\delta t^{\xi_i
    }}-e^{\delta t} ) 
  $.

\smallskip

Ainsi, 

\smallskip

$ -L(w_i^{\xi_i}-w_i) \leq (w_i^{\xi_i})^{\alpha}[ (A\,  {w_i^{\xi_i}}^{N-1-\alpha}+ B)  \,  (e^t-e^{ t^{\xi_i }})+ c\,\,  (e^{\delta t^{\xi_i
    }}-e^{\delta t} ) ].
$
\smallskip

Puisque $ w_i^{\xi_i} \leq \bar c $, on obtient:

\smallskip

$ -L(w_i^{\xi_i}-w_i)\leq  (w_i^{\xi_i})^{\alpha } [(A {\bar c}
^{N-1-\alpha}+ B)  \,  (e^t-e^{ t^{\xi_i }})+ c\,\,  (e^{\delta t^{\xi_i
    }}-e^{\delta t} ) ]. \,\,\,(1)$

\smallskip

D{\'e}terminons le signe de 
 $ \bar Z=[( A{\bar c}
^{N-1-\alpha}+ B)  \, (e^t-e^{ t^{\xi_i }})+ c\,\,  (e^{\delta t^{\xi_i
    }}-e^{\delta t} ) ]. $

\smallskip

Comme $ \alpha \in ] \dfrac{n}{n-2}, \dfrac{n+2}{n-2}[ $, $
\delta=\dfrac{n+2-(n-2)\alpha}{2} \in ]0,1[ $. 

\smallskip

On d{\'e}duit que pour $ t \leq  t_0 <0 $:

$$ e^t \leq   e^{(1-\delta
  )t_0} e ^{\delta t} ,\,\,\, {\rm pour \,\, tout } \,\,\, t\leq t_0 .$$

Comme $ t^{\xi_i}\leq t $ $(\xi_i \leq t )$, en int{\'e}grant les deux membres, on
obtient:

$$ e^t-e^{ t^{\xi_i }} \leq \dfrac{ e^{(1-\delta
  )t_0}}{\delta }(e^{\delta t}-e^{\delta t^{\xi_i
    }}), \,\,\, {\rm pour\,\, tout } \,\,\, t\leq t_0 ,$$

ce qui s'{\'e}crit

$$ (e^{\delta t^{\xi_i
    }}-e^{\delta t} ) \leq \dfrac{ \delta} { e^{(1-\delta
  )t_0} }\,  ( e^{ t^{\xi_i }}-e^t). $$

L'in{\'e}galit{\'e} $ (1) $ devient alors :

$$ -L(w_i^{\xi_i}-w_i) \leq (w_i^{\xi_i})^{\alpha}[-\dfrac{\delta \, c }{ e^{(1-\delta
  )t_0} }+ A \,{\bar c}^{N-1-\alpha}+B]( e^t-e^{ t^{\xi_i }}). $$

Pour $ t_0 < 0 $, assez petit, la quantit{\'e} $ \dfrac{\delta \, c} {
  e^{(1-\delta
  )t_0}}  - A \,{\bar c}^{N-1-\alpha}-B  $ devient
  positive et le r{\'e}sultat cherch{\'e} est obtenu dans l'intervalle $
  [\xi_i,t_0] $.

\smallskip

Le fait de prendre l'intervalle  $  [\xi_i,t_0] $ au lieu de  $
[\xi_i, \log 2] $, n'est pas g{\^e}nant, au contraire, plus l'intervalle
est petit plus l'infimum est grand. La suite de la d{\'e}monstration est
identique {\'a} celle de la fin du Th{\'e}or{\`e}me 16.

\smallskip

On pourrait croire que $ t_0 $ d{\'e}pend de $ \xi_i $ ou de $
w_i^{\xi_i} $, mais  $ t_0 $ d{\'e}pend seulement de $
\bar c $, une constante qui ne d{\'e}pend que de $ n $, $ a $ et $ b $.

\smallskip

On calcule $ t_0 $ puis on introduit $ \xi_i \leq \log \eta_i+2 $ comme dans les autres
th{\'e}or{\`e}mes, et on v{\'e}rifie  l'in{\'e}galit{\'e} $ L(w_i^{\xi_i}-w_i )\leq 0 $, d{\`e}s
que $ w_i^{\xi_i}-w_i \leq 0 $ sur $ [\xi_i,t_0].$

\smallskip

Ayant d{\'e}termin{\'e} $ t_0 <0 $ tel que $\dfrac{\delta \, c} { e^{(1-\delta
  )t_0} }- A \,{\bar c}^{N-1-\alpha}-B $ soit positive, on pose:

\smallskip

$ \xi_i =\sup \{ \mu_i \leq \log \eta_i+2, w_i
  ^{\mu_i}(t,\theta)-w_i(t,\theta)  \leq 0, \forall \, (t,\theta) \in
  [\mu_i,t_0]\times {\mathbb S}_{n-1} \}$.

\smallskip

Par d{\'e}finition de $ \xi_i $, $ w_i^{\xi_i}-w_i \leq 0 $. Ensuite, on
 v{\'e}rifie que $ -L(w_i ^{\xi_i}-w_i) \leq 0 $.

\smallskip

Comme dans le Th{\'e}or{\`e}me 16, le principe du maximum, entra{\^\i}ne:

\smallskip

$ \qquad \qquad \min_{\theta \in {\mathbb S}_{n-1}}w_i(t_0,\theta) \leq \max_{\theta
  \in {\mathbb S}_{n-1}} w_i(2\xi_i-t_0) \,.$

\smallskip

Or,

$ \qquad w_i(t_0,\theta)=e^{t_0} u_i(a_i+e^{t_0}\theta)\geq e^{t_0} \min u_i
\,\,{\rm et} \,\, w_i(2\xi_i-t_0)\leq \dfrac{ c_0 }{u_i(a_i)} $,

\smallskip

donc:

$$ u_i(a_i) \times \min u_i \leq c .$$

Ce qui contredit notre hypoth{\`e}se  $ (*) $.

\renewcommand{\bibname}{R\'ef\'erences suppl\'ementaires}

\end{document}